\documentclass[onecolumn]{svjour3}

\smartqed
\pdfoutput=1 
\usepackage{graphicx}
\usepackage{float}
\usepackage{adjustbox}
\usepackage{mathptmx}
\usepackage{epstopdf}
\usepackage{amssymb}
\usepackage{amsmath}
\usepackage{tabularx}
\usepackage{multirow}
\usepackage{booktabs}
\usepackage{microtype} 
\usepackage[sort,numbers]{natbib} 
\usepackage[parfill]{parskip} 
\usepackage{soul}
\usepackage{color}
\usepackage{cancel}
\usepackage{siunitx} 
\usepackage[hidelinks]{hyperref}

\usepackage{subcaption}
\usepackage{tikz}
\usetikzlibrary{patterns}
\usepackage{pgfplots} 
\usepackage{mathrsfs}

\pgfplotsset{every axis/.append style={
                    label style={font=\Large},
                    tick label style={font=\Large}  
                    }}

\journalname{Computational Particle Mechanics}

\begin{document}

\title{Framework for uncertainty quantification of wave-structure interaction in a flume}
\author{Xiaoyuan Luo \and Vijay Nandurdikar \and \\ Sangri-Yi \and  Alistair Revell \and Georgios Fourtakas \and \\ Ajay B. Harish}
\authorrunning{X. Luo et. al.}
\institute{Xiaoyuan Luo, Vijay Nandurdikar, Georgios Fourtakas, Alistair Revell \at
    School of Engineering\\
    University of Manchester (UK) \\
    \email{Xiaoyuan.Luo@postgrad.manchester.ac.uk}
    \and
    Sangri-Yi \at 
    NHERI SimCenter, University of California\\
    Berkeley, CA (USA) 
    \and
    Ajay B. Harish \at
	School of Engineering\\
    University of Manchester (UK) \\
    \email{Ajay.Harish@manchester.ac.uk}
}
\date{Received: date / Accepted: date}

\maketitle

\begin{abstract}
In this paper we propose a numerical procedure for the quantification of uncertainties in wave-structure interaction. We utilise the smoothed particle hydrodynamics (SPH) scheme for modelling the wave mechanics, coupled one-way with a finite element method (FEM) for the structural response. Physical wave flumes are extensively used in the study of hydrodynamics especially in wave-structure interaction (WSI) and prediction of forces to near-shore structures in disaster mitigation and offshore structures in the oil and gas, and more recently renewable energy sector. Over the years, numerical wave flumes have been developed extensively to enable the modelling of complex wave-structure interaction. However, most of these studies are deterministic and limited to using either simple flexible beams or rigid monolithic structures to model the structural part in the WSI. Additionally, uncertainties are commonly observed in both wave and structural parameters and need to be accounted for. This work presents a numerical framework to enable uncertainty quantification for wave-structure interaction problems in terms of the forces experienced by the structure. A one-way coupling between SPH with the FEM and uncertainty quantification (UQ) is proposed. We employ the so-called  Tokyo wave flume geometry, which has potential for future surrogate modelling in WSI. The developed model is validated using numerical and experimental results from the literature and is used to demonstrate the prediction of probabilistic responses of structures under breaking and non-breaking wave scenarios.

\end{abstract}

\section{Introduction}
Water-borne hazards, particularly tsunamis, storm surges, and inland and coastal flooding events, are known for their devastating effects. The 2004 Indian Ocean tsunami caused waves up to 30 \unit{meters} high, resulting in the deaths of over 250,000 people. Similar large-scale destruction was observed with damage to ports and infrastructure inflicted by the 2010 Chile and 2011 Tohoku tsunamis, oil rigs and other infrastructure during the 2005 Hurricane Katrina storm surge, and large-scale flooding events attributed to adverse weather effects in parts of continental Europe and India during the summer of 2023. Many of these events, once considered to occur only once a decade, are happening more frequently. This increase is expected to continue due to climate change and increased human activity in coastal areas. Therefore, understanding the impact of these extreme phenomena on the urban built environment remains paramount to improving structural resilience. Furthermore, there is a clear drive for clean energy in offshore renewable structures such as wind turbines and hybrid wind and current offshore platforms.

The value of experimental wave flume campaigns for coastal engineers in systematically studying wave-structure interaction (WSI) has been well demonstrated. In recent years, a wide range of coastal engineering problems have been investigated using wave flume experiments, including beach erosion and accretion, the design of breakwaters and seawalls, wave impact on offshore structures, and sediment transport by waves \cite{De2008analysis, Ribberink2001near, Cuomo2010breaking}. However, the construction and operation of wave flumes can be expensive and time-consuming, particularly when investigating different scenarios. Additionally, careful design and scaling of flume experiments are necessary to ensure their real-world relevance. To reduce the cost of experimental campaigns, especially during the design phase, numerical models can be employed to assist with the design of wave flume experiments. This approach can also serve to validate the accuracy of the numerical models themselves.

In an effort to replicate the experiments by computer simulations, digital wave flumes were developed using numerical methods. Common methods and tools that have been adopted include potential flow solvers (e.g., OceanWave3D \cite{OceanWave3D}), shallow water solvers (e.g., GeoClaw \cite{Leveq2011, BergerGeorgeLeVequeMandli11}, AdCirc \cite{Adcirc}), finite volume (FVM) solvers (e.g., OpenFOAM \cite{Chen2014-rz, Hu2016-ux, Jacobsen2012-gs}), and mesh-less solvers (e.g., smoothed particle hydrodynamics (SPH) \cite{ALTOMARE201434, ALTOMARE20151, ALTOMARE201737, Tagliafierro2023}). A comprehensive review of developments in WSI, encompassing both numerical methods and experiments, are provided in the recent work of \citet{Huang2022-lj}.

Herein, three numerical methods, namely SPH for wave modelling, the Finite Element Method (FEM) for structural response, and Latin Hypercube Sampling (LHS) for sampling structural parameters, are employed to facilitate forward uncertainty quantification (UQ). While SPH, FEM, LHS, and UQ are all well-established methods with extensive development in their respective fields, the novelty of this work lies in coupling these components to address uncertainties in the system structural response. Moreover, most of the works where SPH has been coupled with a structural solver have been limited to the adoption of monolithic rigid structures (cubes, cylinders etc.) or flexible beam structures. Additionally, within the context of the SPH framework, UQ has seen limited integration previously. More specifically, this work explores UQ in the structural response under wave loading. This work further serves as a precursor to potential developments of a surrogate model-driven digital wave flume that could enable the incorporation of more extensive probabilistic wave conditions in the future.

The paper is structured as follows: Section 2 outlines the basic theoretical formulations related to SPH and UQ, followed by validation with experimental and OpenFOAM results from the literature in Section 3. Section 4 provides a discussion on the force calculation and probabilistic structural response. Finally, the work concludes with a discussion of the conclusions and outlines focus areas for the future.

\section{Numerical methods}
This section provides a brief overview of the numerical methods, governing equations and UQ techniques employed to quantify uncertainties within the complex WSI system.

\subsection{Smoothed-particle hydrodynamics}
This work adopts the weakly compressible Smoothed Particle Hydrodynamics (WCSPH) approach to model wave generation and propagation within the numerical flume. This choice is driven by the versatility offered by the scheme, which leverages a Lagrangian description of motion and eliminates meshing. This mesh-less Lagrangian approach provides significant flexibility when simulating non-linear and fragmented flows, particularly in the presence of a free surface. This is especially advantageous for scenarios involving wave breaking, overtopping, and impact flows on structures. The WCSPH open-source DualSPHysics solver \cite{dom2021} with graphic processing unit (GPU) acceleration, is utilised in this work.

\subsubsection{Governing equations}
The continuous Lagrangian form of the mass and momentum conservation equations considered herein are,
\begin{equation}
\begin{split}
{\frac{D\rho}{D t}} &=-\rho\nabla\cdot\mathbf{u} \\
{\frac{D\mathbf{u}}{D t}} &=-{\frac{1}{\rho}}\nabla P+\nu \Delta \mathbf{u} +\mathbf{g}
\end{split}
\end{equation}
where $\rho$ is the density, $\mathbf{u}$ is the velocity vector, $\mathbf{g}$ are the body forces, $P$ is the pressure, and $\nu $ is the kinematic viscosity evolved in time $t$.

\subsubsection{SPH discretisation}
In SPH, the continuous approximation (or integral representation) of any generic, sufficiently smooth spatial function $\phi\left(\mathbf{r}\right)$ is expressed as the convolution of the function with a smoothing kernel,
\begin{equation}
\langle \phi({\bf r}) \rangle =\int_{\Omega}\phi({\bf r}^{\prime})W({\bf r}-{\bf r}^{\prime},h) \ d{\bf r}^{\prime}
\label{equation-integral}
\end{equation}
where $\mathbf{r}$ is the position vector, $\mathbf{r}^{\prime}$ is the position within the volume $\Omega$ bounded by $h$ the characteristic smoothing length and $W({\bf r}-{\bf r}^{\prime},h)$ is the smoothing kernel function. SPH uses the smoothing kernel function or simply kernel to calculate the weighted and averaged contributions of the particles within a cut-off distance of $\kappa h$ where $\kappa$ is a kernel specific property and defines the extent of the kernel support domain \cite{DeAl2012}. This volume averaging approximation is represented by the angled brackets $\langle \cdot \rangle$.

In SPH, the domain is discretised into computational nodes, or particles, that are assigned fluid properties such as mass, density, velocity, etc. The integral representation of the field function of \autoref{equation-integral} is discretised by summing the weighted contributions of the nearest neighbour particles within the $\kappa h$ support radius. This implies that for any arbitrary particle $a$, the approximation reads,
\begin{equation}
\langle \phi({\bf r}_{a}) \rangle=\sum_{b=1}^{N_{p}}\frac{m_{b}}{\rho_{b}}\phi({\bf r}_{b})W_{a b}
\end{equation}
where $a$ and $b$ are the interpolating and neighbouring particles, $N_p$ is the number of particles within the domain with $b = 1, \cdots, N_{p}$ neighbours, $m_{b}$ is the mass and $\rho_{b}$ is the density of the neighbouring particle $b$ that defines the volume of each particle as $V_{b}=\frac{m_{b}}{\rho_{b}}$. Herein, ${W}_{a b}\,\equiv\,{W} \left({\bf r}_{a}-{\bf r}_{b},h \right)$ is the kernel function. Note that the brackets $<.>$ that denote approximation that will be dropped in the rest of the paper for brevity. 

This work in particular, uses the Wendland $C^2$ kernel function implemented in DualSPHysics \cite{Wendland1995} that reads,
\begin{equation}
W(r,h)=\alpha_{d}\Biggl\{\left(1-\frac{q}{2}\right)^{4} \left(1+2q\right)\;\;0\leq\,q\,\leq2\,
\end{equation}
where 
\begin{equation}
q={\frac{|{\mathbf{r}}-{\mathbf{r}}^{\prime}|}{h}} = {\frac{|r_{ab}|}{h}} 
\end{equation}
and $\alpha_{d}$ is a normalisation constant that depends on the spatial dimensions of the problem with $\kappa=2$. 

Following the above discretisation technique, the work uses the weakly compressible SPH formulation implemented in DualSPHysics \cite{dom2021} to obtain the following form of the governing equations,
\begin{equation}
\begin{split}
\frac{D\rho_{a}}{D t} &=-\rho_{a}\sum_{b=1}^{N_{p}}\frac{m_{b}}{\rho_{b}}\mathbf{u}_{a b}\cdot\nabla_{a}W_{a b} + \mathfrak{D}\\
{\frac{D{\mathbf{u}}_{a}}{D t}} &=-\sum_{b=1}^{N_{p}}m_{b}\left({\frac{P_{b}+P_{a}}{\rho_{a}\rho_{b}}}+\Pi_{a b}\right)\cdot\nabla_{a}W_{a b} + \mathbf{g}
\end{split}
\label{wc1}
\end{equation}
where $\mathfrak{D}$ is a density diffusion term added to dissipate checkerboarding, defined as a second order derivative of the hydrostatic density \cite{FOURTAKAS2019346}, and $\Pi_{a b}$ is an artificial viscosity \cite{Mo1992} term given by, 
\begin{equation}
\Pi_{a b}= 
  \begin{cases} 
   -\frac{ \alpha \overline{c_{ab}}\mu_{ab} }{\overline{\rho_{ab}} } & \text{if } x \geq 0 \\
   0      & \text{if } x < 0
  \end{cases}
\end{equation}
and $\alpha$ is a  parameter which can be tuned to provide sufficient numerical dissipation. In this work, following \citet{ALTOMARE20151}, the value of $\alpha$ is set to 0.01 to provide sufficient dissipation for numerical stability in wave and wave loading flume applications. 

A barotropic equation of state (EOS) \cite{monaghan1999} closes the system of \autoref{wc1} by relating pressure to density as
\begin{equation}
P = {\frac{c_{0}^{2}\rho_{0}}{\gamma}}\left[\left({\frac{\rho}{\rho_{0}}}\right)^{\gamma}-1\right]
\label{eqstate}
\end{equation}
where the polytropic index $\gamma = 7$, the initial fluid density is $\rho_{0} = 1000$ \unit{\kilogram\per\cubic\meter} and the initial speed of sound is defined as $c_{0} = c \left(\rho_{0} \right)$ within the incompressible limit of $Ma=0.1$ by imposing $c_{0} =10\sqrt{gH}$ where H is initial water height at rest.

The symplectic position Verlet scheme \cite{leimkuhler2015molecular} is employed for time integration in this work. Due to the Lagrangian nature of SPH a second order accurate time integration scheme is required, thus the velocity step is updated at $n+1/2$ using the velocity Verlet half step, as given by:
\begin{equation}
    \mathbf{u}^{n+1/2} = \mathbf{u}^{n} + \frac{1}{2} \Delta t \ \mathbf{a}^{n}
\end{equation}
The position is updated according to symplectic position Verlet scheme,
\begin{equation}
    \mathbf{r}^{n+1} = \mathbf{r}^{n} + \Delta t \ \mathbf{u}^{n} + \frac{1}{2} \Delta t^2 \ \mathbf{a}
\end{equation}
with $\mathbf{a}=\left( \frac{D\mathbf{u}}{D t} \right)^{n}$. Furthermore, the density evolution is calculated using the half time steps of the symplectic position Verlet scheme \cite{parshikov2000improvements},
\begin{equation}
    \rho^{n+1/2} = \rho^{n} + \frac{1}{2} \Delta t  \ \mathscr{D}
\end{equation}
where $\mathscr{D}=\left( \frac{D\rho}{D t} \right)^{n}$ represents the rate of change of density at step $n$.

In this paper, the Dynamic Boundary Condition (DBC) \cite{cmc2007} is utilised to impose solid wall boundary. Boundary particles are subject to the continuity equation to evolve the density and thus impose a no-penetration condition. 

Boundary particles are allowed to move according to externally prescribed motion such as a prescribed motion or wave maker. The piston wave-maker generates waves by imposing the following specific motion function derived from the Rayleigh approximation,
\begin{equation}
    \eta(x_s, t) = H \, \text{sech}^2 \left[ k \left( c \left( t - \frac{T_f}{2} \right) + 2 \sqrt{\frac{H(H+h)}{3}} - x_s \right) \right]
\end{equation}
where, $2\sqrt{H(H+h)/3}$ represents half of the wavemaker stroke based on the assumption that wavemeker's motion starts from $x=0$; $k$ represents the outskirt coefficient, $c$ denotes the wave celerity which is a function of gravity ($g$), wave height ($H$), and water depth ($h$); $T_f$ indicates the time required to generate the solitary wave and $x_s$ the displacements of the wavemaker.

Over the years, several topical reviews \cite{liu_smoothed_2010, lyu_review_2022, hiermaier_review_2009, ye_smoothed_2019} and books \cite{Liu2003book, Liu2003book2,Liu2003book3} on the development and applications of SPH in coastal engineering have been written and reader is directed to the above for further information.

\subsection{Force exerted on the structure}
An accurate estimation of forces is crucial for assessing the structural response within a wave flume. This work compares force calculations obtained from experimental results, SPH simulations, ASCE standards, and analytical expressions. The ASCE standards and semi-empirical employ wave velocities and heights to determine the time history of wave loading on the structure.

\subsubsection{Force calculation in SPH}
The post-processing tool within the DualSPHysics is employed to extract the forces exerted on the structures. The force tool allows for the extraction of horizontal forces acting on the column, as
\begin{equation}
F_{x,SPH} = \sum_{i} {P_i \ dx^2} \frac{\left( \mathbf{x}_{ci} - \mathbf{x}_i \right)}{ | \mathbf{x}_{ci} | }
\end{equation}
where $\mathbf{x}_{ci} = \left(x_c,y_x\right)$ is the column center coordinates, $P_i$ is the pressure and $dx^2$ represents summation over all infinitesimal areas. More detailed discussions, including relevant formulations on force time-history extraction and comparison of SPH with experimental results \cite{Goda1974New, TAKAHASHI94} during impulsive provoked breaking waves can be found in \citet{ALTOMARE20151}. 

\subsubsection{American Society of Civil Engineers (ASCE)\label{asce-standard}}
\citet{asce} discuss a practical method for estimating the maximum wave breaking force (or shock pressure) on vertical walls of buildings in coastal areas. The ASCE relations, also considered in this work, for calculating the maximum combined dynamic and static wave pressures $P_{\mathrm{max}}$, are
\begin{equation}
P_{\mathrm{max}}= \left( C_{p} + 1.2 \right) \gamma_{w}d_{s}
\label{eq:asce-pressure}
\end{equation}
The net breaking wave force per unit length of the structure, also known as the wave impact force $\left(F_{\text{ASCE}}\right)$ is given to be
\begin{equation}
F_{\text{ASCE}} = \left( 1.1 C_p + 2.4 \right) \gamma_w d_s^2
\label{eq:asce}
\end{equation}
where, $C_{p}$ is the pressure coefficient which varies from 1.6 - 3.5 according to the category of structure, $\gamma_{w}$ is the unit weight of water and $d_{s}$ is the still water depth at base of structure.

\subsubsection{Semi-empirical calculation\label{analytical-calculation}}
The semi-empirical approach considered splits the overall wave-induced force into a hydrostatic and a dynamic component, i.e. $F_{\text{total}} = F_{\text{static}} + F_{\text{dynamic}}$. According to ASCE standards \cite{american2017minimum}, the static force calculation considers water density ($\rho = 1000\text{kg/m}^3$) and gravitational acceleration ($g = 9.81 \text{m/s}^2$), and computed as
\begin{equation}
F_{\text{static}} =\frac{1}{2} w \left(h_{b}-0.75 \right) \left[ \rho \times g \left(h_{b}-0.75 \right) \right]  
\label{eq:fstat}
\end{equation}
where, $h_{b}$ represents the gauge value obtained from the near structure wave gauge above the initial free surface; $w$ represents width of the structure. This static force component represents the force exerted on the vertical wall of the structure by a presumably standing water at equilibrium. In addition, the dynamic force caused by the water column is obtained by integrating the velocity data over the height range. Thus, the dynamic force $F_{\text{dynamic}}$ is dependent on the effective velocity $V_{\text{eff}}$ and is given to be
\begin{equation}
F_{\text{dynamic}} =\frac{1}{2} \rho  C_{d} A V_{\text{eff}}^{2}
\label{eq:fdyn}
\end{equation}
where, $C_{d}$ represents the drag coefficient, and $A$ indicates the surface area of the object normal to the flow.

\subsection{One-way coupling method}
\begin{figure}[!htb]
\centering
\includegraphics[width=0.95\textwidth]{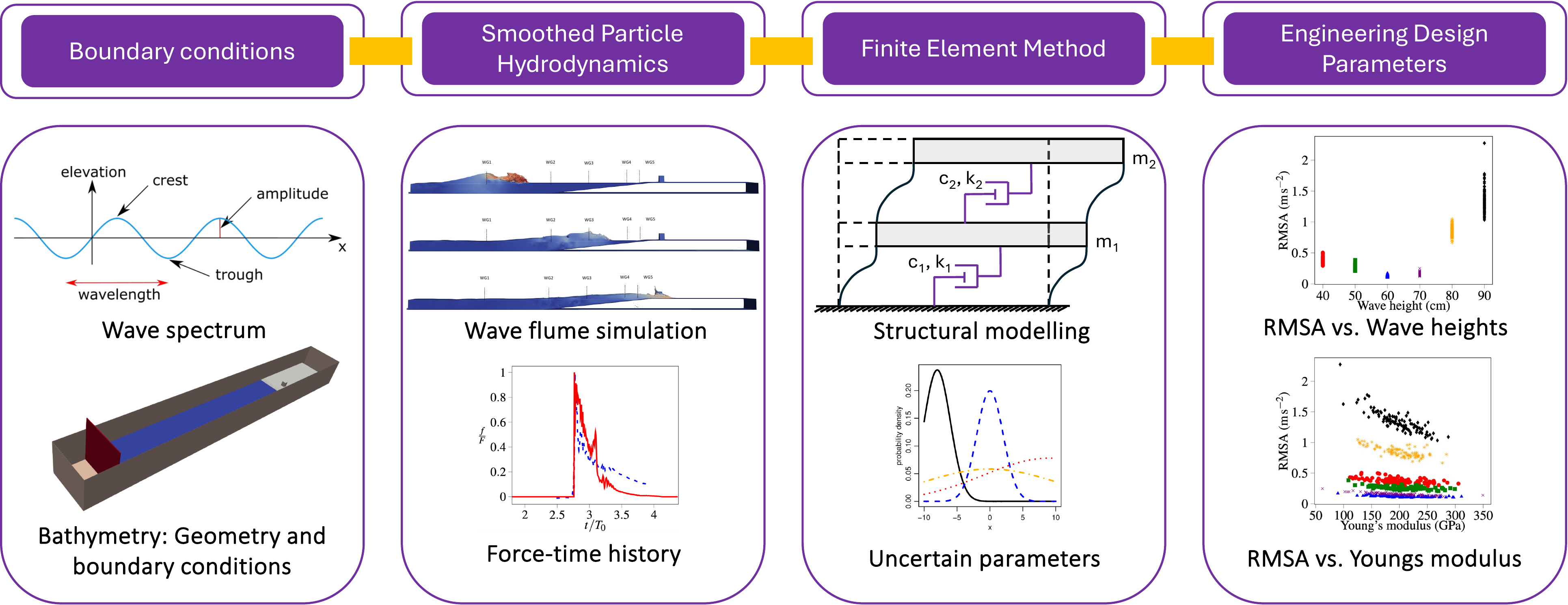}
\caption{Illustration of the one-way coupling}
\label{fig:OneWay}
\end{figure}

The one-way coupling presented in this work is visually illustrated in \autoref{fig:OneWay}. The work demonstrates the applicability for solitary waves but provides a general framework for other types of wave scenarios (regular, irregular, focused waves, etc.) with structures where the response of the structure does not significantly alter the hydrodynamics, i.e., with small displacements. 

Herein, due to the very small deformations of the structure, we consider an off-line coupling approach. The SPH solver propagates the wave and estimates forces that are obtained on the structure using the above mentioned implementations in DualSPHysics. Thus, the temporal dynamics onto the structure by the wave (forces in our work) are recorded on an external data file. The FEM structural solver obtains these temporal forces exerted on the structure by the data file and imposes additional constrains. Note that, the geometrical characteristics of the structure in both solvers are identical and the zero displacement constrain at the base of the structure is only imposed in the FEM solver. As the response of the structure is two orders of magnitude smaller that the SPH discretisation, there is no requirement to exchange the displacement with the SPH solver, resulting in a weak off line coupling. This approach has been applied in SPH successfully in \citet{BenLikeArticle}. The wave forces along with uncertainties in the material parameters are further used in the determination of structural response using FEM. The overall outputs of the numerical framework leads to engineering design parameters (EDP) such as root mean square acceleration (RMSA).

\subsection{Uncertainty quantification (UQ)}
UQ is combined with structural dynamic analyses to estimate the probabilistic response of structures. This work considers the widely adopted  Latin Hypercube Sampling (LHS) as it facilitates an efficient exploration of the parameter space by systematically partitioning the parameter ranges into equi-probable intervals and selecting samples such that one sample is placed in each interval of each parameter dimension \cite{LHSsampling}. The resulting LHS samples ensure a more uniform and representative distribution of samples across all dimensions, effectively preventing sample clustering.

Without loss of generality, the $n$ input parameters can be given to be $p_{i} \ \forall \ i=1,2,..,n$ and characterised by a specified range and divided into $q$ equi-probable intervals or bins. The interval width can be given to be
\begin{equation}
\mathrm{Interval~width}={\frac{\mathrm{Max}(p_{i})-\mathrm{Min}(p_{i})}{q}}
\label{eq:interval_width}
\end{equation}
where, $\mathrm{Max}(p_{i})$ and $\mathrm{Min}(p_{i})$ are maximum and minimum value possible for the parameter $p_{i}$ respectively, $q$ denotes the desired number of samples. The $q-$random permutations of the integers from 1 to $q$ (inclusive), for each parameter $p_{i}$, are generated using a random permutation function, namely the Fisher-Yates shuffle algorithm \cite{Fisher_Yates-AFP}.

The LHS matrix $L$ is further created, which contains the LHS samples in domain of $[0,1]^n$. The dimension of $L$ is $q \times n$ and each row corresponds to a sample point in the $n$-dimensional input space, with $L_{i}=(L_{i1},L_{i2},..,L_{i n})$, where $L_{i j}$ is the value of the $j$-th parameter in the $i$-th sample. The LHS sample matrix is constructed using the generated random permutations and the intervals. The elements of the LHS matrix $L$ are mapped to the actual parameter values based on the intervals created using
\begin{equation}
p_{i j}=\mathrm{Min}(p_{j})+(L_{i j}-1)\times\mathrm{Interval~width}
\label{eq:mappedvalue}
\end{equation}
where, $p_{i j}$ is the mapped value of the $j$-th parameter for the $i$-th sample. 

Given the randomness in the individual structural properties (here depicted as a vector of properties) $\mathbf{p}_{\rm{str}}$ and external forces $\mathbf{F}$, the probability density function (PDF) of the structural responses of interest like Root Mean Squared Acceleration (RMSA), $\mathbf{A}$, is computed as
\begin{equation}
\xi(\mathbf{A}) = \int \xi(\mathbf{A}| \mathbf{p}_{\rm{str}},\mathbf{F}) \ \xi(\mathbf{p}_{\rm{str}})\xi(\mathbf{F}) \ \rm{d}\mathbf{p}_{\rm{str}} \ \rm{d}\mathbf{F}
\label{eq:forward propagation}
\end{equation}
where $\xi(\cdot)$ denotes the PDF and $\xi\left(\mathbf{Y} | \mathbf{X}\right)$ denote the conditional PDF of any $\mathbf{Y}$ given $\mathbf{X}$. Note here that $\mathbf{p}$ is used to denote the set of all properties where uncertainties can exist. Some examples include Youngs' modulus, Poisson ratio etc. Without additional randomness considered, the first term in the integral becomes a Dirac-delta function, $\xi(\mathbf{A}|\mathbf{p}_{\rm{str}},\mathbf{F}) = \delta(h(\mathbf{p}_{\rm{str}},\mathbf{F})=\mathbf{A})$, where $h(\mathbf{p}_{\rm{str}},\mathbf{F})$ represents the system equation, i.e. a combination of wave flume and structural dynamic simulation model. The PDF is approximated by LHS and kernel density estimation (KDE) as
\begin{equation}
\hat{\xi}(\mathbf{A}) = \frac{1}{q} \sum^m_{i=1} K_h(\mathbf{A}-A_i) = \frac{1}{q} \sum^m_{i=1} K_h(\mathbf{A}-h(p_{\rm{str},i},F_{i}))
\label{eq:KDE}
\end{equation}
where $[{p_{\rm{str},i},F_{i}}]$ represents $i$-th sample obtained from LHS and $A_i$ is corresponding structural analysis outcome. $K_h(\cdot)$ is a kernel function and this work uses, the commonly used, Gaussian basis function \cite{Deierlein2020cloud}.

\section{Validation and convergence study}

The SPH results from this work are validated by comparison with experiments conducted at the hybrid tsunami open flume (HyTOFU) in Ujigawa laboratory at the Kyoto University (Japan) and simulations performed using the OpenFOAM volume of fluids (VoF) solver \cite{Trybook2011, DeAnTr2012, RoBrJa2016}. The experimental and OpenFOAM results used for validation are reported in \citet{Joetal2021}.

\subsection{Experimental setup}
\begin{figure}[!htb]
\centering
\includegraphics[width=0.95\textwidth]{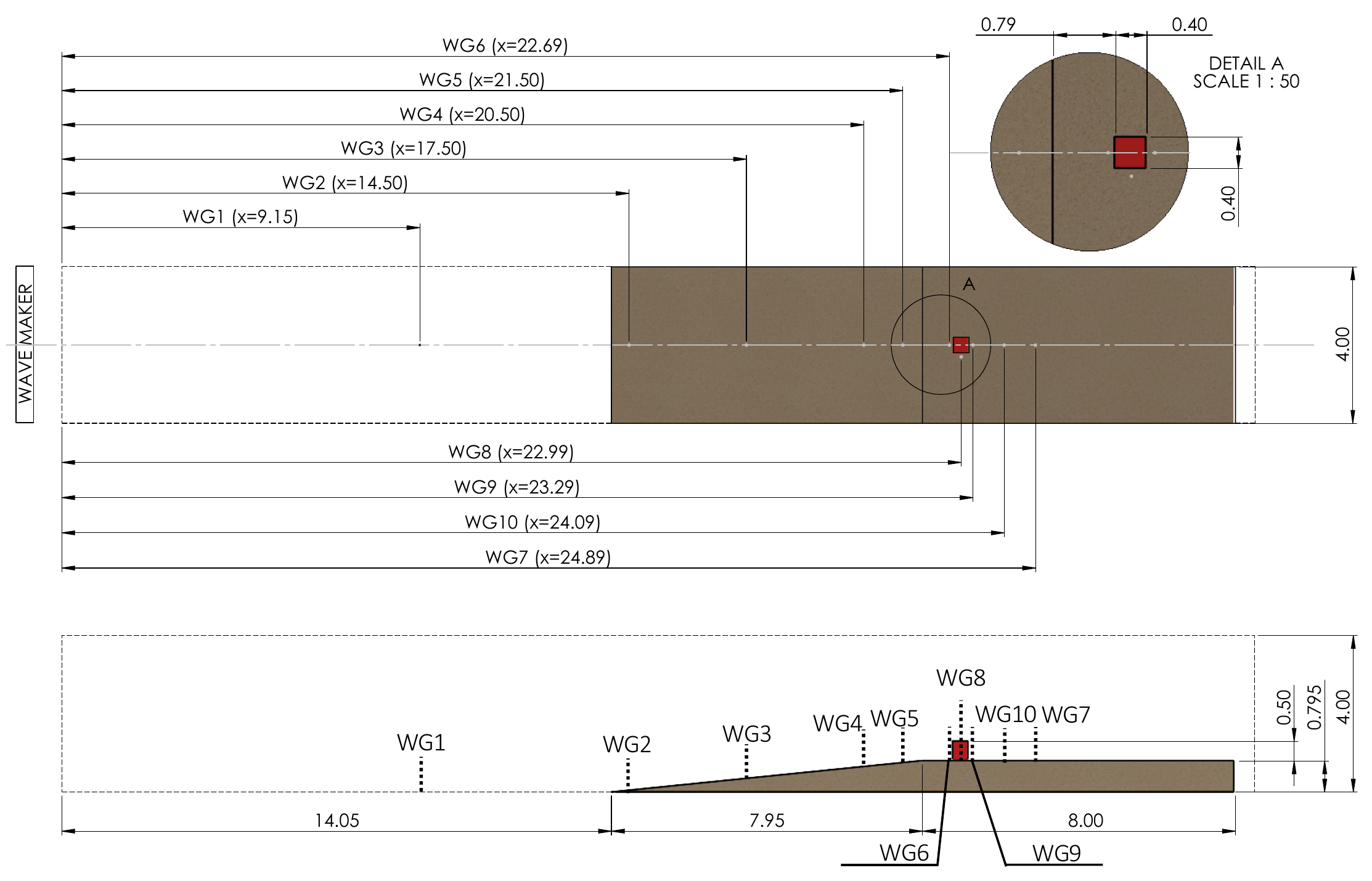}
\caption{HyTOFU experimental set up including the location of wave gauges and the structure of interest. The dimensions are shown in SI units.}
\label{fig:HyTOFUDim}
\end{figure}

The schematic and dimensions of the experimental setup is as shown in \autoref{fig:HyTOFUDim}. The wave flume geometry considered was 32 \unit{\meter} long, 4 \unit{\meter} wide, and 4 \unit{\meter} deep. The wave flume comprised of a flat bed spanning 14.05 \unit{\meter} in length followed by a beach segment of constant slope ratio of 1:10 ($\theta$ = 5.71\unit{\degree}), extending horizontally for another 7.95 \unit{\meter}. Following the slope, the flume setup had an additional 8 \unit{\meter} stretch of an even terrain situated at the height of 0.795 \unit{\meter}. A single structure was situated on that flat terrain, 0.79 \unit{\meter} beyond the end of the slope. The structure considered was of dimensions 0.4 \unit{\meter} $\times$ 0.4 \unit{\meter} $\times$ 0.5 \unit{\meter}. 

The free-surface elevations at ten different locations were measured using wave gauges (WG). The WG 1-6 were located upstream; WG 7, 9, and 10 were located downstream and in close proximity to the back of building; WG 8 was located adjacent to the building. The detailed information on the experimental setup, including information about the instrumentation and working, can be found in \citet{Joetal2021}.

\subsection{SPH simulation setup}
\begin{figure}[!htb]
\centering
\includegraphics[width=0.8\textwidth]{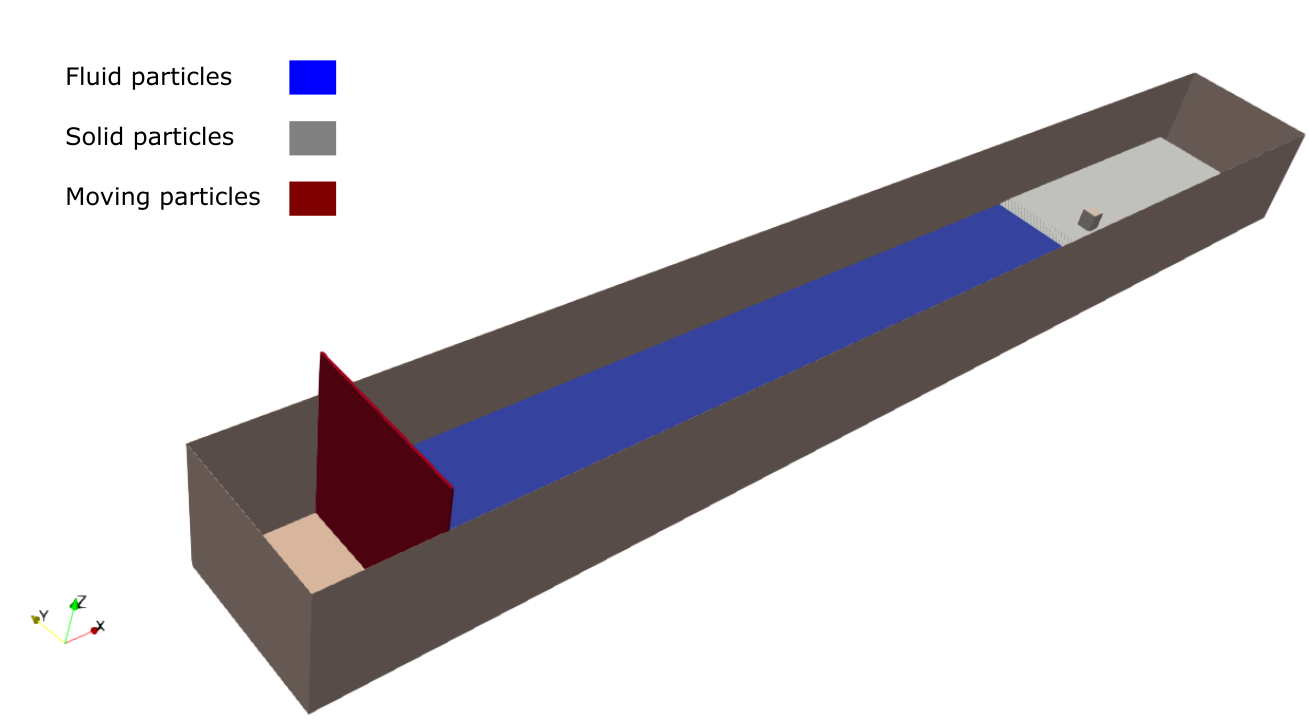}
\caption{DualSPHysics setup of the HyTOFU wave flume. A depiction of the moving wall (red), solid walls (grey) and water (blue).}
\label{fig:HyTOFUSetup}
\end{figure}

A digital twin of the experimental flume has been reconstructed for numerical the simulation using DualSPHysics. The domain, shown in \autoref{fig:HyTOFUSetup} is discretized into wall and fluid particles. The wall particles (grey) use the dynamic boundary conditions (DBC) \cite{cmc2007} and represent the walls and the structure, while the red particles represent the moving particles associated with the piston wavemaker. The blue particles represent the water domain.

The initial water depth is 0.75 \unit{meter} and a solitary wave of 0.40 \unit{\meter} is considered for validation. Further on, solitary waves of 0.40 - 0.90 \unit{meter} are considered for more detailed analysis. For the SPH simulations, the numerical speed of sound is automatically calculated by $c_{0}=10\sqrt{gH_{init}}$ where $H_{init}$ is the water depth at rest. The smoothing length is set to $h=2\mathrm{dp}$ where $\mathrm{dp}$ is the initial particle spacing. Herein, the density diffusion term of \citet{FOURTAKAS2019346} is employed to reduce pressure oscillations and an artificial viscosity of \citet{Mo1992} with $\alpha=0.01$ is used to ensure numerical stability. The initial water density is set to $\rho=1000$ \unit{\kilogram\per\cubic\meter}. A detailed discussion on the method of solitary wave generation can be found in \citet{solitarySPH_2016}.

The upstream of the flume consists of a moving wall and is discretised in SPH using the dynamic boundary conditions (DBC) \cite{cmc2007} with a prescribed boundary velocity. The fixed sidewalls are represented by a set of fixed particles by solving the continuity equation to impose no penetration , and obtain a pressure from \autoref{eqstate}. This approach is computationally efficient, with the density and pressure computed simultaneously and thus resulting in substantial computational time savings. DBC has been used widely in solving coastal engineering problems by discretising complex 3-D geometries without the need for complex mirroring techniques or semi-analytical wall boundary conditions \cite{ALTOMARE201434,zhang2018}. All the walls and building surfaces in this setup are treated with DBC.

\subsection{Numerical convergence}
A convergence study is performed for the initial wave height of $\mathrm{H}=0.4$ \unit{\meter}. The convergence study is considered by varying the initial particle spacing $(\mathrm{dp})$ used in the SPH simulations. Following the work of \citet{ALTOMARE201737, rota2018}, a minimum of four particles is deemed necessary to accurately represent wave height $(\mathrm{H})$, and thus aiming for $\mathrm{H/dp} \ge 4$. The convergence study is considered for $\mathrm{H/dp} = 4, 16, 32$ and the resulting errors and simulation times are tabulated in \autoref{tab:convergence}.
\begin{table}[!htb]
\caption{Numerical convergence}
\centering
\resizebox{\textwidth}{!}{
\begin{tabular}{p{1.25cm} p{1cm} p{2cm} p{2.25cm} p{2.5cm} p{2cm}}
\toprule
\textbf{dp (\unit{\meter})} & \textbf{H/dp} & \textbf{Number of particles} & \textbf{Error in peak height (\%)} & \textbf{Error in wave arrival time (\%)} & \textbf{Simulation time} \\
\midrule
0.1 & 4 & 91791 & 8.87 & -2.20 & 78.51 \unit{\second}\\
0.025 & 16 & 4000681 & 1.96 & -3.14 & 1.66 \unit{\hour}\\
0.0125 & 32 & 29818849 & 1.86 & 1.57 & 30.7 \unit{\hour} \\
\bottomrule
\end{tabular}
}
\label{tab:convergence}
\end{table}

\begin{figure}[!htb]
    \begin{subfigure}{0.44\textwidth}
    \centering
    \adjustbox{width=\linewidth}{
    \input{Section3/images/convergence/wg1c}} 
    \caption{WG1}
    \end{subfigure}
    \begin{subfigure}{0.44\textwidth}
        \centering
        \adjustbox{width=\linewidth}{
\begin{tikzpicture}

\definecolor{darkgray176}{RGB}{176,176,176}
\definecolor{darkorange2491156}{RGB}{249,115,6}

\begin{axis}[
    tick align=outside,
    tick pos=left,
    x grid style={darkgray176},
    xlabel={$\mathrm{dp}$},
    xmin=-0.31, xmax=2.66,
    xtick style={color=black},
    xtick={0.175,1.175,2.175},
    xticklabels={0.1,0.025,0.0125},
    y grid style={darkgray176},
    ylabel={\% error},
    ymin=-15, ymax=15,
    ytick style={color=black}
]
\draw[draw=black,fill=blue] (axis cs:-0.175,0) rectangle (axis cs:0.175,7.68402090537446);
\addlegendimage{ybar,ybar legend,draw=black,fill=blue}

\draw[draw=black,fill=blue] (axis cs:0.825,0) rectangle (axis cs:1.175,-1.37632977034071);
\draw[draw=black,fill=blue] (axis cs:1.825,0) rectangle (axis cs:2.175,-1.37632977034071);
\draw[draw=black,fill=darkorange2491156] (axis cs:0.175,0) rectangle (axis cs:0.525,-6.73637403125397);
\addlegendimage{ybar,ybar legend,draw=black,fill=darkorange2491156}

\draw[draw=black,fill=darkorange2491156] (axis cs:1.175,0) rectangle (axis cs:1.525,-5.70855037479355);
\draw[draw=black,fill=darkorange2491156] (axis cs:2.175,0) rectangle (axis cs:2.525,1.91716427391691);
\end{axis}

\end{tikzpicture}} 
        \caption{Error at WG1}
    \end{subfigure}
    \begin{subfigure}{0.44\textwidth}
        \centering
        \adjustbox{width=1\linewidth}{
        \input{Section3/images/convergence/wg2c}} 
        \caption{WG2}
    \end{subfigure}
    \begin{subfigure}{0.44\textwidth}
        \centering
        \adjustbox{width=\linewidth}{
\begin{tikzpicture}

\definecolor{darkgray176}{RGB}{176,176,176}
\definecolor{darkorange2491156}{RGB}{249,115,6}

\begin{axis}[
tick align=outside,
tick pos=left,
x grid style={darkgray176},
xlabel={$\mathrm{dp}$},
xmin=-0.31, xmax=2.66,
xtick style={color=black},
xtick={0.175,1.175,2.175},
xticklabels={0.1,0.025,0.0125},
y grid style={darkgray176},
ylabel={\% error},
ymin=-15, ymax=15,
ytick style={color=black}
]
\draw[draw=black,fill=blue] (axis cs:-0.175,0) rectangle (axis cs:0.175,13.7384066587396);
\addlegendimage{ybar,ybar legend,draw=black,fill=blue}

\draw[draw=black,fill=blue] (axis cs:0.825,0) rectangle (axis cs:1.175,5.11296076099879);
\draw[draw=black,fill=blue] (axis cs:1.825,0) rectangle (axis cs:2.175,5.11296076099879);
\draw[draw=black,fill=darkorange2491156] (axis cs:0.175,0) rectangle (axis cs:0.525,-2.08868778280542);
\addlegendimage{ybar,ybar legend,draw=black,fill=darkorange2491156}

\draw[draw=black,fill=darkorange2491156] (axis cs:1.175,0) rectangle (axis cs:1.525,-3.53085972850679);
\draw[draw=black,fill=darkorange2491156] (axis cs:2.175,0) rectangle (axis cs:2.525,1.5379185520362);
\end{axis}

\end{tikzpicture}} 
        \caption{Error at WG2}
    \end{subfigure}
    \begin{subfigure}{0.44\textwidth}
        \centering
        \adjustbox{width=\linewidth}{
        \input{Section3/images/convergence/wg3c}} 
        \caption{WG3}
    \end{subfigure}
    \begin{subfigure}{0.44\textwidth}
        \centering
        \adjustbox{width=\linewidth}{
\begin{tikzpicture}

\definecolor{darkgray176}{RGB}{176,176,176}
\definecolor{darkorange2491156}{RGB}{249,115,6}

\begin{axis}[
tick align=outside,
tick pos=left,
x grid style={darkgray176},
xlabel={$\mathrm{dp}$},
xmin=-0.31, xmax=2.66,
xtick style={color=black},
xtick={0.175,1.175,2.175},
xticklabels={0.1,0.025,0.0125},
y grid style={darkgray176},
ylabel={\% error},
ymin=-15, ymax=15,
ytick style={color=black}
]
\draw[draw=black,fill=blue] (axis cs:-0.175,0) rectangle (axis cs:0.175,8.48648648648651);
\addlegendimage{ybar,ybar legend,draw=black,fill=blue}

\draw[draw=black,fill=blue] (axis cs:0.825,0) rectangle (axis cs:1.175,1.96560196560196);
\draw[draw=black,fill=blue] (axis cs:1.825,0) rectangle (axis cs:2.175,1.96560196560196);
\draw[draw=black,fill=darkorange2491156] (axis cs:0.175,0) rectangle (axis cs:0.525,-5.01630094043888);
\addlegendimage{ybar,ybar legend,draw=black,fill=darkorange2491156}

\draw[draw=black,fill=darkorange2491156] (axis cs:1.175,0) rectangle (axis cs:1.525,-3.44827586206896);
\draw[draw=black,fill=darkorange2491156] (axis cs:2.175,0) rectangle (axis cs:2.525,1.56614420062696);
\end{axis}

\end{tikzpicture}} 
        \caption{Error at WG3}
    \end{subfigure}
    \begin{subfigure}{\textwidth}
        \centering
        \includegraphics[width=0.9\columnwidth]{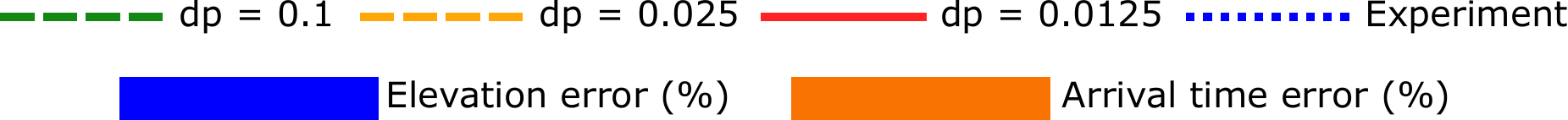}
        \label{fig3:legend} 
    \end{subfigure}
    \caption{Mesh convergence study of SPH results for WG 1-3: normalised free surface elevation obtained from SPH simulations are compared with experimental results from \citet{Joetal2021} for different inter-particle distances ($\mathrm{dp}$) (left); resulting error, in comparison to experiments (right). $\eta$ and $\eta_0$ represent the wave height and characteristic wave height ($\eta_0$ = 0.4 \unit{\meter}), respectively. $t$ represents simulation time (12 \unit{\second} in total) and $T_0$ represents the characteristic time ($T_0 = 2.747$ \unit{\second}).}
    \label{fig:convg}
\end{figure}

Initially, a value of $\mathrm{dp} = 0.1$ \unit{\meter} is chosen, resulting in four particles discretising the wave height. This results in an error of 8.87\% in the peak wave height and -2.2\% in the wave arrival time. Subsequently, the inter-particle distance was reduced to $\mathrm{dp} = 0.025$ \unit{\meter}, resulting in a $\mathrm{H/dp} = 16$. This refinement was reflected in a reduced error of 1.96\% in peak height. However, the wave arrival time continued to exhibit a deviation of -3.14\% compared to the experimental data. Further reduction in the inter-particle distance to $\mathrm{dp} = 0.0125$ \unit{\meter} results in a $\mathrm{H/dp} = 32$. This level of refinement employing $\mathrm{dp} = 0.0125$ \unit{\meter} demonstrated superior agreement with the experimental data. It exhibited an error of 1.86\% in peak height and 1.57\% in wave arrival time. The total simulation time for this refined setup was 30.7 \unit{\hour}. While there is not substantial change in the accuracy between $\rm{dp} = 0.025$ \unit{\meter} vs. $\rm{dp} = 0.0125$ \unit{\meter}, the latter results in at least 32 particles across the dimension of the structure and thus facilitates more accurate evaluation of the forces on the structure. The SPH results, i.e. free surface elevation, for varying inter particle distances, are compared with experimental results \cite{Joetal2021} in \autoref{fig:convg}. As illustrated in \autoref{fig:convg}, the dimensionless wave heights measured at WG1-3 are presented, where $\eta$ represents the wave height, and $\eta_0$ represents the characteristic wave height, which was selected as 0.40 \unit{\meter}. $t$ represents the time, and $T_0$ indicates the characteristic time, determined by the ratio of characteristic wavelength ($L$) to characteristic velocity ($C$), i.e., $T=L/C$. When the wave paddle generates 0.40 \unit{\meter} wave height, the corresponding wavelength and velocity are selected as the characteristic wavelength and velocity and provided in \autoref{table:wave_maker_generation}. \autoref{fig:convg} also shows the errors in free surface elevation and wave arrival times compared against wave flume experiments \cite{Joetal2021}. As evident, the reduction in the inter-particle distance significantly improves the accuracy in the prediction of the wave arrival times. 

\subsection{Validation}
The validation of the developed model has been conducted through comparison of the SPH results, discussed here, with experiments and the volume of fluid (VOF) simulations. The OpenFOAM (using VOF) simulations and experiments are not conducted as a part of this work but the results presented by \citet{Joetal2021} are used. For validation, the solitary wave of initial height of 0.40 \unit{\meter} is considered.

\begin{figure}[!htb]
\centering
  \begin{subfigure}{0.44\textwidth}
    \centering
    \adjustbox{width=\linewidth}{ \input{Section3/images/validation/wg1n}}
    \caption{WG1}
    \label{fig:wg1n-validation}
  \end{subfigure}
  \begin{subfigure}{0.44\textwidth}
    \centering
\adjustbox{width=\linewidth}{ \input{Section3/images/validation/wg2n}}
    \caption{WG2}
    \label{fig:wg2n}
  \end{subfigure}
  \begin{subfigure}{0.44\textwidth}
    \centering
    \adjustbox{width=\linewidth}{ \input{Section3/images/validation/wg3n}}
    \caption{WG3}
    \label{fig:wg3n}
  \end{subfigure}
 \begin{subfigure}{0.44\textwidth}
    \centering
    \adjustbox{width=\linewidth}{ \input{Section3/images/validation/wg4n}}
    \caption{WG 4}
    \label{fig:wg4n}
  \end{subfigure}
  \begin{subfigure}{0.44\textwidth}
    \centering
    \adjustbox{width=\linewidth}{ \input{Section3/images/validation/wg5n}}
    \caption{WG 5}
    \label{fig:wg5n}
  \end{subfigure}
  \begin{subfigure}{0.44\textwidth}
    \centering
\adjustbox{width=\linewidth}{\input{Section3/images/validation/wg8n}}
      \caption{WG 8}
      \label{fig:wg8n}
    \end{subfigure}  
    \newline
    \begin{subfigure}{0.55\textwidth}
    \centering
    \includegraphics[width=\linewidth]{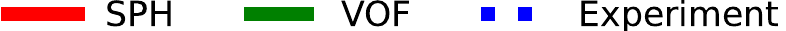}
  \end{subfigure}
  \caption{Validation of SPH simulation through comparison of normalised free surface elevation at the wave gauges (a) WG 1, (b) WG 2, (c) WG 3 (d) WG 4 (e) WG 5 (f) WG 8. $\eta$ and $\eta_0$ represent the wave height and characteristic wave height ($\eta_0$ = 0.4 \unit{\meter}). $t$ represents simulation time (12 \unit{\second} in total) and $T_0$ represents the characteristic time ($T_0 = 2.747$ \unit{\second}). Here the comparison is between SPH (current work) with results from \citet{Joetal2021}, namely VOF (OpenFOAM) and wave flume experiments. The results are compared only for the wage gauges upstream of the structure.}
  \label{fig:validation1}
\end{figure}

\begin{figure}[!htb]
\centering
  \begin{subfigure}{0.44\textwidth}
    \adjustbox{width=\linewidth}{ \input{Section3/images/validation/forces_compare}}
    \label{fig:force-valid}
    \caption{}
  \end{subfigure} 
    \begin{subfigure}{0.44\textwidth}
    \adjustbox{width=\linewidth}{ \input{Section3/images/validation/Force_comparison_normalised_exp_SPH}}
    \label{fig:force_exp_SPH}
    \caption{}
  \end{subfigure}
    \begin{subfigure}{0.44\textwidth}
    \adjustbox{width=\linewidth}{ \input{Section3/images/validation/Force_comparison_normalised_ASCE_Analyse_exp}}
    \label{fig:force_analy_ASCE}
    \caption{}
  \end{subfigure}
    \begin{subfigure}{0.8\textwidth}
    \centering
    \includegraphics[width=\linewidth]{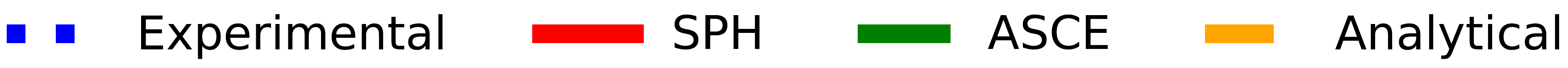}
  \end{subfigure}
\caption{Comparison of non-dimensionalised force-time history acting on the structure. The peak force recorded in experiments \cite{Joetal2021} is used for normalisation. $t$ represents simulation time (12 \unit{\second} in total) and $T_0$ represents the characteristic time ($T_0 = 2.747$ \unit{\second}). Here the comparison is between SPH (current work) with results from \citet{Joetal2021}, namely wave flume experiments, ASCE standard and semi-empirical relations.}
\label{fig:validation3}
\end{figure}

As illustrated in \autoref{fig:validation1}, the normalised free surface elevation in the current work using SPH exhibits excellent agreement with the WG measurements in the initial flat and sloping section (WG 1-3) and satisfactory (WG 4-5) upstream of the building. The wave gauge positioned at the side of the building (WG8) exhibited better agreement with the experiments than VOF simulations. However, a limitation of the SPH approach used here, is that it does not provide variable resolution. A higher resolution can improve the wave gauge readings near to the structure.
\autoref{fig:validation3} shows the comparison of the normalised forces. The comparison uses measured forces during experiments \cite{TOMICZEK201697} with the force calculations outlined earlier from ASCE standards, semi-empirical calculations and the SPH simulations in the current work. 
The force, represented as $F$, is non-dimensionalised using experimental peak force, represented by $F_{\mathrm{exp}}$. The ASCE standards significantly overestimate the forces while both semi-empirical and SPH results underestimate the force-time history and the peak force acting on the structure. However, it is pertinent to note here that the validation is for a case without wave breaking effects.

The overall agreement in terms of wave arrival time and wave elevation with the majority of wave gauges is deemed as a basic validation for the numerical model. Deficiencies of significance are observed in the force calculation. The semi-empirical expressions appear to have lower error, the case considered does not account for wave-breaking effects. Alongside overestimating the peak forces, ASCE standard also provides erroneous rate of reduction from peak value. Considering the flat peak, this is nearly a quasi-static rather than a dynamic loading. Further, the load on the structure does not reduce to zero as rapidly as observed in experiments and other methods. Thus, despite ASCE providing the best probable estimate of the peak force, the overall time-history needs to be considered with caution. SPH method, while largely underestimating the peak force, still captures the rate of change of force (or momentum) reasonably accurately.

Further, it is also pertinent to note that the forces are calculated at one point on the structure and could be further improved through calculation at multiple points and integration. This work is not aimed at improving the SPH method but at coupling with UQ to enable the quantification of uncertainties. The large variation also lends credence to the current work that aims to quantify the uncertainties. Considering the large variation between the forces calculated using the ASCE standard and those obtained from SPH simulations, quantifying these uncertainties is deemed more important than ever.

\section{Results and Discussions}
This section outlines the results and utilises the validated SPH model. The section further outlines the discussions for varying wave heights, the resulting forces and the probabilistic structural response.

\subsection{Solitary wave propagation: Wave heights}
\begin{figure}[!htb]
\centering
  \begin{subfigure}{0.44\textwidth}
    \centering
    \adjustbox{width=\linewidth}{ \input{Section4/images/waveheights/wg1hn}}
    \caption{WG 1}
    \label{fig:wg1hn}
  \end{subfigure}
  \begin{subfigure}{0.44\textwidth}
    \centering
    \adjustbox{width=\linewidth}{ \input{Section4/images/waveheights/wg2hn}}
    \caption{WG 2}
    \label{fig:wg2hn}
  \end{subfigure}
  \begin{subfigure}{0.44\textwidth}
    \centering
    \adjustbox{width=\linewidth}{ \input{Section4/images/waveheights/wg3hn}}
    \caption{WG 3}
    \label{fig:wg3hn}
  \end{subfigure}
  \begin{subfigure}{0.44\textwidth}
    \centering
    \adjustbox{width=\linewidth}{ \input{Section4/images/waveheights/wg4hn}}
    \caption{WG 4}
    \label{fig:wg4hn}
  \end{subfigure}
  \begin{subfigure}{0.44\textwidth}
    \centering
    \adjustbox{width=\linewidth}{ \input{Section4/images/waveheights/wg5hn}}
    \caption{WG 5}
    \label{fig:wg5hn}
  \end{subfigure}
  \begin{subfigure}{0.44\textwidth}
    \centering
    \adjustbox{width=\linewidth}{ \input{Section4/images/waveheights/wg8hn}}
    \caption{WG 8}
    \label{fig:wg8hn}
  \end{subfigure}  
  \begin{subfigure}{\textwidth}
    \centering
    \includegraphics[width=\columnwidth]{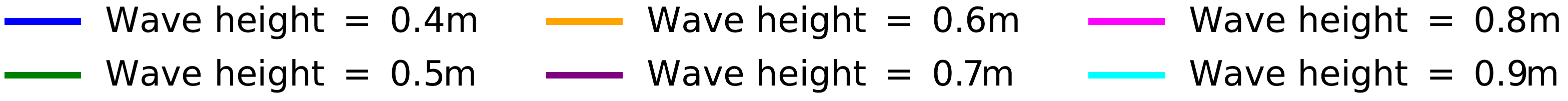}
    \label{fig:legend_wgh01} 
  \end{subfigure}
  \caption{Normalised free surface elevation observed at the wave gauge locations (a) WG 1 (b) WG 2 (c) WG 3 (d) WG 4 (e) WG 5 (f) WG 8 for varying initial wave heights, where $\eta$ represent the calculated wave height at the particular WG, $\eta_0$ represent the characteristic wave height ($\eta_0$ = 0.40 \unit{\meter}), $t$ indicates the simulation time, and $T_0$ represent a characteristic time ($T_0 = 2.747$ \unit{\second}).}
  \label{fig:wgheightn1}
\end{figure}

The SPH setup discussed in the earlier section has been used for the simulation of solitary wave propagation in a wave flume. In this regard, considering same conditions outlined in the validation example, SPH simulations are setup with solitary waves of various initial wave heights, i.e 0.40 - 0.90 \unit{\meter} at an interval of 0.10 \unit{\meter}. Solitary waves at each height are listed in \autoref{table:wave_maker_generation} with the corresponding wave length and celerity information. These solitary waves at different wave heights are generated based on \citet{Rayleigh} wave theory by using the piston-type wavemaker. The normalised wave heights measured at the wave gauges 1 - 5 and 8 are shown in \autoref{fig:wgheightn1}. It is pertinent to note that the solitary wave with initial wave height of 0.40 \unit{\meter} has the least loss to wave height as it moves through the wave flume. The wave with initial wave height of 0.60 \unit{\meter} or more nearly lose about 50\% of wave height by the time the wave reaches the location of WG 5.
\begin{table}[!htb]
  \centering
  \caption{Wave paddles configuration}
  \resizebox{\textwidth}{!}{%
  \begin{tabular}{cccccc}
    \hline
    \text{Wave height (\unit{\meter})} & \text{Wave length (\unit{\meter})} & \text{Celerity (\unit{\meter\per\second})} & \text{Wave height (\unit{\meter})} & \text{Wave length (\unit{\meter})} & \text{Celerity (\unit{\meter\per\second})} \\
    \hline
    0.4 & 9.22634 & 3.35879 & 0.7 & 7.83151 & 3.77154 \\
    0.5 & 8.60361 & 3.50179 & 0.8 & 7.57411 & 3.89942 \\
    0.6 & 8.16210 & 3.63916 & 0.9 & 7.36769 & 4.02324 \\
    \hline
  \end{tabular}}
  \label{table:wave_maker_generation}
\end{table}

Appreciable wave heights are not measured at the wave gauges 7, 9 and 10 and thus not discussed here. However, only the waves with initial wave height of 0.80 and 0.90 \unit{\meter} result in substantial over topping.

\subsection{Forces on structures \label{sec:forces}}
\begin{figure}[!htb]
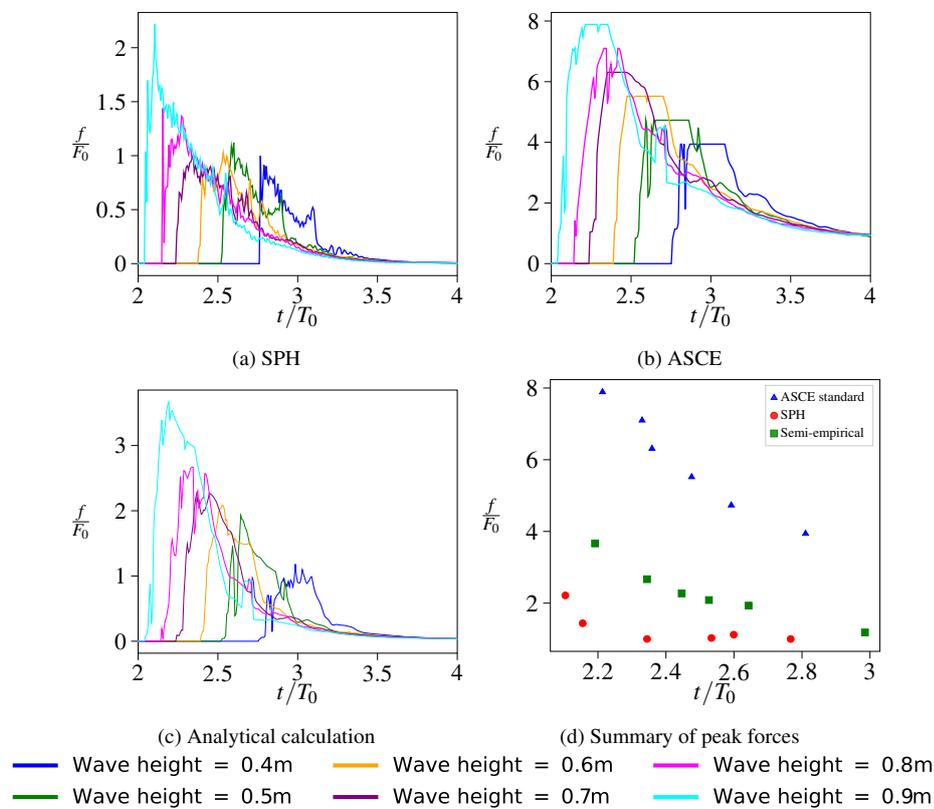

  \centering
   \begin{subfigure}{0.44\textwidth}
    \centering
    \adjustbox{width=\linewidth}{ \input{Section4/images/force1/forces_wh}}
    \caption{SPH}
    \label{fig:forcesph}
  \end{subfigure} 
  \begin{subfigure}{0.44\textwidth}
    \centering
    \adjustbox{width=\linewidth}{ \input{Section4/images/force1/forces_asce}}
    \caption{ASCE}
    \label{fig:forceasce}
  \end{subfigure} 
  \begin{subfigure}{0.44\textwidth}
    \centering
    \adjustbox{width=\linewidth}{ \input{Section4/images/force1/forces_analytical}}
    \caption{Analytical calculation}
    \label{fig:forceanalytical}
  \end{subfigure}
    \begin{subfigure}{0.44\textwidth}
    \centering
    \adjustbox{width=\linewidth}{
\begin{tikzpicture}

\definecolor{darkgray176}{RGB}{176,176,176}
\definecolor{green}{RGB}{0,128,0}
\definecolor{lightgray204}{RGB}{204,204,204}

\begin{axis}[
legend cell align={left},
legend style={fill opacity=0.8, draw opacity=1, text opacity=1, draw=lightgray204},
tick align=outside,
tick pos=left,
x grid style={darkgray176},
xlabel={$t/T_0$},
xmin=2.0601373131143, xmax=3.02921255825831,
xtick style={color=black},
y grid style={darkgray176},
ylabel style={rotate=-90},
ylabel={$\frac{f}{F_0}$},
ymin=0.655574676258993, ymax=8.23293179856115,
ytick style={color=black}
]
\addplot [draw=blue, fill=blue, mark=triangle*, only marks]
table{%
x  y
2.81042348504499 3.94030503597122
2.59201734929127 4.72994532374101
2.47549827926892 5.51958561151079
2.35900978890871 6.30922589928057
2.32989540224509 7.09886618705036
2.21340618379769 7.88850647482015
};
\addlegendentry{ASCE standard}
\addplot [draw=red, fill=red, mark=*, only marks]
table{%
x  y
2.76673388489152 1
2.59927455839585 1.12086330935252
2.53375799630189 1.02877697841727
2.34445532590388 1.00143884892086
2.15515010720069 1.43884892086331
2.10418618789357 2.2158273381295
};
\addlegendentry{SPH}
\addplot [draw=green, fill=green, mark=square*, only marks]
table{%
x  y
2.98516368347904 1.1826025901833
2.64298490903435 1.93128927916118
2.52646875136078 2.08634828895761
2.446395906044 2.27110673694666
2.3444411282036 2.66836928408702
2.19156283989534 3.66635912602927
};
\addlegendentry{Semi-empirical}
\end{axis}

\end{tikzpicture}}
    \caption{Summary of peak forces}
    \label{fig:force-points}
  \end{subfigure}
    \begin{subfigure}{\textwidth}
    \centering
    \includegraphics[width=\columnwidth]{Section4/images/waveheights/wgh_legend.pdf}
    \label{fig:legend_forces} 
  \end{subfigure}
  \caption{Normalised force time history obtained from (a) SPH (b) ASCE standards (c) semi-empirical expressions (d) summary of peak forces. The summary extracts only the peak values from the SPH, ASCE and semi-empirical results. $f$ represents the force measured/calculated, $F_0$ represent the characteristic force  ($F_0 = 695$ \unit{\newton}), $t$ indicates the simulation time, and $T_0$ represent the characteristic time ($T_0 = 2.747$ \unit{\second}).}
  \label{fig:forceComparison}
\end{figure}

The forces are obtained through post-processing from SPH simulations and compared with the ASCE regulations and analytical expressions. The non-dimensionalised drag force time histories for varying initial wave heights from 0.40 - 0.90 \unit{\meter} are shown in \autoref{fig:forceComparison}. The peak force exerted on the structure where the wave height is 0.40 \unit{\meter}, is chosen as the characteristic force ($F_0$ = 695 \unit{\newton}). It is pertinent to note here that the ASCE standard and semi-empirical calculations consider the wave heights and velocities at the front of the structure. However, in contrast, SPH calculations consider all the faces of the structure in the calculation of forces and thus better represents the averaged overall force on the building \cite{dom2021}. However, in the current form, used in this work, SPH uses an offset distance to calculate these forces causing slight, but acceptable, inaccuracy in the peak forces. 

As seen in \autoref{fig:forceComparison}, the forces predicted by ASCE standards are on the higher end, while SPH is on the lower end. This has also been discussed in the earlier section on validation in the discussion related to \autoref{fig:validation3}. The forces predicted using ASCE standards are about four times that predicted by SPH. This can be attributed to the overall area on which the forces act. While semi-empirical expressions demonstrate a better match with experiments, the same might not be translatable to all scenarios, in particular to breaking waves. Thus, semi-empirical solutions have limited applicability and needs to be considered carefully.

Further, the resulting normalised force-time plot, depicted in \autoref{fig:wgheightn1}, reveals a behaviour in contrast to the general assumption of direct proportionality between wave height and maximum drag force. The SPH simulations indicate that the wave loading is not strictly proportional to wave height, as seen through ASCE or semi-empirical calculations. While the wave arrival times are as expected across the board, with the 0.40 \unit{\meter} wave arriving much later than a 0.90 \unit{\meter} solitary wave, the maximum force does not increase proportionally as wave height changes from 0.40 - 0.90 \unit{\meter}. As the wave height increases from 0.40 to 0.50 \unit{\meter}, the peak forces increase; further it decreases for wave heights 0.60 and 0.70 \unit{\meter}; to drastically increase again for wave heights 0.80 and 0.90 \unit{\meter}.

From \autoref{fig:wgheightn1}, it was observed that the solitary wave of 0.40 \unit{\meter} continues to preserve most of the wave shape as it approaches WG 5; however, the wave heights still demonstrate a linear variation at WG 5. This is further substantiated in \autoref{fig:waveheight4} where the wave of height 0.40 \unit{\meter} can be characterised as a surging or collapsing wave. Similarly, \autoref{fig:waveheight7} shows the wave as a spilling breaker and \autoref{fig:waveheight9} as potentially a plunging breaker. The waves with wave heights of 0.80 \unit{\meter} and 0.90 \unit{\meter} also show significant over-topping effects. Further, \autoref{fig:waveheight4} - \autoref{fig:waveheight9} also demonstrate that the wave with height 0.40 \unit{\meter} is a non-breaking wave; 0.70 \unit{\meter} has a spilling effect; 0.90 \unit{\meter} breaks as early as WG1-WG2 but transports a large volume of water. The varied behaviour can, thus, be attributed to the nonlinear response of observed force vs. wave heights.

\subsection{Uncertainty quantification}
Probabilistic structural analysis is performed to identify the range of structural responses considering the uncertainties in both structural properties and the wave heights. While the probabilistic analysis is computationally expensive compared to deterministic simulations, it provides more reliable estimates for practical engineering applications \cite{Judd2018Windstorm, Mohammadi2019Performance, Ghaffary2021Performance, Maraveas2019Assessment}.

Uncertainties in structural responses typically arise from two primary sources: the parameters associated with the structural behaviour and those related to the wave mechanics. The former includes basic attributes such as mass, stiffness, damping of the structure and more intricate characteristics such as the modulus of elasticity and yield strength of the building's individual components. The latter can be attributed as a variation in the wave forces, i.e wave heights, itself. This work leverages the developed NHERI SimCenter UQ engine \cite{Deierlein2020cloud, McKenna2022NHERI}. 

\begin{figure}[hbt!]
    \centering
    \includegraphics[width=0.5\columnwidth]{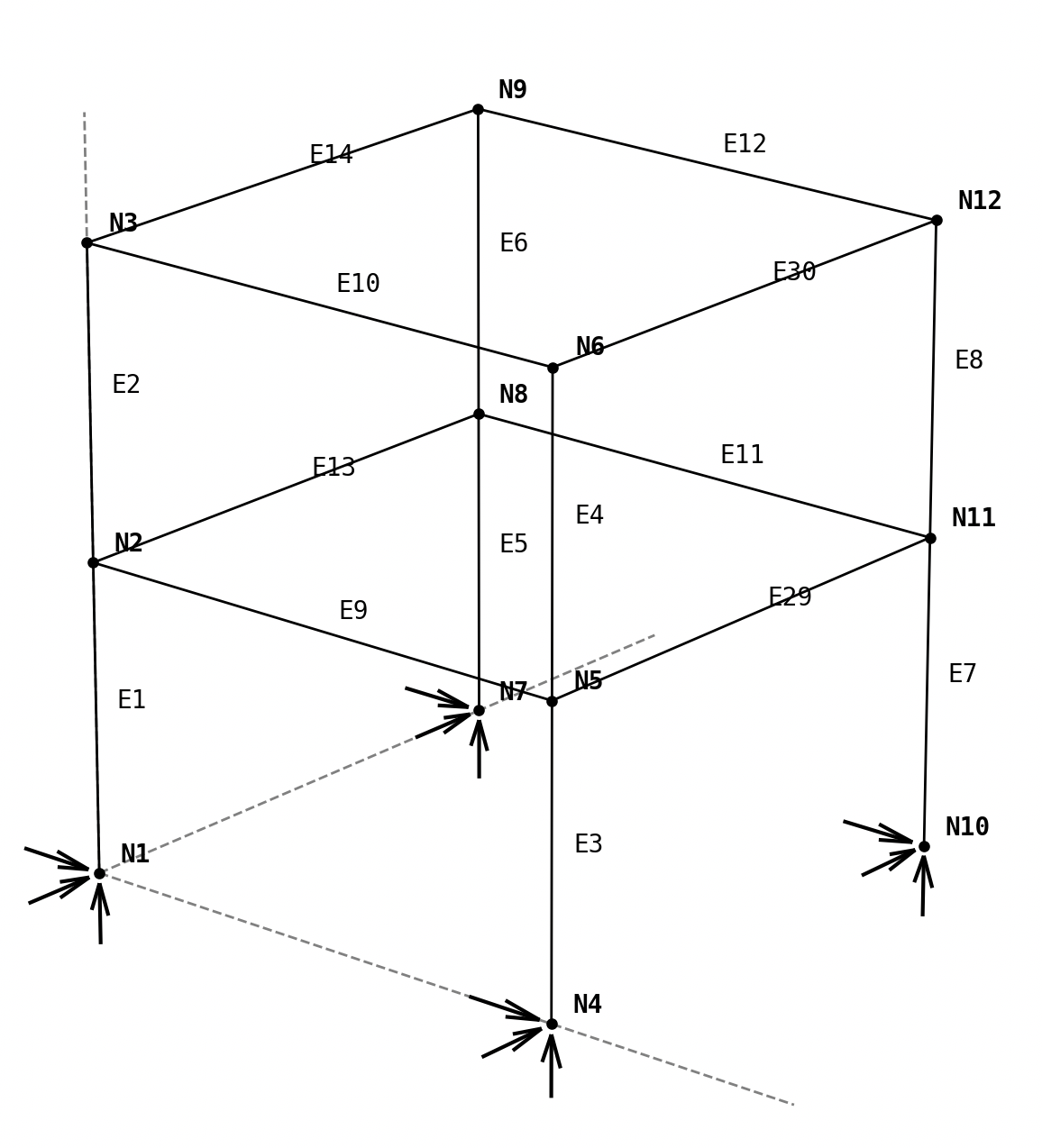}
    \caption{OpenSees beam-column model, constructed by W-section steel, and used in the structural analysis}
    \label{fig:W-section model}
\end{figure}
An OpenSees beam-column model, as shown in \autoref{fig:W-section model}, is used to represent a two-storey structure consisting of column, beam, and girder sections, constructed using steel W-sections. The loading on the structure is derived from the wave loading history applied at each centroid node (i.e. storey level). This force history is obtained from earlier discussed methods: ASCE standard, semi-empirical expressions and SPH. 

\subsubsection{Sample generation}
Considering the inputs as a random variable, a forward UQ method, namely LHS \cite{quoFEM}, is used for sample generation and aid in quantifying the propagation of the uncertainties. The engineering demand parameters (EDP) define the quantities of interest. In this work, peak floor displacement and root mean square acceleration (RMSA) are considered while additional user-defined EDP's can also be added.

\begin{table}[hbt!]
  \caption{List of random variables for two-storey structure}
  \centering
  \resizebox{\textwidth}{!}{%
  \begin{tabular}{p{6cm}>{\centering\arraybackslash}p{2cm}>{\centering\arraybackslash}p{3cm}>{\centering\arraybackslash}p{2cm}}
    \toprule
    \textbf{Variable} & \textbf{Distribution type} & \textbf{Mean} & \textbf{Standard deviation} \\
    \midrule
    Yield strength (\unit{\mega\pascal}) & Normal & 413.685 & 82 \\
    W-section weight per length of column (\unit{\newton\per\meter}) & Normal & 173.4 & 34 \\
    W-section weight per length of beam (\unit{\newton\per\meter}) & Normal & 133.554 & 26 \\
    W-section weight per length of girder (\unit{\newton\per\meter}) & Normal & 133.554 & 26 \\
    Young's modulus for steel (\unit{\giga\pascal}) & Normal & 200 & 40 \\
    \bottomrule
  \end{tabular}}
  \label{tab:random_variables_2stories}
\end{table}
The LHS method is utilised to generate 600 samples. The structural random variables used in this work are specified in \autoref{tab:random_variables_2stories}. A 20\% coefficient of variation is considered to evaluate the probabilistic structural response. Further, uniform distribution is assumed for the wave parameters, i.e. wave heights over the discretised domain of 0.40 - 0.90 \unit{\meter}, at intervals of 0.10 \unit{\meter}.

\begin{figure}[hbt!]
 \centering
   \begin{subfigure}{0.46\textwidth}
     \centering
     \adjustbox{width=\linewidth}{
\begin{tikzpicture}

\definecolor{darkgray148163195}{RGB}{148,163,195}
\definecolor{darkgray176}{RGB}{176,176,176}
\definecolor{darksalmon232149117}{RGB}{232,149,117}
\definecolor{dimgray88}{RGB}{88,88,88}
\definecolor{mediumaquamarine113182160}{RGB}{113,182,160}

\begin{axis}[
tick align=outside,
tick pos=left,
x grid style={darkgray176},
xmin=-0.5, xmax=2.5,
xtick style={color=black},
xtick={0,1,2},
xticklabels={Column,Beam,Girder},
y grid style={darkgray176},
ylabel={\(\displaystyle \mathrm{Weights\ (\unit{\newton\per\meter})}\)},
ymin=40.246485, ymax=278.225215,
ytick style={color=black}
]
\path [draw=dimgray88, fill=mediumaquamarine113182160, very thick]
(axis cs:-0.4,152.64)
--(axis cs:0.4,152.64)
--(axis cs:0.4,193.92375)
--(axis cs:-0.4,193.92375)
--(axis cs:-0.4,152.64)
--cycle;
\addplot [very thick, dimgray88]
table {%
0 152.64
0 113.905
};
\addplot [very thick, dimgray88]
table {%
0 193.92375
0 244.092
};
\addplot [very thick, dimgray88]
table {%
-0.2 113.905
0.2 113.905
};
\addplot [very thick, dimgray88]
table {%
-0.2 244.092
0.2 244.092
};
\addplot [black, mark=o, mark size=3, mark options={solid,fill opacity=0,draw=dimgray88}, only marks]
table {%
0 80.7739
0 267.408
};
\path [draw=dimgray88, fill=darksalmon232149117, very thick]
(axis cs:0.6,114.7495)
--(axis cs:1.4,114.7495)
--(axis cs:1.4,148.8445)
--(axis cs:0.6,148.8445)
--(axis cs:0.6,114.7495)
--cycle;
\addplot [very thick, dimgray88]
table {%
1 114.7495
1 72.0999
};
\addplot [very thick, dimgray88]
table {%
1 148.8445
1 197.268
};
\addplot [very thick, dimgray88]
table {%
0.8 72.0999
1.2 72.0999
};
\addplot [very thick, dimgray88]
table {%
0.8 197.268
1.2 197.268
};
\addplot [black, mark=o, mark size=3, mark options={solid,fill opacity=0,draw=dimgray88}, only marks]
table {%
1 62.6645
};
\path [draw=dimgray88, fill=darkgray148163195, very thick]
(axis cs:1.6,115.22625)
--(axis cs:2.4,115.22625)
--(axis cs:2.4,153.65725)
--(axis cs:1.6,153.65725)
--(axis cs:1.6,115.22625)
--cycle;
\addplot [very thick, dimgray88]
table {%
2 115.22625
2 67.8814
};
\addplot [very thick, dimgray88]
table {%
2 153.65725
2 196.106
};
\addplot [very thick, dimgray88]
table {%
1.8 67.8814
2.2 67.8814
};
\addplot [very thick, dimgray88]
table {%
1.8 196.106
2.2 196.106
};
\addplot [black, mark=o, mark size=3, mark options={solid,fill opacity=0,draw=dimgray88}, only marks]
table {%
2 51.0637
};
\addplot [very thick, dimgray88]
table {%
-0.4 172.232
0.4 172.232
};
\addplot [very thick, dimgray88]
table {%
0.6 137.698
1.4 137.698
};
\addplot [very thick, dimgray88]
table {%
1.6 131.9125
2.4 131.9125
};
\end{axis}

\end{tikzpicture}}
     \label{fig:9a}
     \caption{}
   \end{subfigure}
   \hfill
   \begin{subfigure}{0.46\textwidth}
     \centering
     \adjustbox{width=\linewidth}{
\begin{tikzpicture}

\definecolor{darkgray176}{RGB}{176,176,176}
\definecolor{darksalmon232149117}{RGB}{232,149,117}
\definecolor{dimgray88}{RGB}{88,88,88}
\definecolor{mediumaquamarine113182160}{RGB}{113,182,160}

\begin{axis}[
    tick align=outside,
    tick pos=left,
    x grid style={darkgray176},
    xmin=-0.5, xmax=1.5,
    xtick style={color=black},
    xtick={0,1},
    xticklabels={Yield strength,Young's modulus},
    y grid style={darkgray176},
    ymin=0, ymax=2,
    ylabel style={rotate=-90},
    ylabel={\shortstack{$\frac{\sigma}{\sigma_0}$ \\ $\frac{E}{E_0}$}},
    ytick style={color=black}
]
\path [draw=dimgray88, fill=mediumaquamarine113182160, very thick]
(axis cs:-0.4,0.878907864679648)
--(axis cs:0.4,0.878907864679648)
--(axis cs:0.4,1.14120284757726)
--(axis cs:-0.4,1.14120284757726)
--(axis cs:-0.4,0.878907864679648)
--cycle;
\addplot [very thick, dimgray88]
table {%
0 0.878907864679648
0 0.575271039559085
};
\addplot [very thick, dimgray88]
table {%
0 1.14120284757726
0 1.39377545717152
};
\addplot [very thick, dimgray88]
table {%
-0.2 0.575271039559085
0.2 0.575271039559085
};
\addplot [very thick, dimgray88]
table {%
-0.2 1.39377545717152
0.2 1.39377545717152
};
\addplot [black, mark=o, mark size=3, mark options={solid,fill opacity=0,draw=dimgray88}, only marks]
table {%
0 0.4584285144494
};
\path [draw=dimgray88, fill=darksalmon232149117, very thick]
(axis cs:0.6,0.86855)
--(axis cs:1.4,0.86855)
--(axis cs:1.4,1.15576625)
--(axis cs:0.6,1.15576625)
--(axis cs:0.6,0.86855)
--cycle;
\addplot [very thick, dimgray88]
table {%
1 0.86855
1 0.55714
};
\addplot [very thick, dimgray88]
table {%
1 1.15576625
1 1.528185
};
\addplot [very thick, dimgray88]
table {%
0.8 0.55714
1.2 0.55714
};
\addplot [very thick, dimgray88]
table {%
0.8 1.528185
1.2 1.528185
};
\addplot [very thick, dimgray88]
table {%
-0.4 0.982502387081958
0.4 0.982502387081958
};
\addplot [very thick, dimgray88]
table {%
0.6 0.9800275
1.4 0.9800275
};
\end{axis}

\end{tikzpicture}}
     \label{fig:9b}
     \caption{}
   \end{subfigure}
 \caption{Distribution of structural properties considered: (a) column weight, beam weight, girder weight; (b) normalised yield strength ($\sigma_0 = 413.685$ \unit{\mega\pascal}), and normalised Young's modulus ($E_0 = 200$ \unit{\giga\pascal}).}
 \label{fig:Root mean square accelerations_2stories_diffscaling}
\end{figure}
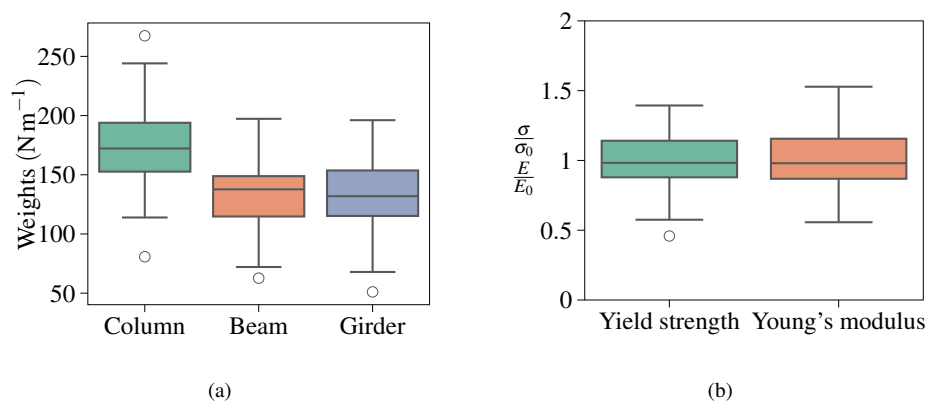

The boxplots in \autoref{fig:Root mean square accelerations_2stories_diffscaling} depict the distribution of sample values generated using LHS method. It compares the distributions of column weight, beam weight, and girder weight and can be observed that the median of column weight is the largest, while the differences between the medians of beam weight and girder weight are relatively small. This suggests that columns may contribute significantly within the structure.

The relationships between the peak displacements, root mean square accelerations (RMSA) and the random variables, outlined in \autoref{tab:random_variables_2stories}, are illustrated in \autoref{fig:Root mean square accelerations_2stories} and \autoref{fig:Peak displacements_2stories} respectively. 

\subsubsection{Peak floor displacement}
\autoref{fig:myEsl_peakdisp} shows a nonlinear and strong negative correlation between Youngs' modulus and peak displacement in $x$-direction; a strong nonlinear trend in relation to initial wave heights. A general trend demonstrates little correlation between peak displacement and other parameters. 

However, as shown earlier, the increase in wave height does not correspond to a proportional increase in the forces. Thus, the multi-dimensionality of the problem renders it hard to ascertain exact correlations across the spectrum by just using the peak displacement.

\subsubsection{Root mean square acceleration (RMSA)}
\begin{figure}[hbt!]
 \centering
 \begin{subfigure}{0.44\textwidth}
   \centering
   \adjustbox{width=\linewidth}{ \input{Section4/images/uq/OpenSees/myQdl_rmsa}}
   \caption{}
   \label{fig:myQdl_rmsa}
 \end{subfigure}
 \begin{subfigure}{0.44\textwidth}
   \centering
   \adjustbox{width=\linewidth}{ \input{Section4/images/uq/OpenSees/myQBeam_rmsa}}
   \caption{}
   \label{fig:myQBeam_rmsa}
 \end{subfigure}
 \begin{subfigure}{0.44\textwidth}
   \centering
   \adjustbox{width=\linewidth}{ \input{Section4/images/uq/OpenSees/myQGrid_rmsa}}
   \caption{}
   \label{fig:myQGrid_rmsa}
 \end{subfigure}
  \begin{subfigure}{0.44\textwidth}
   \centering
   \adjustbox{width=\linewidth}{ \input{Section4/images/uq/OpenSees/myfy_rmsa}}
   \caption{}
   \label{fig:myfy_rmsa}
 \end{subfigure}
 \begin{subfigure}{0.44\textwidth}
   \centering
   \adjustbox{width=\linewidth}{ \input{Section4/images/uq/OpenSees/myEs_rmsa}}
   \caption{}
   \label{fig:myEsl_rmsa}
 \end{subfigure}
 \begin{subfigure}{0.44\textwidth}
   \centering
   \adjustbox{width=\linewidth}{ 
\begin{tikzpicture}

\definecolor{darkgray176}{RGB}{176,176,176}
\definecolor{green01270}{RGB}{0,127,0}
\definecolor{orange}{RGB}{255,165,0}
\definecolor{purple}{RGB}{128,0,128}

\begin{axis}[
tick align=outside,
tick pos=left,
x grid style={darkgray176},
xlabel={Wave height (\unit{\centi\meter})},
xmin=0.375, xmax=0.925,
xtick style={color=black},
y grid style={darkgray176},
ylabel={RMSA (\unit{\meter\per\second\squared})},
ymin=-0.008939395, ymax=2.382889495,
ytick style={color=black}
]
\addplot [draw=red, fill=red, mark=*, only marks]
table{%
x  y
0.4 0.336456
0.4 0.361967
0.4 0.355815
0.4 0.378438
0.4 0.374188
0.4 0.443858
0.4 0.340281
0.4 0.295202
0.4 0.38174
0.4 0.379262
0.4 0.339699
0.4 0.376042
0.4 0.376962
0.4 0.397902
0.4 0.373767
0.4 0.402882
0.4 0.338422
0.4 0.329943
0.4 0.362398
0.4 0.500213
0.4 0.3368
0.4 0.448734
0.4 0.363075
0.4 0.330938
0.4 0.382568
0.4 0.335546
0.4 0.428779
0.4 0.366051
0.4 0.328084
0.4 0.32661
0.4 0.353735
0.4 0.359008
0.4 0.386627
0.4 0.413564
0.4 0.364333
0.4 0.379774
0.4 0.33496
0.4 0.373965
0.4 0.370676
0.4 0.345859
0.4 0.363078
0.4 0.306503
0.4 0.417837
0.4 0.429729
0.4 0.34451
0.4 0.381555
0.4 0.337831
0.4 0.423804
0.4 0.345547
0.4 0.432836
0.4 0.417289
0.4 0.374005
0.4 0.346956
0.4 0.396686
0.4 0.418967
0.4 0.348647
0.4 0.371111
0.4 0.363327
0.4 0.449813
0.4 0.344406
0.4 0.367013
0.4 0.33368
0.4 0.344803
0.4 0.323807
0.4 0.4219
0.4 0.361129
0.4 0.40015
0.4 0.412762
0.4 0.416392
0.4 0.381624
0.4 0.40331
0.4 0.375106
0.4 0.356734
0.4 0.403824
0.4 0.404669
0.4 0.490042
0.4 0.414797
0.4 0.39854
0.4 0.40594
0.4 0.393493
0.4 0.338883
0.4 0.389392
0.4 0.380307
0.4 0.41215
0.4 0.347813
0.4 0.479093
0.4 0.332459
0.4 0.504986
0.4 0.437118
0.4 0.422163
0.4 0.399567
0.4 0.482189
0.4 0.349746
0.4 0.364102
0.4 0.331184
0.4 0.338823
0.4 0.351694
0.4 0.291977
0.4 0.383615
0.4 0.386461
};
\addplot [draw=green01270, fill=green01270, mark=square*, only marks]
table{%
x  y
0.5 0.281814
0.5 0.289098
0.5 0.254519
0.5 0.242593
0.5 0.278937
0.5 0.2194
0.5 0.274418
0.5 0.229266
0.5 0.271456
0.5 0.384616
0.5 0.244504
0.5 0.252969
0.5 0.325136
0.5 0.235804
0.5 0.250963
0.5 0.248841
0.5 0.214799
0.5 0.246112
0.5 0.25129
0.5 0.288401
0.5 0.264
0.5 0.29428
0.5 0.286394
0.5 0.264031
0.5 0.279251
0.5 0.231256
0.5 0.298054
0.5 0.285824
0.5 0.238966
0.5 0.26672
0.5 0.27992
0.5 0.229294
0.5 0.280901
0.5 0.251645
0.5 0.290384
0.5 0.233435
0.5 0.243231
0.5 0.245071
0.5 0.316947
0.5 0.280868
0.5 0.266078
0.5 0.269916
0.5 0.254178
0.5 0.300843
0.5 0.264209
0.5 0.252545
0.5 0.242448
0.5 0.245385
0.5 0.280553
0.5 0.208024
0.5 0.256204
0.5 0.233189
0.5 0.302034
0.5 0.246844
0.5 0.240203
0.5 0.307495
0.5 0.249933
0.5 0.256905
0.5 0.26119
0.5 0.280406
0.5 0.293444
0.5 0.21242
0.5 0.249059
0.5 0.221924
0.5 0.239172
0.5 0.23235
0.5 0.285214
0.5 0.264837
0.5 0.26906
0.5 0.246381
0.5 0.263091
0.5 0.244851
0.5 0.235179
0.5 0.287773
0.5 0.227642
0.5 0.261991
0.5 0.250942
0.5 0.281405
0.5 0.26037
0.5 0.296117
0.5 0.26784
0.5 0.243469
0.5 0.275904
0.5 0.245527
0.5 0.290517
0.5 0.247028
0.5 0.297376
0.5 0.305324
0.5 0.221984
0.5 0.255537
0.5 0.290342
0.5 0.276906
0.5 0.269959
0.5 0.23555
0.5 0.28675
0.5 0.265233
0.5 0.295783
0.5 0.284404
0.5 0.253504
0.5 0.27981
};
\addplot [draw=blue, fill=blue, mark=triangle*, only marks]
table{%
x  y
0.6 0.147658
0.6 0.115919
0.6 0.114039
0.6 0.118135
0.6 0.126304
0.6 0.147202
0.6 0.112009
0.6 0.117556
0.6 0.123721
0.6 0.133444
0.6 0.134823
0.6 0.125321
0.6 0.104455
0.6 0.121496
0.6 0.13289
0.6 0.122999
0.6 0.134879
0.6 0.108341
0.6 0.0997801
0.6 0.122161
0.6 0.115941
0.6 0.130292
0.6 0.110384
0.6 0.146363
0.6 0.128568
0.6 0.126741
0.6 0.115614
0.6 0.108601
0.6 0.119769
0.6 0.134767
0.6 0.120652
0.6 0.105429
0.6 0.14225
0.6 0.118969
0.6 0.115511
0.6 0.108954
0.6 0.10911
0.6 0.125701
0.6 0.138221
0.6 0.143023
0.6 0.142863
0.6 0.122719
0.6 0.110214
0.6 0.120225
0.6 0.11845
0.6 0.112564
0.6 0.123094
0.6 0.107352
0.6 0.146137
0.6 0.107148
0.6 0.120806
0.6 0.10805
0.6 0.119303
0.6 0.126506
0.6 0.111607
0.6 0.104165
0.6 0.13655
0.6 0.105375
0.6 0.103351
0.6 0.135696
0.6 0.121843
0.6 0.136447
0.6 0.123057
0.6 0.129685
0.6 0.108948
0.6 0.108492
0.6 0.118429
0.6 0.106745
0.6 0.111536
0.6 0.132135
0.6 0.113001
0.6 0.120807
0.6 0.114725
0.6 0.118832
0.6 0.125101
0.6 0.130583
0.6 0.134727
0.6 0.12944
0.6 0.11487
0.6 0.107254
0.6 0.131769
0.6 0.119124
0.6 0.104477
0.6 0.157555
0.6 0.120556
0.6 0.1176
0.6 0.12609
0.6 0.124735
0.6 0.133282
0.6 0.134858
0.6 0.14332
0.6 0.126311
0.6 0.171071
0.6 0.121405
0.6 0.145196
0.6 0.109329
0.6 0.118806
0.6 0.105804
0.6 0.155591
0.6 0.120227
};
\addplot [draw=purple, fill=purple, mark=x, only marks]
table{%
x  y
0.7 0.171227
0.7 0.149932
0.7 0.194578
0.7 0.154717
0.7 0.146981
0.7 0.152024
0.7 0.197845
0.7 0.17268
0.7 0.18606
0.7 0.155665
0.7 0.184304
0.7 0.164052
0.7 0.140104
0.7 0.14343
0.7 0.124572
0.7 0.148495
0.7 0.164992
0.7 0.180158
0.7 0.149359
0.7 0.146678
0.7 0.130834
0.7 0.173279
0.7 0.173895
0.7 0.170225
0.7 0.130979
0.7 0.165475
0.7 0.14774
0.7 0.167138
0.7 0.17546
0.7 0.150243
0.7 0.130739
0.7 0.166671
0.7 0.151519
0.7 0.193849
0.7 0.155186
0.7 0.137208
0.7 0.183101
0.7 0.178579
0.7 0.180608
0.7 0.18236
0.7 0.186866
0.7 0.146603
0.7 0.161196
0.7 0.13808
0.7 0.148063
0.7 0.155779
0.7 0.150898
0.7 0.142261
0.7 0.175519
0.7 0.133252
0.7 0.147603
0.7 0.15414
0.7 0.143385
0.7 0.248712
0.7 0.146814
0.7 0.177581
0.7 0.15688
0.7 0.198864
0.7 0.142801
0.7 0.145029
0.7 0.152859
0.7 0.150314
0.7 0.141487
0.7 0.165299
0.7 0.150458
0.7 0.16541
0.7 0.177482
0.7 0.133631
0.7 0.143921
0.7 0.188555
0.7 0.133695
0.7 0.214295
0.7 0.162598
0.7 0.19855
0.7 0.166345
0.7 0.146557
0.7 0.206097
0.7 0.125642
0.7 0.191521
0.7 0.140143
0.7 0.133434
0.7 0.145104
0.7 0.138016
0.7 0.160858
0.7 0.166512
0.7 0.166304
0.7 0.143832
0.7 0.176218
0.7 0.169096
0.7 0.118484
0.7 0.158941
0.7 0.14762
0.7 0.154951
0.7 0.151826
0.7 0.17625
0.7 0.153487
0.7 0.163186
0.7 0.142873
0.7 0.129001
0.7 0.153767
};
\addplot [draw=orange, fill=orange, mark=asterisk, only marks]
table{%
x  y
0.8 0.831745
0.8 0.997629
0.8 1.04982
0.8 0.845004
0.8 0.780417
0.8 0.838524
0.8 0.701682
0.8 0.691703
0.8 0.886181
0.8 0.894868
0.8 0.801079
0.8 0.957774
0.8 0.791206
0.8 1.02642
0.8 0.784795
0.8 0.811101
0.8 0.933403
0.8 0.956456
0.8 0.807454
0.8 0.777095
0.8 0.799141
0.8 0.883791
0.8 0.80501
0.8 0.775204
0.8 0.749938
0.8 0.84614
0.8 0.772902
0.8 0.784396
0.8 0.980268
0.8 0.759617
0.8 0.832897
0.8 0.837917
0.8 0.83964
0.8 0.894826
0.8 0.916253
0.8 0.945846
0.8 0.759592
0.8 0.742759
0.8 1.00712
0.8 0.88529
0.8 0.814531
0.8 0.764979
0.8 0.900708
0.8 0.891934
0.8 0.812405
0.8 0.788213
0.8 0.791658
0.8 0.924115
0.8 0.781437
0.8 0.827403
0.8 0.81103
0.8 0.811781
0.8 0.780213
0.8 0.876358
0.8 0.778518
0.8 0.789781
0.8 0.814243
0.8 0.741943
0.8 0.823418
0.8 0.831838
0.8 0.921328
0.8 1.02044
0.8 0.887944
0.8 0.787302
0.8 0.968179
0.8 0.759241
0.8 0.792388
0.8 0.981307
0.8 0.936641
0.8 0.795378
0.8 0.928979
0.8 0.873232
0.8 0.824432
0.8 0.879511
0.8 0.667232
0.8 0.788421
0.8 0.752554
0.8 0.881028
0.8 0.851902
0.8 0.736984
0.8 0.756594
0.8 0.872828
0.8 1.00467
0.8 0.81011
0.8 0.933977
0.8 0.959623
0.8 0.73989
0.8 0.925413
0.8 0.852402
0.8 0.791026
0.8 0.772126
0.8 0.77433
0.8 0.893206
0.8 0.894201
0.8 0.833704
0.8 0.955707
0.8 0.945571
0.8 0.960062
0.8 0.834105
0.8 0.886112
};
\addplot [draw=black, fill=black, mark=diamond*, only marks]
table{%
x  y
0.9 1.31402
0.9 1.43743
0.9 1.12645
0.9 1.24544
0.9 1.28001
0.9 1.23376
0.9 1.39385
0.9 1.48043
0.9 1.41883
0.9 1.65154
0.9 1.30553
0.9 1.37737
0.9 1.25751
0.9 1.385
0.9 1.37615
0.9 1.20039
0.9 1.24747
0.9 1.69855
0.9 1.28922
0.9 1.56919
0.9 1.51543
0.9 1.29779
0.9 1.48585
0.9 1.09099
0.9 1.41589
0.9 1.24146
0.9 1.42202
0.9 1.03454
0.9 2.27417
0.9 1.27839
0.9 1.07117
0.9 1.15309
0.9 1.24555
0.9 1.34794
0.9 1.27967
0.9 1.43029
0.9 1.52169
0.9 1.32976
0.9 1.4148
0.9 1.30244
0.9 1.21546
0.9 1.49994
0.9 1.22374
0.9 1.17082
0.9 1.15035
0.9 1.45097
0.9 1.38179
0.9 1.26135
0.9 1.34684
0.9 1.34175
0.9 1.48827
0.9 1.49324
0.9 1.36779
0.9 1.26735
0.9 1.6265
0.9 1.75662
0.9 1.17914
0.9 1.31367
0.9 1.30112
0.9 1.20882
0.9 1.18847
0.9 1.32476
0.9 1.48771
0.9 1.05779
0.9 1.37028
0.9 1.13086
0.9 1.57602
0.9 1.16807
0.9 1.3163
0.9 1.35469
0.9 1.42982
0.9 1.15382
0.9 1.40905
0.9 1.27393
0.9 1.77689
0.9 1.36421
0.9 1.15924
0.9 1.16321
0.9 1.50016
0.9 1.27432
0.9 1.34333
0.9 1.30963
0.9 1.42774
0.9 1.08325
0.9 1.35006
0.9 1.22821
0.9 1.22899
0.9 1.21796
0.9 1.35536
0.9 1.3415
0.9 1.50756
0.9 1.20133
0.9 1.22856
0.9 1.52911
0.9 1.3354
0.9 1.13552
0.9 1.27977
0.9 1.71831
0.9 1.45814
0.9 1.11104
};
\end{axis}

\end{tikzpicture}}
   \caption{}
   \label{fig:MultipleEvent_rmsa}
 \end{subfigure}
 \begin{subfigure}{\textwidth}
    \centering
    \includegraphics[width=\columnwidth]{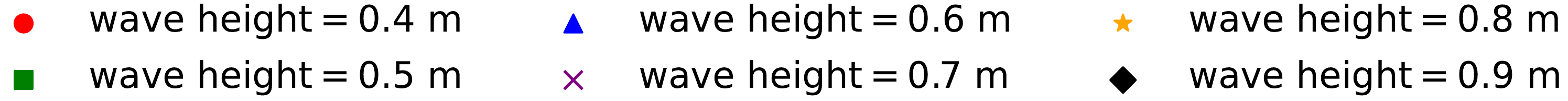}
    \label{fig:legend_RMSA} 
  \end{subfigure}
\caption{Two-storey structure's root mean square accelerations (RMSA) in $x$-direction obtained using OpenSees plotted vs. (a) W-section weight per length for columns (column weight), (b) W-section weight per length for beams (beam weight), (c) W-section weight per length for girder (girder weight), (d) Yield strength, (e) Young's modulus, (f) Wave height}
\label{fig:Root mean square accelerations_2stories}
\end{figure}

Considering the correlations shown between RMSA and the random variables in \autoref{fig:Root mean square accelerations_2stories}, the trend discussed for the Youngs' modulus and wave height with respect to RMSA is consistent with the observations with peak displacement.

However, unlike as observed with peak displacements, the variations in RMSA versus random variables show a clear distinct behaviour for each of the different wave heights considered. This further establishes that the initial wave heights significantly influence RMSA. The column, beam and girder stiffness, yield strength all show a weak nonlinear relationship with the RMSA. As expected, the overall RMSA decreases as Youngs' modulus value increases.

In reality, the loads exerted on a structure are influenced by different aspects, including the dynamics of the environment, the configuration of the structure, and the material composition. This variability introduces complexity and uncertainty into the load-bearing dynamics of the structure. The RMSA analysis shown in \autoref{fig:Root mean square accelerations_2stories} is based on the assumption of constant loads, whereas other distribution types reflect different structural responses that may occur. For example, static loads or uniformly distributed loads show characteristics of a normal distribution or uniform distribution, especially when the loads are symmetrically distributed. However, when considering irregularities in material properties or variations in foundation depths, modelling with other distributions, for example like beta distribution, is promising.

Furthermore realistic structural models should be explored to enable more accurate structural modelling and benchmarking. It is possible that the relationship between these structural parameters and RMSA can be strongly nonlinear for taller structures and warrants further investigation. 

\subsubsection{Influence of probability distribution}
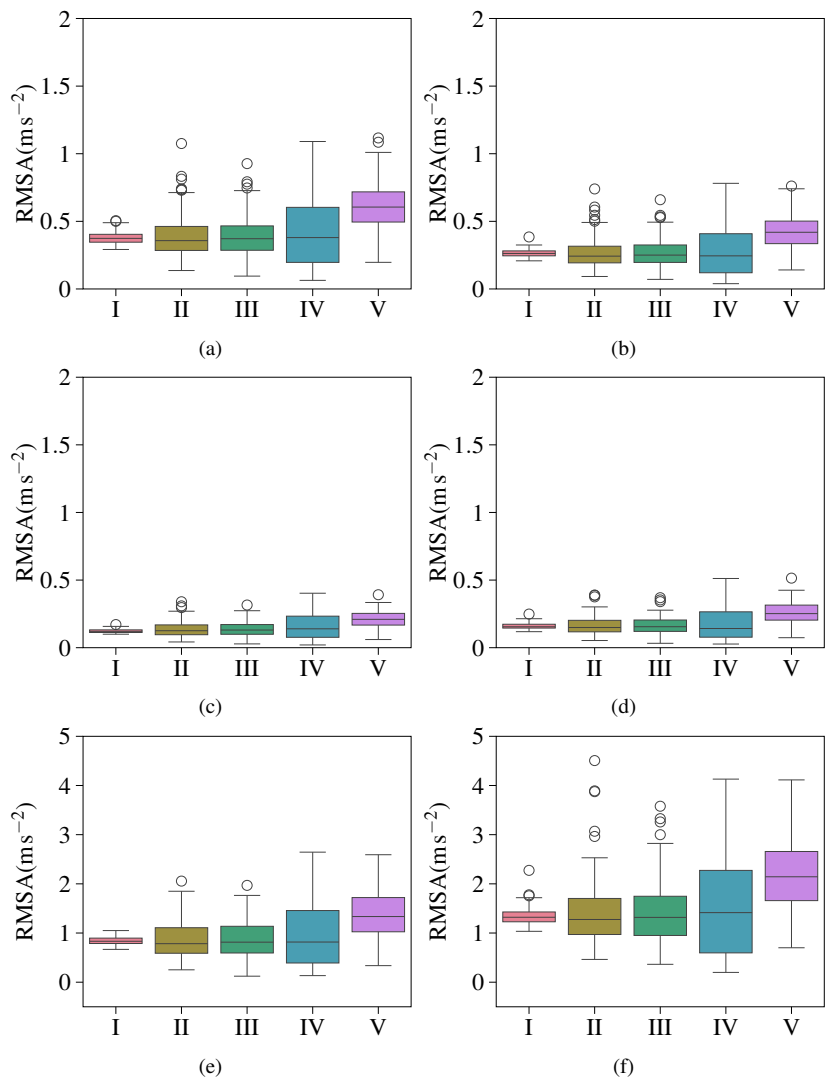
\begin{figure}[!htb]
 \centering
 \begin{subfigure}{0.44\textwidth}
   \centering
   \adjustbox{width=\linewidth}{
\begin{tikzpicture}

\definecolor{darkgray176}{RGB}{176,176,176}
\definecolor{darkslategray66}{RGB}{66,66,66}
\definecolor{mediumseagreen66160120}{RGB}{66,160,120}
\definecolor{palevioletred230129147}{RGB}{230,129,147}
\definecolor{peru15814564}{RGB}{158,145,64}
\definecolor{plum198136228}{RGB}{198,136,228}
\definecolor{steelblue73158179}{RGB}{73,158,179}

\begin{axis}[
tick align=outside,
tick pos=left,
x grid style={darkgray176},
xmin=-0.5, xmax=4.5,
xtick style={color=black},
xtick={0,1,2,3,4},
xticklabels={\uppercase\expandafter{\romannumeral 1}, \uppercase\expandafter{\romannumeral 2}, \uppercase\expandafter{\romannumeral 3}, \uppercase\expandafter{\romannumeral 4}, \uppercase\expandafter{\romannumeral 5}},
y grid style={darkgray176},
ylabel={RMSA(\unit{\meter\per\second\squared})},
ymin=0, ymax=2,
ytick style={color=black}
]
\path [draw=darkslategray66, fill=palevioletred230129147]
(axis cs:-0.4,0.345781)
--(axis cs:0.4,0.345781)
--(axis cs:0.4,0.40403525)
--(axis cs:-0.4,0.40403525)
--(axis cs:-0.4,0.345781)
--cycle;
\addplot [darkslategray66]
table {%
0 0.345781
0 0.291977
};
\addplot [darkslategray66]
table {%
0 0.40403525
0 0.490042
};
\addplot [darkslategray66]
table {%
-0.2 0.291977
0.2 0.291977
};
\addplot [darkslategray66]
table {%
-0.2 0.490042
0.2 0.490042
};
\addplot [black, mark=o, mark size=3, mark options={solid,fill opacity=0,draw=darkslategray66}, only marks]
table {%
0 0.500213
0 0.504986
};
\path [draw=darkslategray66, fill=peru15814564]
(axis cs:0.6,0.28474425)
--(axis cs:1.4,0.28474425)
--(axis cs:1.4,0.4622935)
--(axis cs:0.6,0.4622935)
--(axis cs:0.6,0.28474425)
--cycle;
\addplot [darkslategray66]
table {%
1 0.28474425
1 0.13628
};
\addplot [darkslategray66]
table {%
1 0.4622935
1 0.712695
};
\addplot [darkslategray66]
table {%
0.8 0.13628
1.2 0.13628
};
\addplot [darkslategray66]
table {%
0.8 0.712695
1.2 0.712695
};
\addplot [black, mark=o, mark size=3, mark options={solid,fill opacity=0,draw=darkslategray66}, only marks]
table {%
1 0.729986
1 0.735496
1 0.739907
1 1.07547
1 0.8328
1 0.81201
};
\path [draw=darkslategray66, fill=mediumseagreen66160120]
(axis cs:1.6,0.28608)
--(axis cs:2.4,0.28608)
--(axis cs:2.4,0.4662895)
--(axis cs:1.6,0.4662895)
--(axis cs:1.6,0.28608)
--cycle;
\addplot [darkslategray66]
table {%
2 0.28608
2 0.0950133
};
\addplot [darkslategray66]
table {%
2 0.4662895
2 0.727195
};
\addplot [darkslategray66]
table {%
1.8 0.0950133
2.2 0.0950133
};
\addplot [darkslategray66]
table {%
1.8 0.727195
2.2 0.727195
};
\addplot [black, mark=o, mark size=3, mark options={solid,fill opacity=0,draw=darkslategray66}, only marks]
table {%
2 0.746714
2 0.9269
2 0.774994
2 0.792616
};
\path [draw=darkslategray66, fill=steelblue73158179]
(axis cs:2.6,0.19653925)
--(axis cs:3.4,0.19653925)
--(axis cs:3.4,0.60343425)
--(axis cs:2.6,0.60343425)
--(axis cs:2.6,0.19653925)
--cycle;
\addplot [darkslategray66]
table {%
3 0.19653925
3 0.0635136
};
\addplot [darkslategray66]
table {%
3 0.60343425
3 1.09034
};
\addplot [darkslategray66]
table {%
2.8 0.0635136
3.2 0.0635136
};
\addplot [darkslategray66]
table {%
2.8 1.09034
3.2 1.09034
};
\path [draw=darkslategray66, fill=plum198136228]
(axis cs:3.6,0.494829)
--(axis cs:4.4,0.494829)
--(axis cs:4.4,0.717973)
--(axis cs:3.6,0.717973)
--(axis cs:3.6,0.494829)
--cycle;
\addplot [darkslategray66]
table {%
4 0.494829
4 0.196848
};
\addplot [darkslategray66]
table {%
4 0.717973
4 1.00978
};
\addplot [darkslategray66]
table {%
3.8 0.196848
4.2 0.196848
};
\addplot [darkslategray66]
table {%
3.8 1.00978
4.2 1.00978
};
\addplot [black, mark=o, mark size=3, mark options={solid,fill opacity=0,draw=darkslategray66}, only marks]
table {%
4 1.1167
4 1.08488
};
\addplot [darkslategray66]
table {%
-0.4 0.3740965
0.4 0.3740965
};
\addplot [darkslategray66]
table {%
0.6 0.3573115
1.4 0.3573115
};
\addplot [darkslategray66]
table {%
1.6 0.371619
2.4 0.371619
};
\addplot [darkslategray66]
table {%
2.6 0.3799015
3.4 0.3799015
};
\addplot [darkslategray66]
table {%
3.6 0.6055535
4.4 0.6055535
};
\end{axis}

\end{tikzpicture}}
   \caption{}
   \label{fig:wv04_multidistribution}
 \end{subfigure}
 \begin{subfigure}{0.44\textwidth}
   \centering
   \adjustbox{width=\linewidth}{
\begin{tikzpicture}

\definecolor{darkgray176}{RGB}{176,176,176}
\definecolor{darkslategray66}{RGB}{66,66,66}
\definecolor{mediumseagreen66160120}{RGB}{66,160,120}
\definecolor{palevioletred230129147}{RGB}{230,129,147}
\definecolor{peru15814564}{RGB}{158,145,64}
\definecolor{plum198136228}{RGB}{198,136,228}
\definecolor{steelblue73158179}{RGB}{73,158,179}

\begin{axis}[
tick align=outside,
tick pos=left,
x grid style={darkgray176},
xmin=-0.5, xmax=4.5,
xtick style={color=black},
xtick={0,1,2,3,4},
xticklabels={\uppercase\expandafter{\romannumeral 1}, \uppercase\expandafter{\romannumeral 2}, \uppercase\expandafter{\romannumeral 3}, \uppercase\expandafter{\romannumeral 4}, \uppercase\expandafter{\romannumeral 5}},
y grid style={darkgray176},
ylabel={RMSA(\unit{\meter\per\second\squared})},
ymin=0, ymax=2,
ytick style={color=black}
]
\path [draw=darkslategray66, fill=palevioletred230129147]
(axis cs:-0.4,0.245016)
--(axis cs:0.4,0.245016)
--(axis cs:0.4,0.28150725)
--(axis cs:-0.4,0.28150725)
--(axis cs:-0.4,0.245016)
--cycle;
\addplot [darkslategray66]
table {%
0 0.245016
0 0.208024
};
\addplot [darkslategray66]
table {%
0 0.28150725
0 0.325136
};
\addplot [darkslategray66]
table {%
-0.2 0.208024
0.2 0.208024
};
\addplot [darkslategray66]
table {%
-0.2 0.325136
0.2 0.325136
};
\addplot [black, mark=o, mark size=3, mark options={solid,fill opacity=0,draw=darkslategray66}, only marks]
table {%
0 0.384616
};
\path [draw=darkslategray66, fill=peru15814564]
(axis cs:0.6,0.19295625)
--(axis cs:1.4,0.19295625)
--(axis cs:1.4,0.316096)
--(axis cs:0.6,0.316096)
--(axis cs:0.6,0.19295625)
--cycle;
\addplot [darkslategray66]
table {%
1 0.19295625
1 0.0922883
};
\addplot [darkslategray66]
table {%
1 0.316096
1 0.49187
};
\addplot [darkslategray66]
table {%
0.8 0.0922883
1.2 0.0922883
};
\addplot [darkslategray66]
table {%
0.8 0.49187
1.2 0.49187
};
\addplot [black, mark=o, mark size=3, mark options={solid,fill opacity=0,draw=darkslategray66}, only marks]
table {%
1 0.739576
1 0.544739
1 0.58296
1 0.501141
1 0.513392
1 0.605523
};
\path [draw=darkslategray66, fill=mediumseagreen66160120]
(axis cs:1.6,0.19629875)
--(axis cs:2.4,0.19629875)
--(axis cs:2.4,0.32526525)
--(axis cs:1.6,0.32526525)
--(axis cs:1.6,0.19629875)
--cycle;
\addplot [darkslategray66]
table {%
2 0.19629875
2 0.0708574
};
\addplot [darkslategray66]
table {%
2 0.32526525
2 0.493901
};
\addplot [darkslategray66]
table {%
1.8 0.0708574
2.2 0.0708574
};
\addplot [darkslategray66]
table {%
1.8 0.493901
2.2 0.493901
};
\addplot [black, mark=o, mark size=3, mark options={solid,fill opacity=0,draw=darkslategray66}, only marks]
table {%
2 0.659863
2 0.527579
2 0.533193
2 0.5443
};
\path [draw=darkslategray66, fill=steelblue73158179]
(axis cs:2.6,0.11983)
--(axis cs:3.4,0.11983)
--(axis cs:3.4,0.40844725)
--(axis cs:2.6,0.40844725)
--(axis cs:2.6,0.11983)
--cycle;
\addplot [darkslategray66]
table {%
3 0.11983
3 0.0391474
};
\addplot [darkslategray66]
table {%
3 0.40844725
3 0.781043
};
\addplot [darkslategray66]
table {%
2.8 0.0391474
3.2 0.0391474
};
\addplot [darkslategray66]
table {%
2.8 0.781043
3.2 0.781043
};
\path [draw=darkslategray66, fill=plum198136228]
(axis cs:3.6,0.335061)
--(axis cs:4.4,0.335061)
--(axis cs:4.4,0.5018355)
--(axis cs:3.6,0.5018355)
--(axis cs:3.6,0.335061)
--cycle;
\addplot [darkslategray66]
table {%
4 0.335061
4 0.140588
};
\addplot [darkslategray66]
table {%
4 0.5018355
4 0.740759
};
\addplot [darkslategray66]
table {%
3.8 0.140588
4.2 0.140588
};
\addplot [darkslategray66]
table {%
3.8 0.740759
4.2 0.740759
};
\addplot [black, mark=o, mark size=3, mark options={solid,fill opacity=0,draw=darkslategray66}, only marks]
table {%
4 0.761013
};
\addplot [darkslategray66]
table {%
-0.4 0.262541
0.4 0.262541
};
\addplot [darkslategray66]
table {%
0.6 0.242656
1.4 0.242656
};
\addplot [darkslategray66]
table {%
1.6 0.2503475
2.4 0.2503475
};
\addplot [darkslategray66]
table {%
2.6 0.244796
3.4 0.244796
};
\addplot [darkslategray66]
table {%
3.6 0.4191925
4.4 0.4191925
};
\end{axis}

\end{tikzpicture}}
   \caption{}
   \label{fig:wv05_multidistribution}
 \end{subfigure}
 \begin{subfigure}{0.44\textwidth}
   \centering
   \adjustbox{width=\linewidth}{
\begin{tikzpicture}

\definecolor{darkgray176}{RGB}{176,176,176}
\definecolor{darkslategray66}{RGB}{66,66,66}
\definecolor{mediumseagreen66160120}{RGB}{66,160,120}
\definecolor{palevioletred230129147}{RGB}{230,129,147}
\definecolor{peru15814564}{RGB}{158,145,64}
\definecolor{plum198136228}{RGB}{198,136,228}
\definecolor{steelblue73158179}{RGB}{73,158,179}

\begin{axis}[
tick align=outside,
tick pos=left,
x grid style={darkgray176},
xmin=-0.5, xmax=4.5,
xtick style={color=black},
xtick={0,1,2,3,4},
xticklabels={\uppercase\expandafter{\romannumeral 1}, \uppercase\expandafter{\romannumeral 2}, \uppercase\expandafter{\romannumeral 3}, \uppercase\expandafter{\romannumeral 4}, \uppercase\expandafter{\romannumeral 5}},
y grid style={darkgray176},
ylabel={RMSA(\unit{\meter\per\second\squared})},
ymin=0, ymax=2,
ytick style={color=black}
]
\path [draw=darkslategray66, fill=palevioletred230129147]
(axis cs:-0.4,0.11242525)
--(axis cs:0.4,0.11242525)
--(axis cs:0.4,0.1318605)
--(axis cs:-0.4,0.1318605)
--(axis cs:-0.4,0.11242525)
--cycle;
\addplot [darkslategray66]
table {%
0 0.11242525
0 0.0997801
};
\addplot [darkslategray66]
table {%
0 0.1318605
0 0.157555
};
\addplot [darkslategray66]
table {%
-0.2 0.0997801
0.2 0.0997801
};
\addplot [darkslategray66]
table {%
-0.2 0.157555
0.2 0.157555
};
\addplot [black, mark=o, mark size=3, mark options={solid,fill opacity=0,draw=darkslategray66}, only marks]
table {%
0 0.171071
};
\path [draw=darkslategray66, fill=peru15814564]
(axis cs:0.6,0.09603235)
--(axis cs:1.4,0.09603235)
--(axis cs:1.4,0.168583)
--(axis cs:0.6,0.168583)
--(axis cs:0.6,0.09603235)
--cycle;
\addplot [darkslategray66]
table {%
1 0.09603235
1 0.0424623
};
\addplot [darkslategray66]
table {%
1 0.168583
1 0.269745
};
\addplot [darkslategray66]
table {%
0.8 0.0424623
1.2 0.0424623
};
\addplot [darkslategray66]
table {%
0.8 0.269745
1.2 0.269745
};
\addplot [black, mark=o, mark size=3, mark options={solid,fill opacity=0,draw=darkslategray66}, only marks]
table {%
1 0.339764
1 0.309257
1 0.296408
};
\path [draw=darkslategray66, fill=mediumseagreen66160120]
(axis cs:1.6,0.099187525)
--(axis cs:2.4,0.099187525)
--(axis cs:2.4,0.1711735)
--(axis cs:1.6,0.1711735)
--(axis cs:1.6,0.099187525)
--cycle;
\addplot [darkslategray66]
table {%
2 0.099187525
2 0.0281419
};
\addplot [darkslategray66]
table {%
2 0.1711735
2 0.27349
};
\addplot [darkslategray66]
table {%
1.8 0.0281419
2.2 0.0281419
};
\addplot [darkslategray66]
table {%
1.8 0.27349
2.2 0.27349
};
\addplot [black, mark=o, mark size=3, mark options={solid,fill opacity=0,draw=darkslategray66}, only marks]
table {%
2 0.315595
};
\path [draw=darkslategray66, fill=steelblue73158179]
(axis cs:2.6,0.07681295)
--(axis cs:3.4,0.07681295)
--(axis cs:3.4,0.23321275)
--(axis cs:2.6,0.23321275)
--(axis cs:2.6,0.07681295)
--cycle;
\addplot [darkslategray66]
table {%
3 0.07681295
3 0.0203207
};
\addplot [darkslategray66]
table {%
3 0.23321275
3 0.402808
};
\addplot [darkslategray66]
table {%
2.8 0.0203207
3.2 0.0203207
};
\addplot [darkslategray66]
table {%
2.8 0.402808
3.2 0.402808
};
\path [draw=darkslategray66, fill=plum198136228]
(axis cs:3.6,0.16722975)
--(axis cs:4.4,0.16722975)
--(axis cs:4.4,0.25373125)
--(axis cs:3.6,0.25373125)
--(axis cs:3.6,0.16722975)
--cycle;
\addplot [darkslategray66]
table {%
4 0.16722975
4 0.0603421
};
\addplot [darkslategray66]
table {%
4 0.25373125
4 0.334371
};
\addplot [darkslategray66]
table {%
3.8 0.0603421
4.2 0.0603421
};
\addplot [darkslategray66]
table {%
3.8 0.334371
4.2 0.334371
};
\addplot [black, mark=o, mark size=3, mark options={solid,fill opacity=0,draw=darkslategray66}, only marks]
table {%
4 0.391681
};
\addplot [darkslategray66]
table {%
-0.4 0.1208065
0.4 0.1208065
};
\addplot [darkslategray66]
table {%
0.6 0.125549
1.4 0.125549
};
\addplot [darkslategray66]
table {%
1.6 0.130552
2.4 0.130552
};
\addplot [darkslategray66]
table {%
2.6 0.1396955
3.4 0.1396955
};
\addplot [darkslategray66]
table {%
3.6 0.209742
4.4 0.209742
};
\end{axis}

\end{tikzpicture}}
   \caption{}
   \label{fig:wv06_multidistribution}
 \end{subfigure}
  \begin{subfigure}{0.44\textwidth}
   \centering
   \adjustbox{width=\linewidth}{
\begin{tikzpicture}

\definecolor{darkgray176}{RGB}{176,176,176}
\definecolor{darkslategray66}{RGB}{66,66,66}
\definecolor{mediumseagreen66160120}{RGB}{66,160,120}
\definecolor{palevioletred230129147}{RGB}{230,129,147}
\definecolor{peru15814564}{RGB}{158,145,64}
\definecolor{plum198136228}{RGB}{198,136,228}
\definecolor{steelblue73158179}{RGB}{73,158,179}

\begin{axis}[
tick align=outside,
tick pos=left,
x grid style={darkgray176},
xmin=-0.5, xmax=4.5,
xtick style={color=black},
xtick={0,1,2,3,4},
xticklabels={\uppercase\expandafter{\romannumeral 1}, \uppercase\expandafter{\romannumeral 2}, \uppercase\expandafter{\romannumeral 3}, \uppercase\expandafter{\romannumeral 4}, \uppercase\expandafter{\romannumeral 5}},
y grid style={darkgray176},
ylabel={RMSA(\unit{\meter\per\second\squared})},
ymin=0, ymax=2,
ytick style={color=black}
]
\path [draw=darkslategray66, fill=palevioletred230129147]
(axis cs:-0.4,0.14508525)
--(axis cs:0.4,0.14508525)
--(axis cs:0.4,0.173433)
--(axis cs:-0.4,0.173433)
--(axis cs:-0.4,0.14508525)
--cycle;
\addplot [darkslategray66]
table {%
0 0.14508525
0 0.118484
};
\addplot [darkslategray66]
table {%
0 0.173433
0 0.214295
};
\addplot [darkslategray66]
table {%
-0.2 0.118484
0.2 0.118484
};
\addplot [darkslategray66]
table {%
-0.2 0.214295
0.2 0.214295
};
\addplot [black, mark=o, mark size=3, mark options={solid,fill opacity=0,draw=darkslategray66}, only marks]
table {%
0 0.248712
};
\path [draw=darkslategray66, fill=peru15814564]
(axis cs:0.6,0.11731325)
--(axis cs:1.4,0.11731325)
--(axis cs:1.4,0.20244525)
--(axis cs:0.6,0.20244525)
--(axis cs:0.6,0.11731325)
--cycle;
\addplot [darkslategray66]
table {%
1 0.11731325
1 0.0530671
};
\addplot [darkslategray66]
table {%
1 0.20244525
1 0.301709
};
\addplot [darkslategray66]
table {%
0.8 0.0530671
1.2 0.0530671
};
\addplot [darkslategray66]
table {%
0.8 0.301709
1.2 0.301709
};
\addplot [black, mark=o, mark size=3, mark options={solid,fill opacity=0,draw=darkslategray66}, only marks]
table {%
1 0.375776
1 0.387005
1 0.38947
};
\path [draw=darkslategray66, fill=mediumseagreen66160120]
(axis cs:1.6,0.1205875)
--(axis cs:2.4,0.1205875)
--(axis cs:2.4,0.204787)
--(axis cs:1.6,0.204787)
--(axis cs:1.6,0.1205875)
--cycle;
\addplot [darkslategray66]
table {%
2 0.1205875
2 0.0325245
};
\addplot [darkslategray66]
table {%
2 0.204787
2 0.277449
};
\addplot [darkslategray66]
table {%
1.8 0.0325245
2.2 0.0325245
};
\addplot [darkslategray66]
table {%
1.8 0.277449
2.2 0.277449
};
\addplot [black, mark=o, mark size=3, mark options={solid,fill opacity=0,draw=darkslategray66}, only marks]
table {%
2 0.369672
2 0.353373
2 0.340484
};
\path [draw=darkslategray66, fill=steelblue73158179]
(axis cs:2.6,0.077620525)
--(axis cs:3.4,0.077620525)
--(axis cs:3.4,0.26528825)
--(axis cs:2.6,0.26528825)
--(axis cs:2.6,0.077620525)
--cycle;
\addplot [darkslategray66]
table {%
3 0.077620525
3 0.027104
};
\addplot [darkslategray66]
table {%
3 0.26528825
3 0.511811
};
\addplot [darkslategray66]
table {%
2.8 0.027104
3.2 0.027104
};
\addplot [darkslategray66]
table {%
2.8 0.511811
3.2 0.511811
};
\path [draw=darkslategray66, fill=plum198136228]
(axis cs:3.6,0.2039345)
--(axis cs:4.4,0.2039345)
--(axis cs:4.4,0.31495625)
--(axis cs:3.6,0.31495625)
--(axis cs:3.6,0.2039345)
--cycle;
\addplot [darkslategray66]
table {%
4 0.2039345
4 0.0739694
};
\addplot [darkslategray66]
table {%
4 0.31495625
4 0.424857
};
\addplot [darkslategray66]
table {%
3.8 0.0739694
4.2 0.0739694
};
\addplot [darkslategray66]
table {%
3.8 0.424857
4.2 0.424857
};
\addplot [black, mark=o, mark size=3, mark options={solid,fill opacity=0,draw=darkslategray66}, only marks]
table {%
4 0.514382
};
\addplot [darkslategray66]
table {%
-0.4 0.154834
0.4 0.154834
};
\addplot [darkslategray66]
table {%
0.6 0.149177
1.4 0.149177
};
\addplot [darkslategray66]
table {%
1.6 0.1551
2.4 0.1551
};
\addplot [darkslategray66]
table {%
2.6 0.1419525
3.4 0.1419525
};
\addplot [darkslategray66]
table {%
3.6 0.252012
4.4 0.252012
};
\end{axis}

\end{tikzpicture}}
   \caption{}
   \label{fig:wv07_multidistribution}
 \end{subfigure}
 \begin{subfigure}{0.44\textwidth}
   \centering
   \adjustbox{width=\linewidth}{
\begin{tikzpicture}

\definecolor{darkgray176}{RGB}{176,176,176}
\definecolor{darkslategray66}{RGB}{66,66,66}
\definecolor{mediumseagreen66160120}{RGB}{66,160,120}
\definecolor{palevioletred230129147}{RGB}{230,129,147}
\definecolor{peru15814564}{RGB}{158,145,64}
\definecolor{plum198136228}{RGB}{198,136,228}
\definecolor{steelblue73158179}{RGB}{73,158,179}

\begin{axis}[
tick align=outside,
tick pos=left,
x grid style={darkgray176},
xmin=-0.5, xmax=4.5,
xtick style={color=black},
xtick={0,1,2,3,4},
xticklabels={\uppercase\expandafter{\romannumeral 1}, \uppercase\expandafter{\romannumeral 2}, \uppercase\expandafter{\romannumeral 3}, \uppercase\expandafter{\romannumeral 4}, \uppercase\expandafter{\romannumeral 5}},
y grid style={darkgray176},
ylabel={RMSA(\unit{\meter\per\second\squared})},
ymin=-0.5, ymax=5,
ytick style={color=black}
]
\path [draw=darkslategray66, fill=palevioletred230129147]
(axis cs:-0.4,0.78667525)
--(axis cs:0.4,0.78667525)
--(axis cs:0.4,0.896328)
--(axis cs:-0.4,0.896328)
--(axis cs:-0.4,0.78667525)
--cycle;
\addplot [darkslategray66]
table {%
0 0.78667525
0 0.667232
};
\addplot [darkslategray66]
table {%
0 0.896328
0 1.04982
};
\addplot [darkslategray66]
table {%
-0.2 0.667232
0.2 0.667232
};
\addplot [darkslategray66]
table {%
-0.2 1.04982
0.2 1.04982
};
\path [draw=darkslategray66, fill=peru15814564]
(axis cs:0.6,0.58843425)
--(axis cs:1.4,0.58843425)
--(axis cs:1.4,1.10787)
--(axis cs:0.6,1.10787)
--(axis cs:0.6,0.58843425)
--cycle;
\addplot [darkslategray66]
table {%
1 0.58843425
1 0.251771
};
\addplot [darkslategray66]
table {%
1 1.10787
1 1.84931
};
\addplot [darkslategray66]
table {%
0.8 0.251771
1.2 0.251771
};
\addplot [darkslategray66]
table {%
0.8 1.84931
1.2 1.84931
};
\addplot [black, mark=o, mark size=3, mark options={solid,fill opacity=0,draw=darkslategray66}, only marks]
table {%
1 2.05595
};
\path [draw=darkslategray66, fill=mediumseagreen66160120]
(axis cs:1.6,0.59374075)
--(axis cs:2.4,0.59374075)
--(axis cs:2.4,1.136995)
--(axis cs:1.6,1.136995)
--(axis cs:1.6,0.59374075)
--cycle;
\addplot [darkslategray66]
table {%
2 0.59374075
2 0.122258
};
\addplot [darkslategray66]
table {%
2 1.136995
2 1.7648
};
\addplot [darkslategray66]
table {%
1.8 0.122258
2.2 0.122258
};
\addplot [darkslategray66]
table {%
1.8 1.7648
2.2 1.7648
};
\addplot [black, mark=o, mark size=3, mark options={solid,fill opacity=0,draw=darkslategray66}, only marks]
table {%
2 1.96885
};
\path [draw=darkslategray66, fill=steelblue73158179]
(axis cs:2.6,0.3903825)
--(axis cs:3.4,0.3903825)
--(axis cs:3.4,1.4570575)
--(axis cs:2.6,1.4570575)
--(axis cs:2.6,0.3903825)
--cycle;
\addplot [darkslategray66]
table {%
3 0.3903825
3 0.1338
};
\addplot [darkslategray66]
table {%
3 1.4570575
3 2.64451
};
\addplot [darkslategray66]
table {%
2.8 0.1338
3.2 0.1338
};
\addplot [darkslategray66]
table {%
2.8 2.64451
3.2 2.64451
};
\path [draw=darkslategray66, fill=plum198136228]
(axis cs:3.6,1.02505)
--(axis cs:4.4,1.02505)
--(axis cs:4.4,1.7210325)
--(axis cs:3.6,1.7210325)
--(axis cs:3.6,1.02505)
--cycle;
\addplot [darkslategray66]
table {%
4 1.02505
4 0.336748
};
\addplot [darkslategray66]
table {%
4 1.7210325
4 2.59211
};
\addplot [darkslategray66]
table {%
3.8 0.336748
4.2 0.336748
};
\addplot [darkslategray66]
table {%
3.8 2.59211
4.2 2.59211
};
\addplot [darkslategray66]
table {%
-0.4 0.8323675
0.4 0.8323675
};
\addplot [darkslategray66]
table {%
0.6 0.7834655
1.4 0.7834655
};
\addplot [darkslategray66]
table {%
1.6 0.814146
2.4 0.814146
};
\addplot [darkslategray66]
table {%
2.6 0.8167615
3.4 0.8167615
};
\addplot [darkslategray66]
table {%
3.6 1.33508
4.4 1.33508
};
\end{axis}

\end{tikzpicture}}
   \caption{}
   \label{fig:wv08_multidistribution}
 \end{subfigure}
 \begin{subfigure}{0.44\textwidth}
   \centering
   \adjustbox{width=\linewidth}{
\begin{tikzpicture}

\definecolor{darkgray176}{RGB}{176,176,176}
\definecolor{darkslategray66}{RGB}{66,66,66}
\definecolor{mediumseagreen66160120}{RGB}{66,160,120}
\definecolor{palevioletred230129147}{RGB}{230,129,147}
\definecolor{peru15814564}{RGB}{158,145,64}
\definecolor{plum198136228}{RGB}{198,136,228}
\definecolor{steelblue73158179}{RGB}{73,158,179}

\begin{axis}[
tick align=outside,
tick pos=left,
x grid style={darkgray176},
xmin=-0.5, xmax=4.5,
xtick style={color=black},
xtick={0,1,2,3,4},
xticklabels={\uppercase\expandafter{\romannumeral 1}, \uppercase\expandafter{\romannumeral 2}, \uppercase\expandafter{\romannumeral 3}, \uppercase\expandafter{\romannumeral 4}, \uppercase\expandafter{\romannumeral 5}},
y grid style={darkgray176},
ylabel={RMSA(\unit{\meter\per\second\squared})},
ymin=-0.5, ymax=5,
ytick style={color=black}
]
\path [draw=darkslategray66, fill=palevioletred230129147]
(axis cs:-0.4,1.2284725)
--(axis cs:0.4,1.2284725)
--(axis cs:0.4,1.42826)
--(axis cs:-0.4,1.42826)
--(axis cs:-0.4,1.2284725)
--cycle;
\addplot [darkslategray66]
table {%
0 1.2284725
0 1.03454
};
\addplot [darkslategray66]
table {%
0 1.42826
0 1.71831
};
\addplot [darkslategray66]
table {%
-0.2 1.03454
0.2 1.03454
};
\addplot [darkslategray66]
table {%
-0.2 1.71831
0.2 1.71831
};
\addplot [black, mark=o, mark size=3, mark options={solid,fill opacity=0,draw=darkslategray66}, only marks]
table {%
0 2.27417
0 1.75662
0 1.77689
};
\path [draw=darkslategray66, fill=peru15814564]
(axis cs:0.6,0.96974375)
--(axis cs:1.4,0.96974375)
--(axis cs:1.4,1.7030025)
--(axis cs:0.6,1.7030025)
--(axis cs:0.6,0.96974375)
--cycle;
\addplot [darkslategray66]
table {%
1 0.96974375
1 0.4631
};
\addplot [darkslategray66]
table {%
1 1.7030025
1 2.52992
};
\addplot [darkslategray66]
table {%
0.8 0.4631
1.2 0.4631
};
\addplot [darkslategray66]
table {%
0.8 2.52992
1.2 2.52992
};
\addplot [black, mark=o, mark size=3, mark options={solid,fill opacity=0,draw=darkslategray66}, only marks]
table {%
1 3.89035
1 3.0697
1 4.50665
1 2.96363
1 3.87447
};
\path [draw=darkslategray66, fill=mediumseagreen66160120]
(axis cs:1.6,0.95023425)
--(axis cs:2.4,0.95023425)
--(axis cs:2.4,1.7475625)
--(axis cs:1.6,1.7475625)
--(axis cs:1.6,0.95023425)
--cycle;
\addplot [darkslategray66]
table {%
2 0.95023425
2 0.36362
};
\addplot [darkslategray66]
table {%
2 1.7475625
2 2.82209
};
\addplot [darkslategray66]
table {%
1.8 0.36362
2.2 0.36362
};
\addplot [darkslategray66]
table {%
1.8 2.82209
2.2 2.82209
};
\addplot [black, mark=o, mark size=3, mark options={solid,fill opacity=0,draw=darkslategray66}, only marks]
table {%
2 3.25855
2 2.9993
2 3.57823
2 3.32544
};
\path [draw=darkslategray66, fill=steelblue73158179]
(axis cs:2.6,0.596116)
--(axis cs:3.4,0.596116)
--(axis cs:3.4,2.2735725)
--(axis cs:2.6,2.2735725)
--(axis cs:2.6,0.596116)
--cycle;
\addplot [darkslategray66]
table {%
3 0.596116
3 0.199025
};
\addplot [darkslategray66]
table {%
3 2.2735725
3 4.1298
};
\addplot [darkslategray66]
table {%
2.8 0.199025
3.2 0.199025
};
\addplot [darkslategray66]
table {%
2.8 4.1298
3.2 4.1298
};
\path [draw=darkslategray66, fill=plum198136228]
(axis cs:3.6,1.6588025)
--(axis cs:4.4,1.6588025)
--(axis cs:4.4,2.657895)
--(axis cs:3.6,2.657895)
--(axis cs:3.6,1.6588025)
--cycle;
\addplot [darkslategray66]
table {%
4 1.6588025
4 0.700246
};
\addplot [darkslategray66]
table {%
4 2.657895
4 4.11369
};
\addplot [darkslategray66]
table {%
3.8 0.700246
4.2 0.700246
};
\addplot [darkslategray66]
table {%
3.8 4.11369
4.2 4.11369
};
\addplot [darkslategray66]
table {%
-0.4 1.32053
0.4 1.32053
};
\addplot [darkslategray66]
table {%
0.6 1.2752
1.4 1.2752
};
\addplot [darkslategray66]
table {%
1.6 1.31798
2.4 1.31798
};
\addplot [darkslategray66]
table {%
2.6 1.41471
3.4 1.41471
};
\addplot [darkslategray66]
table {%
3.6 2.143615
4.4 2.143615
};
\end{axis}

\end{tikzpicture}}
   \caption{}
   \label{fig:wv09_multidistribution}
 \end{subfigure}
 \caption{RMSA with different distribution including: (\uppercase\expandafter{\romannumeral 1}) constant, (\uppercase\expandafter{\romannumeral 2}) lognormal, (\uppercase\expandafter{\romannumeral 3}) normal, (\uppercase\expandafter{\romannumeral 4}) uniform, (\uppercase\expandafter{\romannumeral 5}) beta under varying wave heights: (a) 0.4 \unit{\meter}, (b) 0.5 \unit{\meter}, (c) 0.6 \unit{\meter}, (d) 0.7 \unit{\meter}, (e) 0.8 \unit{\meter}, (f) 0.9 \unit{\meter}.}
 \label{fig:RMSA with multidistribution}
\end{figure}
This study further analyses the sensitivity of the choice of probability distribution. In this regard, lognormal (mean=1, standard deviation=0.2), normal (mean=1, standard deviation=0.2), uniform (min=0.4, max=1.6), and beta (alpha=5, beta=2, min=0.4, max=1.6) distributions are considered. The probability distribution parameters are considered such that the load values generated fall within the same overall range. This setup allows the distribution shape and characteristics to be the primary variables, enabling a clearer observation of the influence of distribution shapes and characteristics on structural responses. The diverse impacts of different distributions on structural responses are illustrated in \autoref{fig:RMSA with multidistribution}. 

The boxplot shown in \autoref{fig:RMSA with multidistribution} helps statistically visualise the variation of RMSA under different distributions and initial wave heights. It could be observed that the median lines in each distribution's boxplot are quite close, except for the slightly higher median line in the beta distribution. This indicates that the system exhibits strong robustness to moderate random loads (inferred from initial wave heights) under different distributions, as the RMSA response to moderate shifts in magnitude remains consistent.

Further, the interquartile range (IQR) is used to represent the middle $50\%$ of the RMSA data's distribution range, where the range between the upper quartile and lower quartile is the IQR. The upper whisker and lower whisker respectively represent the major distribution within the range of maximum and minimum values.

Under conditions of different wave heights with a uniform distribution applied to the loads in each case, the length from the upper whisker to the upper quartile is significantly longer compared to the length from the lower whisker to the lower quartile for each wave height. This phenomenon may be attributed to the asymmetry in structural response when the uniformly distributed load is applied as a random variable. Such asymmetry becomes particularly evident under higher load conditions, resulting in increased fluctuations in RMSA responses. Moreover, the application of uniformly distributed forces could induce extreme values in structural responses, causing the RMSA to predominantly cluster near the upper whisker.

Further, the circles in \autoref{fig:RMSA with multidistribution} mostly appear in when normal and lognormal distributions are considered. These circles represent outliers, indicating extreme RMSA values that exceed normal operational conditions. These extreme RMSA values need to be analysed to determine whether the structure is critically stressed to the point of damage and/or failure.

\begin{figure}[!htb]
 \centering
 \begin{subfigure}{0.47\textwidth}
   \centering
   \includegraphics[width=\columnwidth]{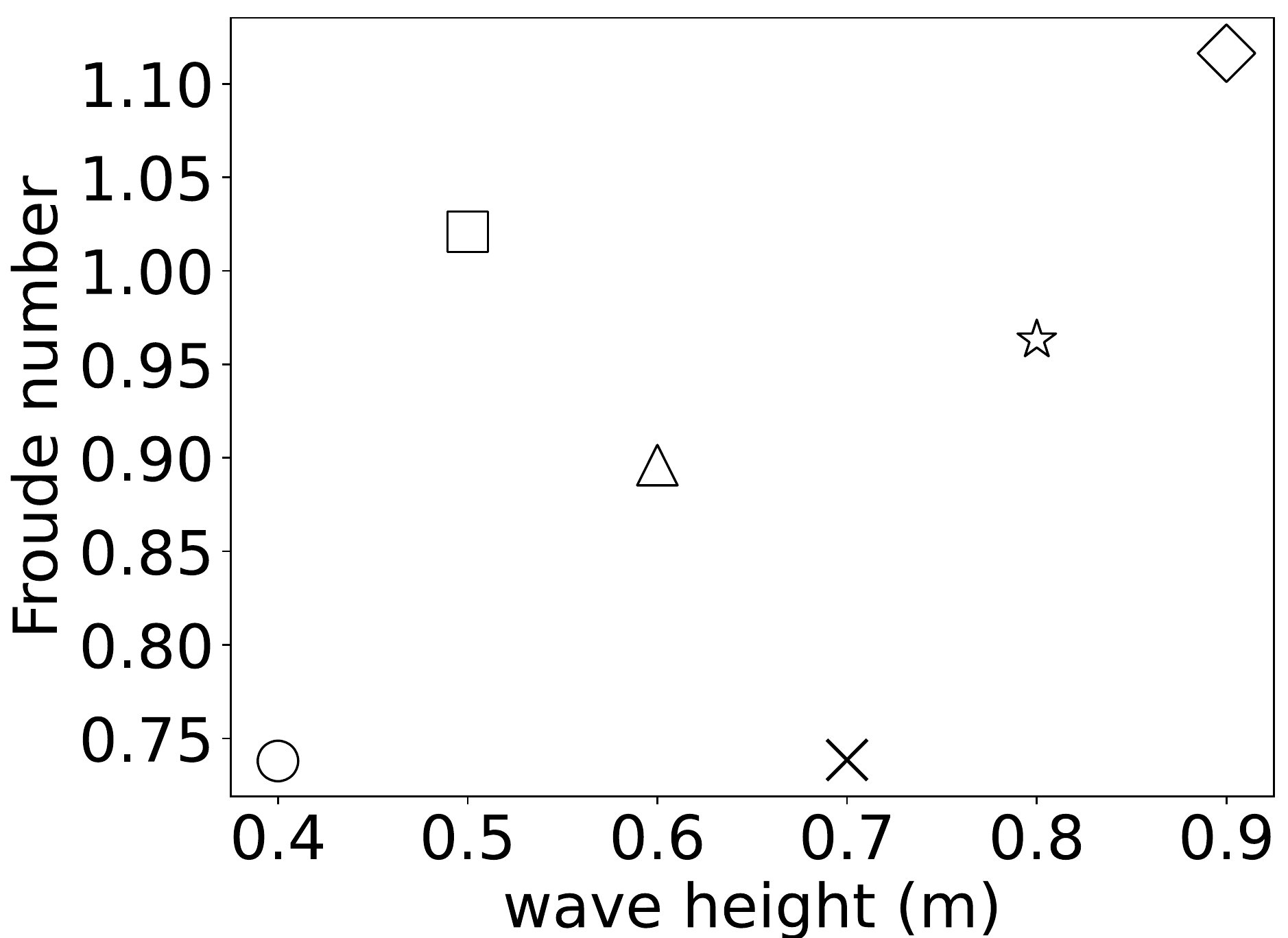}
   \caption{}
   \label{fig:Froud_max}
 \end{subfigure}
 \begin{subfigure}{0.44\textwidth}
   \centering
   \includegraphics[width=\linewidth]{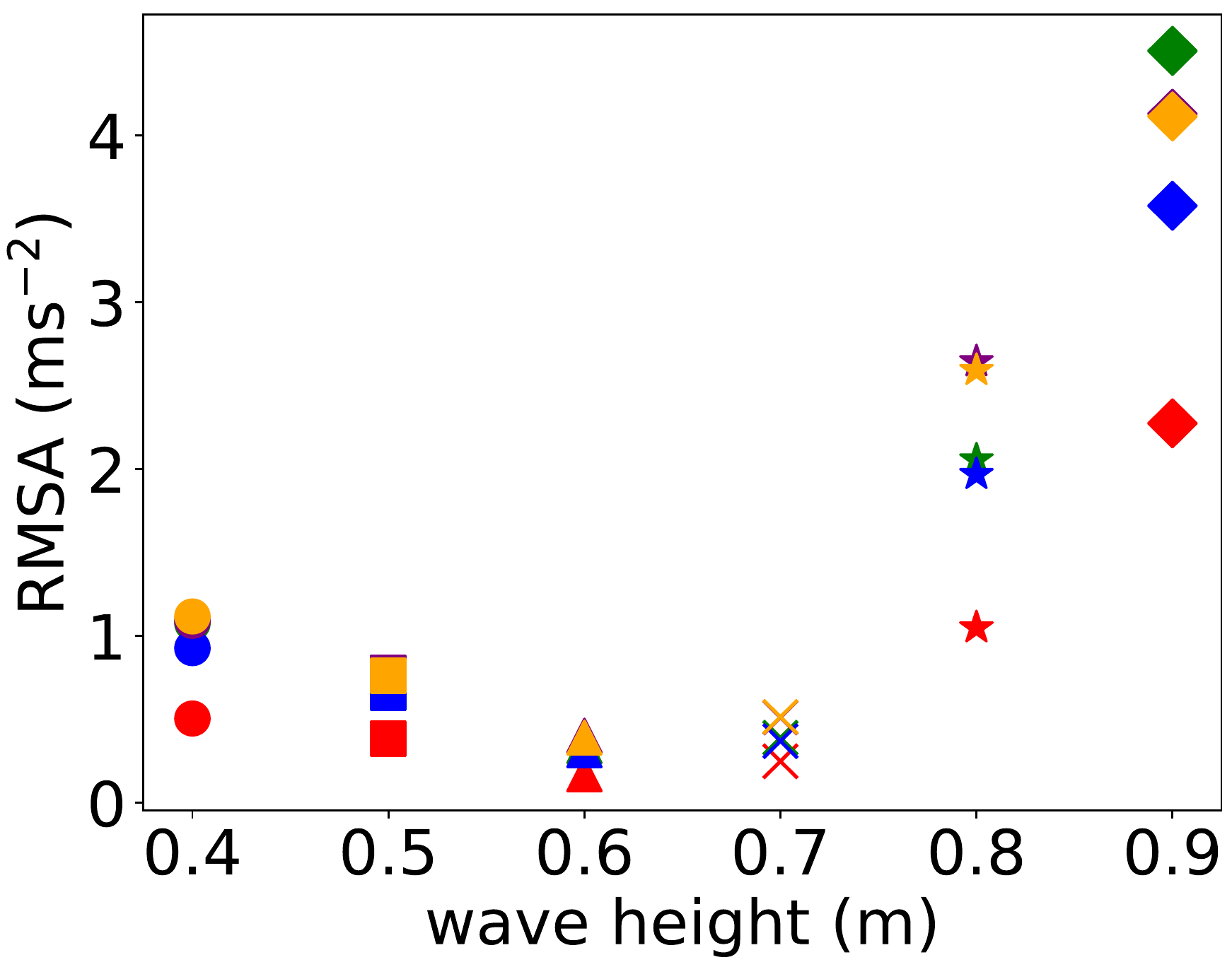}
   \caption{}
   \label{fig:multidistribution_max}
 \end{subfigure}
  \begin{subfigure}{\textwidth}
   \centering
   \includegraphics[width=\columnwidth]{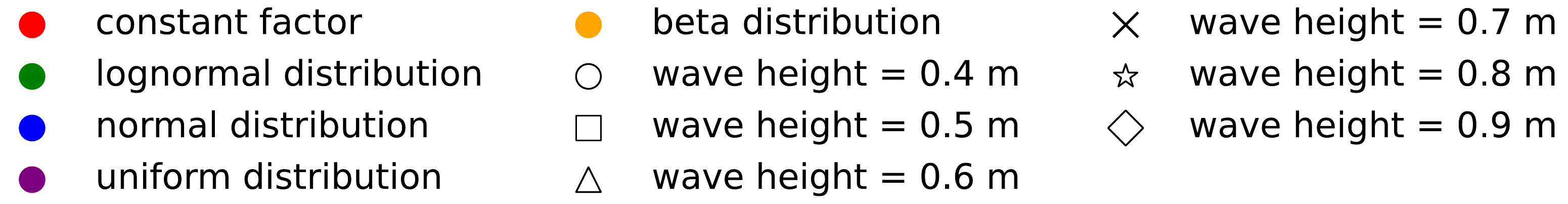}
 \end{subfigure}
\caption{(a) Maximum Froude number under varying wave heights. (b) RMSA with different distribution including: constant, lognormal, normal, uniform, beta under varying wave heights: 0.40 - 0.5 \unit{\meter} at intervals of 0.10 \unit{\meter}.}
\label{fig: Froud_max&RMSA_max}
\end{figure}

Thus, analysing the maximum values of RMSA produced by considering different wave heights and distributions could help to assess the maximum impact force that the building may experience under various scenarios and shown in \autoref{fig:multidistribution_max}. The constant factor results in the smallest maximum RMSA values among all distribution conditions.

When the wave height is 0.80 \unit{\meter}, the maximum value of RMSA generated by beta distribution applying to the wave load is $2.59211$ \unit{\meter\per\second\squared} which is almost 2.5 times of the RMSA generated by the constant factor. 

When the wave height is 0.90 \unit{\meter}, the maximum value of RMSA caused by lognormal distribution on wave load is $4.50665$ \unit{\meter\per\second\squared}, which is also the maximum value of RMSA under six different wave heights with different distribution conditions.

Overall, when the wave heights are 0.40 \unit{\meter}, 0.50 \unit{\meter}, 0.60 \unit{\meter}, 0.70 \unit{\meter}, and 0.80 \unit{\meter}, the maximum RMSA values generated by wave loads on the structure under Beta distribution and uniform distribution show very little discrepancy. These values are higher compared to the maximum RMSA values under normal distribution, lognormal distribution, and constant factor applied to the wave loads.

\subsubsection{Correlation with Froude number}
The Froude number is introduced in the study, which could indicate the correlation between the flow inertia force and the external force. Additionally, different Froude numbers indicate varying wave characteristics, which affect the magnitude and application of wave loads, thereby influencing the RMSA values of the structure.

The $V_{eff}$ from \autoref{eq:fdyn} was used in calculation of Froude number. The $V_{eff}$, here is defined as the average speed within the height range along the structure's centre line for each time step. Then, the maximum $V_{eff}$ value of different wave height are selected to calculate the Froude number, maximum Froude numbers for wave height ranging from 0.40 - 0.90 \unit{\meter}, in intervals of 0.10 \unit{\meter}, are plotted in \autoref{fig:Froud_max}. 

In \autoref{fig:Froud_max}, when wave heights are 0.50 \unit{\meter} and 0.90 \unit{\meter}, the maximum Froude number exceeds 1, indicating supercritical flow conditions for these cases. In all other cases, the maximum Froude number is less than 1, indicating sub-critical flow. It is noteworthy that at the wave height of 0.90 \unit{\meter}, where the maximum Froude number is approximately 1.116, the highest among all conditions depicted in \autoref{fig:Froud_max}. Therefore, the Froude number indicates that wave speeds are highest in this case, resulting in the maximum wave loads and RMSA values. 

Additionally, it is evident that the variation trend of the Froude number with wave height aligns consistently with the trend of force variation with wave height in \autoref{fig:forcesph}, and corresponds to the trend of RMSA variation with wave height in \autoref{fig:wv09_multidistribution}, except when the wave height is 0.40 \unit{\meter}. This discrepancy may be attributed to the lower wave height, resulting in a more stable flow regime.

\section{Conclusions and future work}
The work provides a comprehensive numerical framework to estimate probabilistic structural response under wave loading conditions in a wave flume. The HyTOFU flume configuration is used in this regard as the standard wave flume configuration. However, the developed method is general and applicable to other geometries as well. As shown, the flow conditions and the structural response are reasonably diverse even for a single wave flume geometry and thus, uncertainty quantification remains paramount to expand the interpretation of wave flume results into real-life conditions.  

While SPH for wave loading on a single structure nor FEM themselves are not novel, the novelty of this work lies in the coupling of different numerical methods, albeit one-way. The work demonstrates a nonlinear dependence on the initial wave height vs. the root mean square acceleration, as shown in \autoref{fig:MultipleEvent_rmsa}. This is further substantiated by the wave contours shown in Appendix in \autoref{fig:waveheight4} - \autoref{fig:waveheight9} where the wave breaking effects for medium-sized waves, of height comparable to still water depth, are shown. 

Overall, the work provides a systematic method for uncertainty quantification and data curation for wave-structure interaction problems. The envisioned next steps of this work is in the usage of the developed models for surrogate modelling towards a true digital wave flume. Despite the development of digital wave flumes using CFD models, full 3-D high-fidelity CFD simulations remain computationally expensive. Consequently, digital wave flumes have remained simulation constructs rather than true digital twins. This necessitates the creation of surrogate models capable of providing a sufficiently fast and accurate representation of the system. Furthermore, a range of uncertainties significantly influence  WSI problems, including initial wave conditions, air-water interaction, temperature, bathymetric surface roughness, and structural properties, among others. The computational cost of simulation techniques hinders extensive quantification of these uncertainties within the system. Conversely, the existence of a surrogate model could enable rapid evaluation of potential WSI scenarios. Considering the expense and time associated with experimental setup, such a surrogate model would not only facilitate the incorporation of uncertainties but also lead to more cost-effective and time-efficient experiment design. Further, considering extensive experimental data available with large wave flumes like those in Kyoto (Japan), Hannover (Germany), Oregon (USA) etc., this provides an excellent opportunity to develop digital twins of use to the community.

\section*{Conflict of interest statement}
On behalf of all authors, the corresponding author states that there is no conflict of interest. 

\section*{Acknowledgements}
The authors would like to acknowledge multiple sources of funding that has allowed this work. (a) Isaac Newton Institute for Mathematical Sciences for support and hospitality during the programme Uncertainty Quantification and Stochastic Modelling of Materials when work on this paper was undertaken. This work was supported by EPSRC Grant Number EP/R014604/1; (b) The Engineering and Physical Sciences Research Council (EPSRC) Impact Acceleration Account via The University of Manchester [IAA 391] (c) This project made use of time on ARCHER2 granted via the UK High-End Computing Consortium for Wave Structure Interaction (HEC-WSI) (\url{http://hec-wsi.ac.uk}), supported by EPSRC (grant no. EP/X035751/1).

\section*{CRediT authorship contribution statement}
\textbf{Xiaoyuan Luo:} Software, validation, investigation, data curation, critical analysis, writing - original draft, simulation; \textbf{Vijay Nandurdikar:} Software, simulation, validation, formal analysis, writing - original draft;  \textbf{Sangri-Yi:} Resources, supervision, , writing - editing and review; \textbf{Alistair Revell:} Resources, supervision, writing - editing and review; \textbf{Georgios Fourtakas:} Resources, supervision, writing - editing and review; \textbf{Ajay B Harish:} Conceptualisation, methodology, resources, writing - original draft, supervision, project administration, funding acquisition.

\section*{Data availability statement}
The input files, resulting data and plotting routines with documentation (presented in this work) will be shared through Zenodo and link to be made available in the final accepted manuscript.

\bibliographystyle{ClassFiles/spbasic}
\bibliography{Extras/LiteratureDatabase}

\begin{thebibliography}{54}
\providecommand{\natexlab}[1]{#1}
\providecommand{\url}[1]{{#1}}
\providecommand{\urlprefix}{URL }
\expandafter\ifx\csname urlstyle\endcsname\relax
  \providecommand{\doi}[1]{DOI~\discretionary{}{}{}#1}\else
  \providecommand{\doi}{DOI~\discretionary{}{}{}\begingroup
  \urlstyle{rm}\Url}\fi
\providecommand{\eprint}[2][]{\url{#2}}

\bibitem[{Altomare et~al(2014)Altomare, Crespo, Rogers, Dominguez, Gironella,
  and G{\'o}mez-Gesteira}]{ALTOMARE201434}
Altomare C, Crespo AJC, Rogers BD, Dominguez JM, Gironella X,
  G{\'o}mez-Gesteira M (2014) Numerical modelling of armour block sea
  breakwater with smoothed particle hydrodynamics. Comput Struct 130:34--45

\bibitem[{Altomare et~al(2015)Altomare, Crespo, Dom{\'i}nguez,
  G{\'o}mez-Gesteira, Suzuki, and Verwaest}]{ALTOMARE20151}
Altomare C, Crespo AJC, Dom{\'i}nguez JM, G{\'o}mez-Gesteira M, Suzuki T,
  Verwaest T (2015) Applicability of smoothed particle hydrodynamics for
  estimation of sea wave impact on coastal structures. Coast Eng 96:1--12

\bibitem[{Altomare et~al(2017)Altomare, Dom{\'i}nguez, Crespo,
  Gonz{\'a}lez-Cao, Suzuki, G{\'o}mez-Gesteira, and Troch}]{ALTOMARE201737}
Altomare C, Dom{\'i}nguez JM, Crespo AJC, Gonz{\'a}lez-Cao J, Suzuki T,
  G{\'o}mez-Gesteira M, Troch P (2017) Long-crested wave generation and
  absorption for {SPH}-based dualsphysics model. Coast Eng 127:37--54

\bibitem[{ASCE(2017)}]{american2017minimum}
ASCE (2017) Minimum design loads and associated criteria for buildings and
  other structures,. ASCE, \doi{10.1061/9780784414248}

\bibitem[{ASCE/SEI 7-16(2017)}]{asce}
ASCE/SEI 7-16 (2017) Minimum design loads and associated criteria for buildings
  and other structures

\bibitem[{Berger et~al(2011)Berger, George, LeVeque, and
  Mandli}]{BergerGeorgeLeVequeMandli11}
Berger MJ, George DL, LeVeque RJ, Mandli KT (2011) The {GeoClaw} software for
  depth-averaged flows with adaptive refinement. Adv in Water Res 34:1195--1206

\bibitem[{Chen et~al(2014)Chen, Zang, Hillis, Morgan, and
  Plummer}]{Chen2014-rz}
Chen LF, Zang J, Hillis AJ, Morgan GCJ, Plummer AR (2014) Numerical
  investigation of wave-structure interaction using {OpenFOAM}. Ocean Eng
  88:91--109

\bibitem[{Crespo et~al(2007)Crespo, G{\'o}mez-Gesteira, and
  Dalrymple}]{cmc2007}
Crespo A, G{\'o}mez-Gesteira M, Dalrymple R (2007) Boundary conditions
  generated by dynamic particles in {SPH} methods. Comput Mater Contin
  5:173--184

\bibitem[{Cuomo et~al(2010)Cuomo, Allsop, Bruce, and
  Pearson}]{Cuomo2010breaking}
Cuomo G, Allsop W, Bruce T, Pearson J (2010) Breaking wave loads at vertical
  seawalls and breakwaters. Coast Eng 57:424--439

\bibitem[{Dehnen and Aly(2012)}]{DeAl2012}
Dehnen W, Aly H (2012) Improving convergence in smoothed particle hydrodynamics
  simulations without pairing instability. Mon Not R Astron Soc 425:1068--1082

\bibitem[{Deierlein et~al(2020{\natexlab{a}})Deierlein, McKenna, Zsarnóczay,
  Kijewski-Correa, Kareem, Elhaddad, Lowes, Schoettler, and Govindjee}]{quoFEM}
Deierlein G, McKenna F, Zsarnóczay A, Kijewski-Correa T, Kareem A, Elhaddad W,
  Lowes L, Schoettler M, Govindjee S (2020{\natexlab{a}}) A cloud-enabled
  application framework for simulating regional-scale impacts of natural
  hazards on the built environment. Frontiers in Built Environment 6:558,706

\bibitem[{Deierlein et~al(2020{\natexlab{b}})Deierlein, McKenna,
  Zsarn{\'o}czay, Kijewski-Correa, Kareem, Elhaddad, Lowes, Schoettler, and
  Govindjee}]{Deierlein2020cloud}
Deierlein GG, McKenna F, Zsarn{\'o}czay A, Kijewski-Correa T, Kareem A,
  Elhaddad W, Lowes L, Schoettler MJ, Govindjee S (2020{\natexlab{b}}) A
  cloud-enabled application framework for simulating regional-scale impacts of
  natural hazards on the built environment. Front Built Environ 6:558,706

\bibitem[{Deshpande et~al(2012)Deshpande, Anumolu, and Trujillo}]{DeAnTr2012}
Deshpande SS, Anumolu L, Trujillo MF (2012) Evaluating the performance of the
  two-phase flow solver interfoam. Comput Sci Discov 5:014,016

\bibitem[{Dom{\'i}nguez et~al(2021)Dom{\'i}nguez, Fourtakas, Altomare, Canelas,
  Tafuni, Garc{\'i}a-Feal, Mart{\'i}nez-Est{\'e}vez, Mokos, Vacondio, Crespo,
  Rogers, Stansby, and G{\'o}mez-Gesteira}]{dom2021}
Dom{\'i}nguez J, Fourtakas G, Altomare C, Canelas R, Tafuni A, Garc{\'i}a-Feal
  O, Mart{\'i}nez-Est{\'e}vez I, Mokos A, Vacondio R, Crespo A, Rogers BD,
  Stansby P, G{\'o}mez-Gesteira M (2021) {DualSPHysics}: From fluid dynamics to
  multiphysics problems. Comp Part Mech 9:867--895

\bibitem[{Eberl(2016)}]{Fisher_Yates-AFP}
Eberl M (2016) Fisher–yates shuffle. Archive of Formal Proofs

\bibitem[{Engsig-Karup et~al(2009)Engsig-Karup, Bingham, and
  Lindberg}]{OceanWave3D}
Engsig-Karup AP, Bingham HB, Lindberg O (2009) An efficient flexible-order
  model for {3D} nonlinear water waves. J Comput Phys 228:2100--2118

\bibitem[{Fourtakas et~al(2019)Fourtakas, Dominguez, Vacondio, and
  Rogers}]{FOURTAKAS2019346}
Fourtakas G, Dominguez J, Vacondio R, Rogers BD (2019) Local uniform stencil
  ({LUST}) boundary condition for arbitrary 3-d boundaries in parallel smoothed
  particle hydrodynamics ({SPH}) models. Comp Fluids 190:346--361

\bibitem[{Ghaffary and Moustafa(2021)}]{Ghaffary2021Performance}
Ghaffary A, Moustafa M (2021) Performance-based assessment and structural
  response of 20-storey {SAC} building under wind hazards through collapse. J
  of Struct Eng 147:04020,346

\bibitem[{Goda(1974)}]{Goda1974New}
Goda Y (1974) New wave pressure formulae for composite breakwaters, ASCE, pp
  1702--1720. \doi{10.1061/9780872621138.103}

\bibitem[{Hu et~al(2016)Hu, Greaves, and Raby}]{Hu2016-ux}
Hu ZZ, Greaves D, Raby A (2016) Numerical wave tank study of extreme waves and
  wave-structure interaction using {OpenFoam}. Ocean Eng 126:329--342

\bibitem[{Huang(2022)}]{Huang2022-lj}
Huang Le (2022) A review on the modelling of wave-structure interactions based
  on {OpenFOAM}. OpenFOAM J 2:116--142

\bibitem[{Jacobsen et~al(2012)Jacobsen, Fuhrman, and
  Freds{\o}e}]{Jacobsen2012-gs}
Jacobsen NG, Fuhrman DR, Freds{\o}e J (2012) A wave generation toolbox for the
  open-source {CFD} library: {OpenFoam}. Int J Numer Methods Fluids
  70:1073--1088

\bibitem[{Judd(2018)}]{Judd2018Windstorm}
Judd J (2018) Windstorm resilience of a 10-storey steel frame office building.
  ASCE-ASME J Risk Uncertain Eng Syst A 4:04018,020

\bibitem[{Leimkuhler and Matthews(2015)}]{leimkuhler2015molecular}
Leimkuhler B, Matthews C (2015) Molecular dynamics. Springer

\bibitem[{LeVeque et~al(2011)LeVeque, George, and Berger}]{Leveq2011}
LeVeque RJ, George DL, Berger MJ (2011) Adaptive mesh refinement techniques for
  tsunamis and other geophysical flows over topography. Acta Numer 20:211--289

\bibitem[{Liu(2009)}]{Liu2003book3}
Liu GR (2009) Meshfree methods: {M}oving beyond the finite element method. CRC
  Press

\bibitem[{Liu and Gu(2005)}]{Liu2003book2}
Liu GR, Gu YT (2005) An introduction to meshfree methods and their programming.
  Springer Dordrecht

\bibitem[{Liu et~al(2003)Liu, Liu, and Liu}]{Liu2003book}
Liu GR, Liu MB, Liu GR (2003) Smoothed particle hydrodynamics: {A} meshfree
  particle method. World Scientific Publishing Company

\bibitem[{Liu and Liu(2010)}]{liu_smoothed_2010}
Liu MB, Liu GR (2010) Smoothed prticle hydrodynamics ({SPH}): {A}n overview and
  recent developments. Arch Comput Methods Eng 17:25--76

\bibitem[{Luettich and Westerink(1991)}]{Adcirc}
Luettich RA, Westerink JJ (1991) A solution for the vertical variation of
  stress, rather than velocity, in a three-dimensional circulation model. Int J
  Numer Methods Fluids 12:911--928

\bibitem[{Lyu et~al(2022)Lyu, Sun, Huang, Zhong, Peng, Jiang, and
  Ji}]{lyu_review_2022}
Lyu HG, Sun PN, Huang XT, Zhong SY, Peng YX, Jiang T, Ji CN (2022) A review of
  {SPH} techniques for hydrodynamic simulations of ocean energy devices.
  Energies 15:502

\bibitem[{Maraveas and Tsavdaridis(2019)}]{Maraveas2019Assessment}
Maraveas C, Tsavdaridis K (2019) Assessment and retrofitting of an existing
  steel structure subjected to wind-induced failure analysis. J of Build Eng
  23:53--67

\bibitem[{McKay and et~al.(1979)}]{LHSsampling}
McKay MD, et~al (1979) Comparison of three methods for selecting values of
  input variables in the analysis of output from a computer code. Technometrics
  21:239–245

\bibitem[{McKenna et~al(2022)McKenna, Yi, Satish, Zsarnoczay, Gardner, Zhong,
  and Elhaddad}]{McKenna2022NHERI}
McKenna F, Yi S, Satish A, Zsarnoczay A, Gardner M, Zhong K, Elhaddad W (2022)
  Nheri-{SimCenter}/{quoFEM}: Version 3.0.0. Zenodo

\bibitem[{Mohammadi et~al(2019)Mohammadi, Azizinamini, Griffis, and
  Irwin}]{Mohammadi2019Performance}
Mohammadi A, Azizinamini A, Griffis L, Irwin P (2019) Performance assessment of
  an existing 47-storey high-rise building under extreme wind loads. J of
  Struct Eng 145:04018,232

\bibitem[{Monaghan et~al(1999)Monaghan, Cas, Kos, and Hallworth}]{monaghan1999}
Monaghan J, Cas R, Kos A, Hallworth M (1999) Gravity currents descending a ramp
  in a stratified tank. J Fluid Mech 379:39--69

\bibitem[{Monaghan(1992)}]{Mo1992}
Monaghan JJ (1992) Smoothed particle hydrodynamics. Annu Rev Astron Astr
  30:543--574

\bibitem[{Moris et~al(2021)Moris, Kennedy, and Westerink}]{Joetal2021}
Moris JP, Kennedy AB, Westerink JJ (2021) Tsunami wave run-up load reduction
  inside a building array. Coast Eng 169:103,910

\bibitem[{Parshikov et~al(2000)Parshikov, Medin, Loukashenko, and
  Milekhin}]{parshikov2000improvements}
Parshikov AN, Medin SA, Loukashenko II, Milekhin VA (2000) Improvements in sph
  method by means of interparticle contact algorithm and analysis of
  perforation tests at moderate projectile velocities. International Journal of
  Impact Engineering 24:779--796

\bibitem[{Pringgana et~al(2023)Pringgana, Cunningham, and
  Rogers}]{BenLikeArticle}
Pringgana G, Cunningham LS, Rogers BD (2023) Mitigating tsunami effects on
  buildings via novel use of discrete onshore protection systems. Coast Eng
  65:149--173

\bibitem[{Ribberink et~al(2012)Ribberink, Dohmen-Janssen, Hanes, McLean, and
  Vincent}]{Ribberink2001near}
Ribberink JS, Dohmen-Janssen CM, Hanes DM, McLean SR, Vincent C (2012) Near-bed
  sand transport mechanisms under waves - {A} large-scale flume experiment
  {(Sistex99)}, ASCE, pp 3263--3276. \doi{10.1061/40549(276)254}

\bibitem[{Roenby et~al(2016)Roenby, Bredmose, and Jasak}]{RoBrJa2016}
Roenby J, Bredmose H, Jasak H (2016) A computational method for sharp interface
  advection. R Soc Open Sci 3:160,405

\bibitem[{Roselli et~al(2018)Roselli, Vernengo, Altomare, Brizzolara,
  Bonfiglio, and Guercio}]{rota2018}
Roselli ARR, Vernengo G, Altomare C, Brizzolara S, Bonfiglio L, Guercio R
  (2018) Ensuring numerical stability of wave propagation by tuning model
  parameters using genetic algorithms and response surface methods. Environ
  Model Softw 103:62--73

\bibitem[{Sampath et~al(2016)Sampath, Montanari, Akinci, Prescott, and
  Smith}]{solitarySPH_2016}
Sampath R, Montanari N, Akinci N, Prescott S, Smith C (2016) Large-scale
  solitary wave simulation with implicity incompressible {SPH}. J Ocean Eng Mar
  Energy 2:313--329

\bibitem[{Serre(1953)}]{Rayleigh}
Serre F (1953) Contribution {\`a} l{\'e}tude des {\'e}coulements permanents et
  variables dans les canaux. La Houille B lanche 8:374--388

\bibitem[{Tagliafierro et~al(2023)Tagliafierro, Karimirad, Altomare, Göteman,
  Martínez-Estévez, Capasso, Domínguez, Viccione, Gómez-Gesteira, and
  Crespo}]{Tagliafierro2023}
Tagliafierro B, Karimirad M, Altomare C, Göteman M, Martínez-Estévez I,
  Capasso S, Domínguez JM, Viccione G, Gómez-Gesteira M, Crespo AJ (2023)
  Numerical validations and investigation of a semi-submersible floating
  offshore wind turbine platform interacting with ocean waves using an sph
  framework. Applied Ocean Research 141:103,757

\bibitem[{Takahashi et~al(1994)Takahashi, Tanimoto, and
  Shimosako}]{TAKAHASHI94}
Takahashi S, Tanimoto K, Shimosako K (1994) A proposal of impulsive pressure
  coefficient for design of composite breakwaters. Proc Int Conf Hydro- Tech
  Eng Port Harbor Constr (Hydro-Port '94)

\bibitem[{Tomiczek et~al(2016)Tomiczek, Prasetyo, Mori, Yasuda, and
  Kennedy}]{TOMICZEK201697}
Tomiczek T, Prasetyo A, Mori N, Yasuda T, Kennedy A (2016) Physical modelling
  of tsunami onshore propagation, peak pressures, and shielding effects in an
  urban building array. Coast Eng 117:97--112

\bibitem[{Tryggvason et~al(2011)Tryggvason, Scardovelli, and
  Zaleski}]{Trybook2011}
Tryggvason G, Scardovelli R, Zaleski S (2011) Direct numerical simulations of
  gas–liquid multiphase flows. Cambridge University Press

\bibitem[{Vignjevic and Campbell(2009)}]{hiermaier_review_2009}
Vignjevic R, Campbell J (2009) Review of development of the smooth particle
  hydrodynamics ({SPH}) method. In: Hiermaier S (ed) Predictive modeling of
  dynamic processes, Springer US, Boston, MA, pp 367--396

\bibitem[{van Thiel~de Vries et~al(2008)van Thiel~de Vries, van Gent, Walstra,
  and Reniers}]{De2008analysis}
van Thiel~de Vries JSM, van Gent MRA, Walstra DJR, Reniers AJHM (2008) Analysis
  of dune erosion processes in large-scale flume experiments. Coast Eng
  55:1028--1040

\bibitem[{Wendland(1995)}]{Wendland1995}
Wendland H (1995) Piecewise polynomial, positive definite and compactly
  supported radial functions of minimal degree. Adv in Comp Math 4:389--396

\bibitem[{Ye et~al(2019)Ye, Pan, Huang, and Liu}]{ye_smoothed_2019}
Ye T, Pan D, Huang C, Liu M (2019) Smoothed particle hydrodynamics ({SPH}) for
  complex fluid flows: {Recent} developments in methodology and applications.
  Phys Fluids 31:011,301

\bibitem[{Zhang et~al(2018)Zhang, Crespo, Altomare, Dominguez, Marzeddu, Shang,
  and Gesteira}]{zhang2018}
Zhang F, Crespo A, Altomare C, Dominguez J, Marzeddu A, Shang Sp, Gesteira M
  (2018) {DualSPHysics: A} numerical tool to simulate real breakwaters. J
  Hydrodyn 30:95–105

\end{thebibliography}


\appendix
\section{Appendix A: Wave height contours}
\begin{figure}
 \centering
 \begin{subfigure}{0.9\textwidth}
   \centering
   \includegraphics[width=\textwidth]{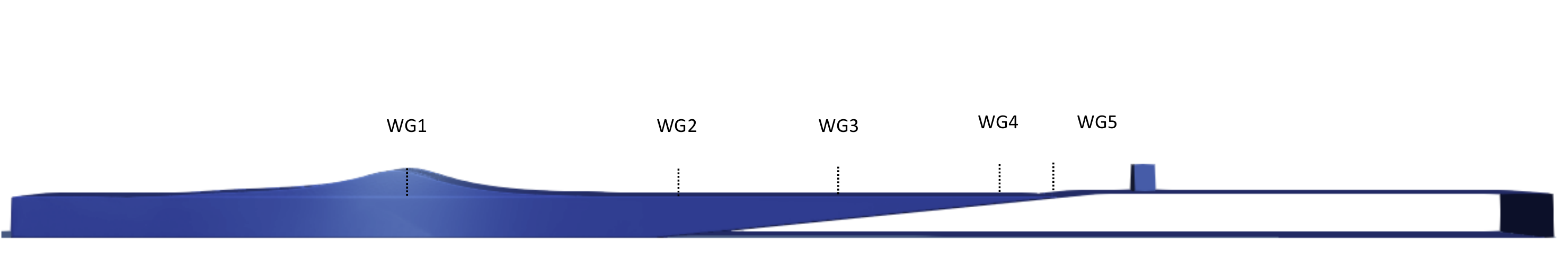}
   \caption{}
   \label{fig:WH4-01}
 \end{subfigure}
 \begin{subfigure}{0.9\textwidth}
   \centering
   \includegraphics[width=\textwidth]{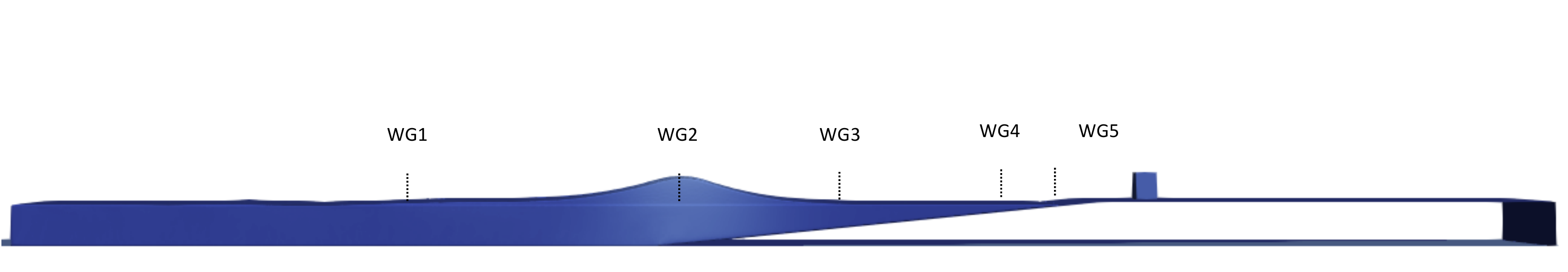}
   \caption{}
   \label{fig:WH4-02}
 \end{subfigure}
 \begin{subfigure}{0.9\textwidth}
   \centering
   \includegraphics[width=\textwidth]{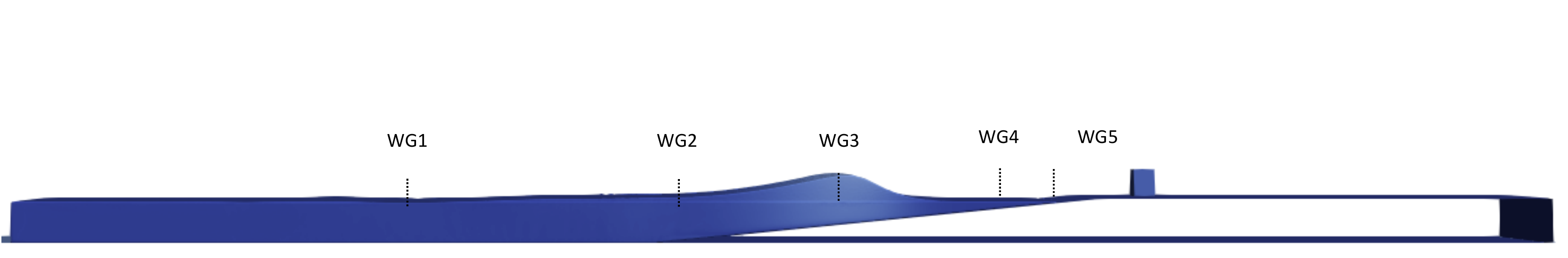}
   \caption{}
   \label{fig:WH4-03}
 \end{subfigure}
 \begin{subfigure}{0.9\textwidth}
   \centering
   \includegraphics[width=\textwidth]{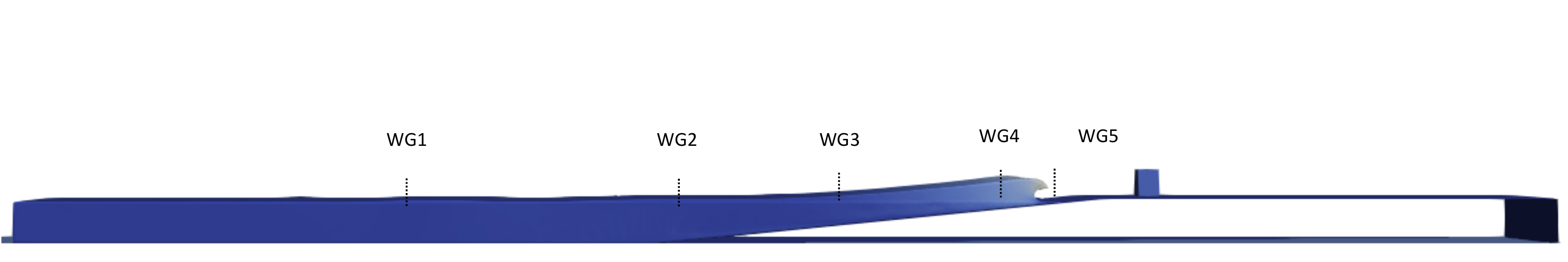}
   \caption{}
   \label{fig:WH4-04}
 \end{subfigure}
  \begin{subfigure}{0.9\textwidth}
   \centering
   \includegraphics[width=\textwidth]{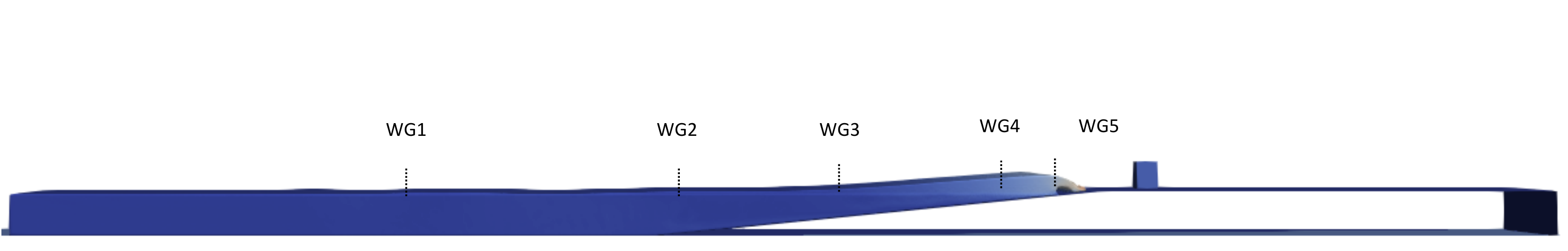}
   \caption{}
   \label{fig:WH4-05}
 \end{subfigure}
  \begin{subfigure}{0.9\textwidth}
   \centering
   \includegraphics[width=\textwidth]{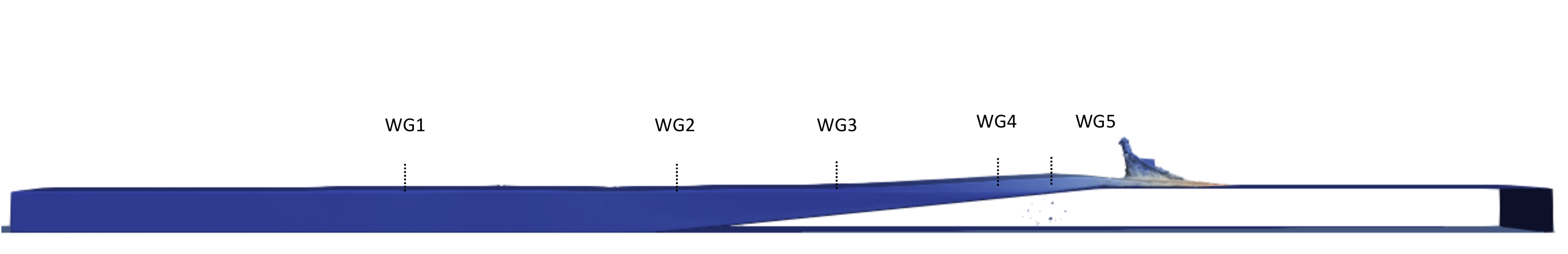}
   \caption{}
   \label{fig:WH4-06}
 \end{subfigure}
 \begin{subfigure}{0.9\textwidth}
   \centering
   \includegraphics[width=0.4\textwidth]{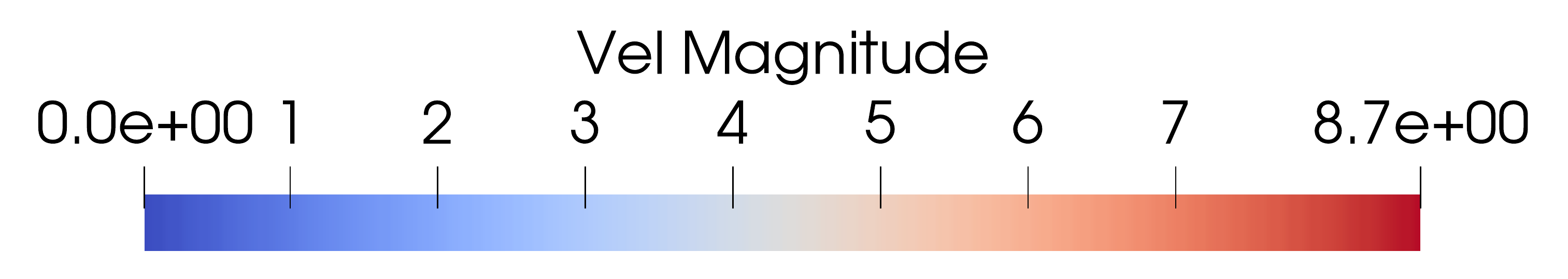}
   \caption{}
   \label{fig:vel-mag04}
 \end{subfigure}
\caption{Propagation of the solitary wave for the initial wave height of 0.40 \unit{\meter}.}
\label{fig:waveheight4}
\end{figure}

\begin{figure}
 \centering
 \begin{subfigure}{0.9\textwidth}
   \centering
   \includegraphics[width=\textwidth]{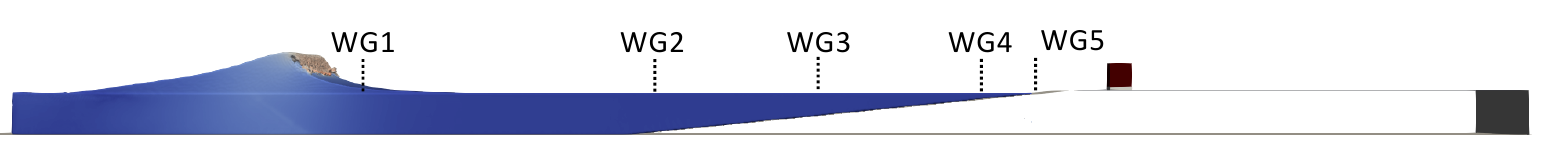}
   \caption{}
   \label{fig:WH6-01}
 \end{subfigure}
 \begin{subfigure}{0.9\textwidth}
   \centering
   \includegraphics[width=\textwidth]{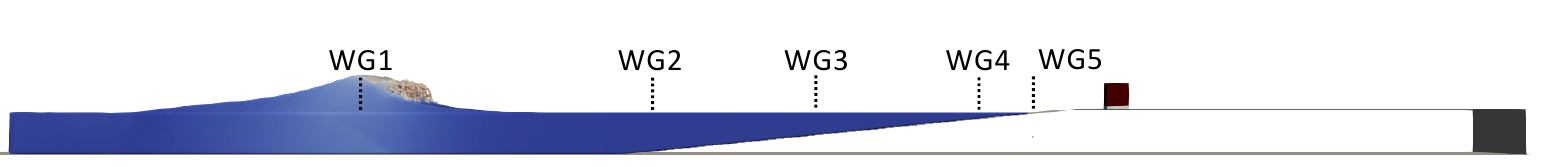}
   \caption{}
   \label{fig:WH6-02}
 \end{subfigure}
 \begin{subfigure}{0.9\textwidth}
   \centering
   \includegraphics[width=\textwidth]{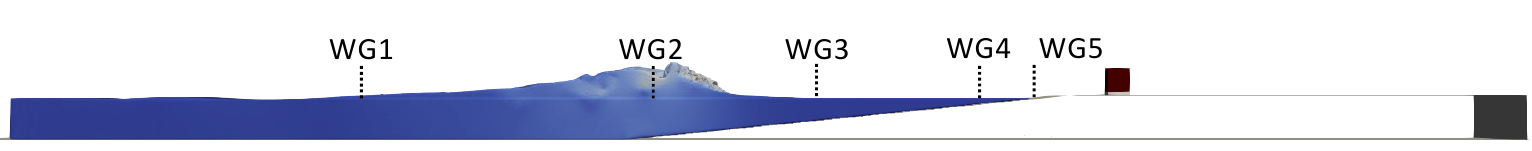}
   \caption{}
   \label{fig:WH6-03}
 \end{subfigure}
 \begin{subfigure}{0.9\textwidth}
   \centering
   \includegraphics[width=\textwidth]{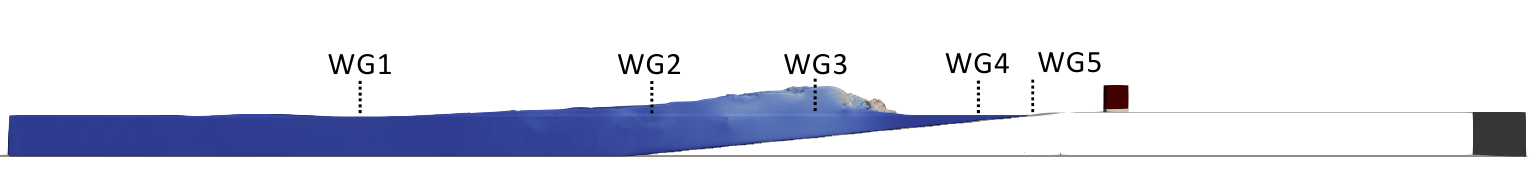}
   \caption{}
   \label{fig:WH6-04}
 \end{subfigure}
 \begin{subfigure}{0.9\textwidth}
   \centering
   \includegraphics[width=\textwidth]{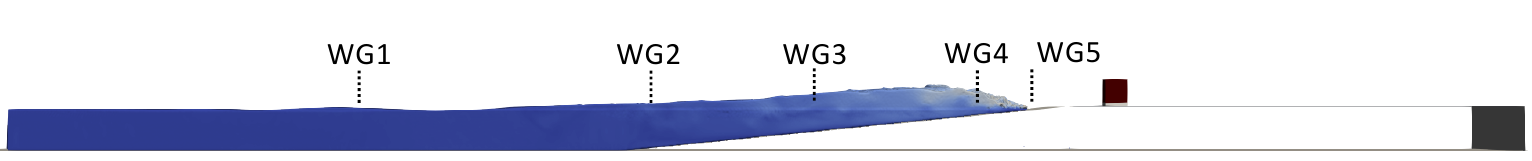}
   \caption{}
   \label{fig:WH6-05}
 \end{subfigure}
  \begin{subfigure}{0.9\textwidth}
   \centering
   \includegraphics[width=\textwidth]{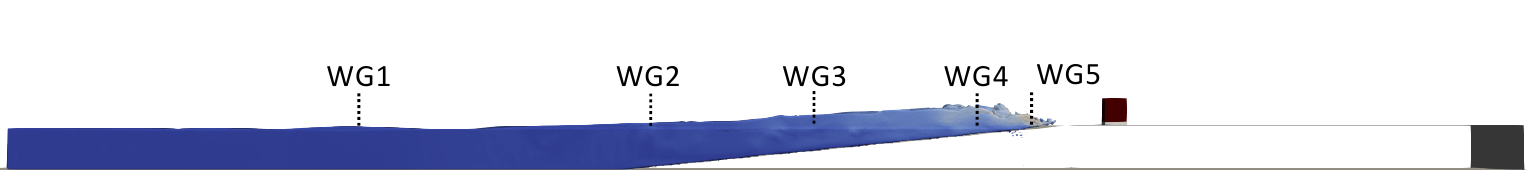}
   \caption{}
   \label{fig:WH6-06}
 \end{subfigure}
  \begin{subfigure}{0.9\textwidth}
   \centering
   \includegraphics[width=\textwidth]{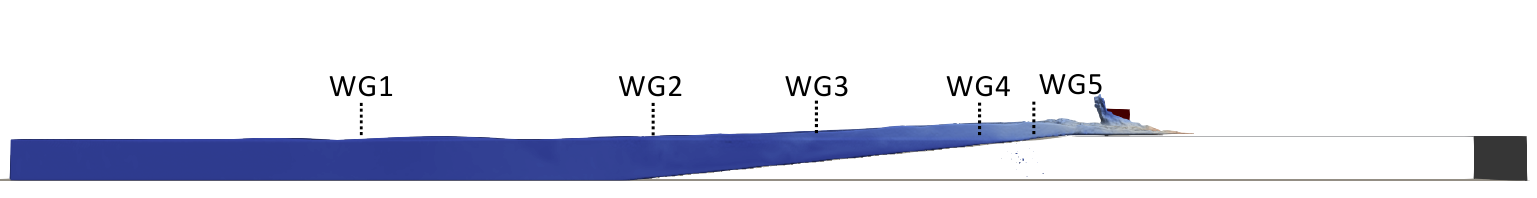}
   \caption{}
   \label{fig:WH6-07}
 \end{subfigure}
 \begin{subfigure}{0.9\textwidth}
   \centering
   \includegraphics[width=0.4\textwidth]{Appendix/images/vel_magnitude}
   \caption{}
   \label{fig:vel-mag06}
 \end{subfigure}
\caption{Propagation of the solitary wave for the initial wave height of 0.60 \unit{\meter}.}
\label{fig:waveheight6}
\end{figure}

\begin{figure}
 \centering
 \begin{subfigure}{0.9\textwidth}
   \centering
   \includegraphics[width=\textwidth]{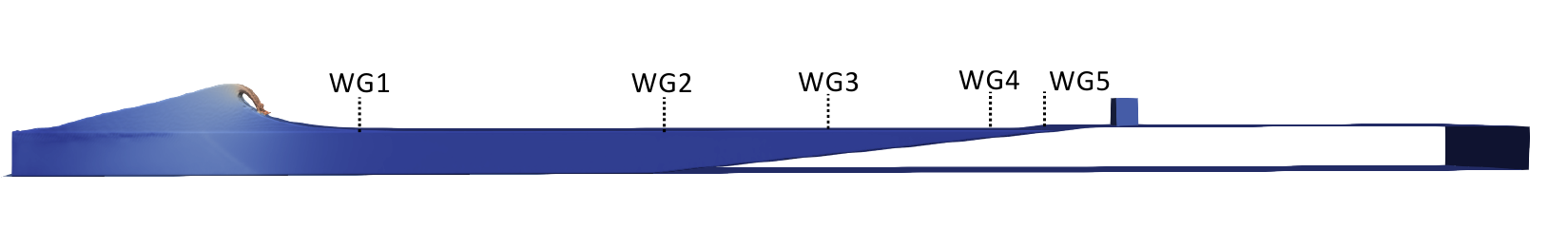}
   \caption{}
   \label{fig:WH7-01}
 \end{subfigure}
 \begin{subfigure}{0.9\textwidth}
   \centering
   \includegraphics[width=\textwidth]{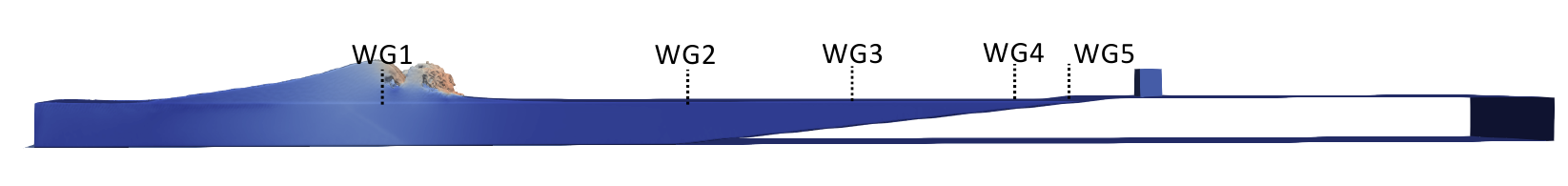}
   \caption{}
   \label{fig:WH7-02}
 \end{subfigure}
 \begin{subfigure}{0.9\textwidth}
   \centering
   \includegraphics[width=\textwidth]{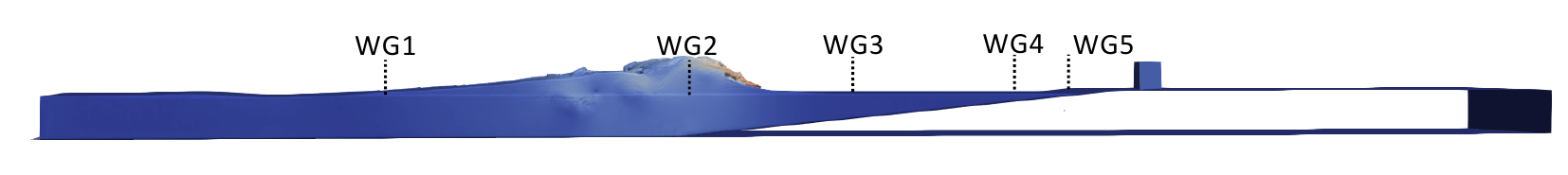}
   \caption{}
   \label{fig:WH7-03}
 \end{subfigure}
 \begin{subfigure}{0.9\textwidth}
   \centering
   \includegraphics[width=\textwidth]{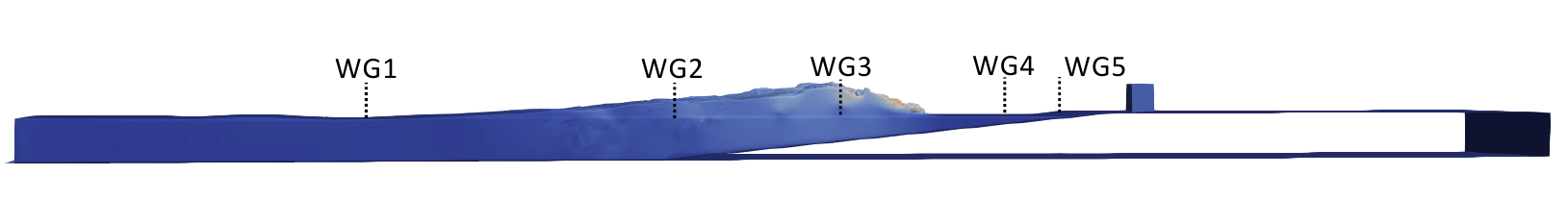}
   \caption{}
   \label{fig:WH7-04}
 \end{subfigure}
  \begin{subfigure}{0.9\textwidth}
   \centering
   \includegraphics[width=\textwidth]{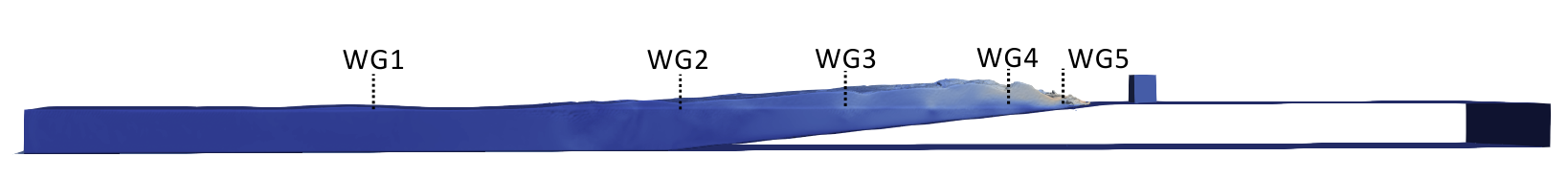}
   \caption{}
   \label{fig:WH7-05}
 \end{subfigure}
  \begin{subfigure}{0.9\textwidth}
   \centering
   \includegraphics[width=\textwidth]{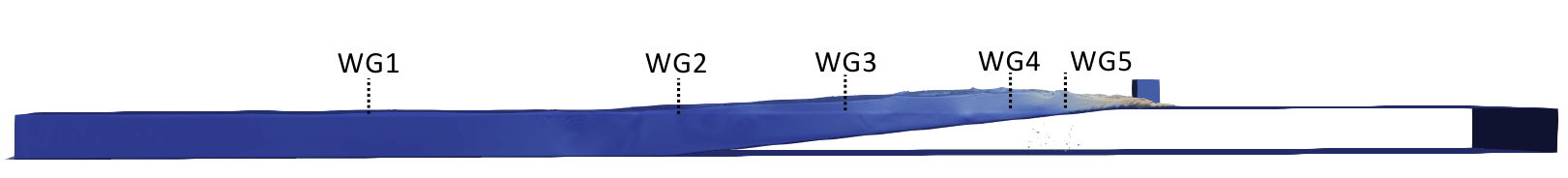}
   \caption{}
   \label{fig:WH7-06}
 \end{subfigure}
   \begin{subfigure}{0.9\textwidth}
   \centering
   \includegraphics[width=\textwidth]{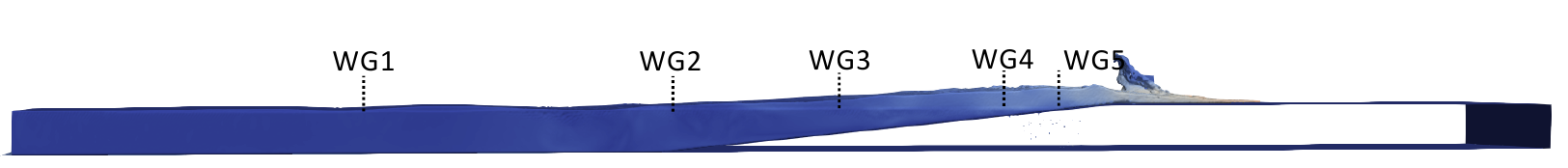}
   \caption{}
   \label{fig:WH7-07}
 \end{subfigure}
 \begin{subfigure}{0.9\textwidth}
   \centering
   \includegraphics[width=0.4\textwidth]{Appendix/images/vel_magnitude}
   \caption{}
   \label{fig:vel-mag07}
 \end{subfigure}
\caption{Propagation of the solitary wave for the initial wave height of 0.70 \unit{\meter}.}
\label{fig:waveheight7}
\end{figure}

\begin{figure}
 \centering
 \begin{subfigure}{0.9\textwidth}
   \centering
   \includegraphics[width=\textwidth]{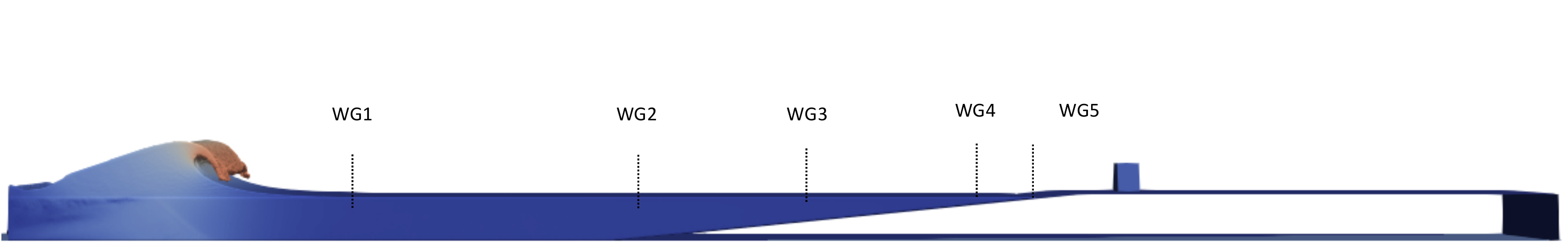}
   \caption{}
   \label{fig:WH9-01}
 \end{subfigure}
 \begin{subfigure}{0.9\textwidth}
   \centering
   \includegraphics[width=\textwidth]{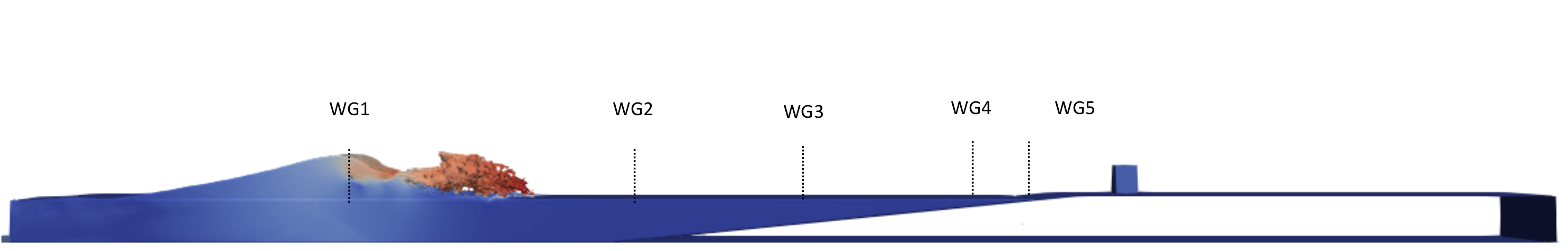}
   \caption{}
   \label{fig:WH9-02}
 \end{subfigure}
 \begin{subfigure}{0.9\textwidth}
   \centering
   \includegraphics[width=\textwidth]{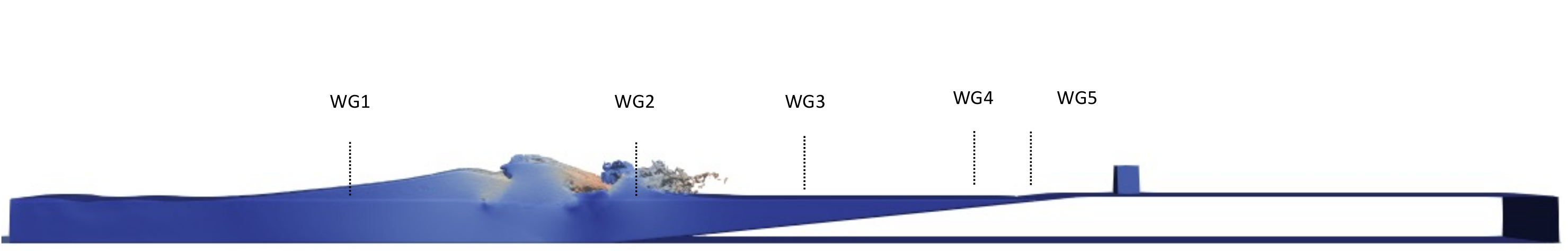}
   \caption{}
   \label{fig:WH9-03}
 \end{subfigure}
 \begin{subfigure}{0.9\textwidth}
   \centering
   \includegraphics[width=\textwidth]{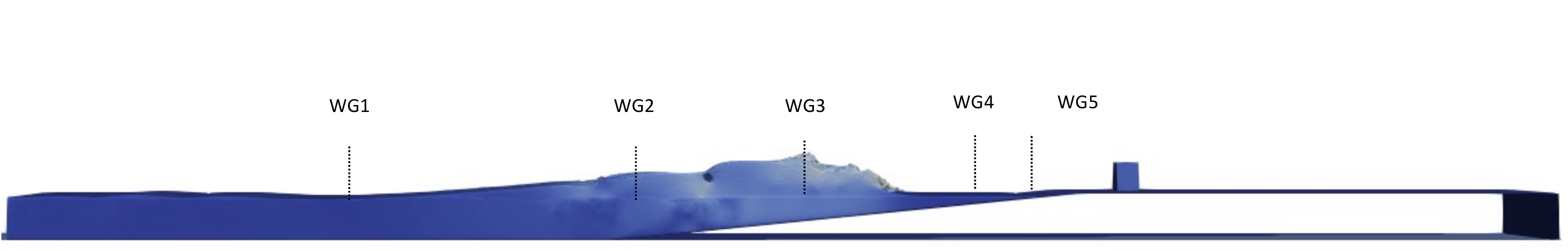}
   \caption{}
   \label{fig:WH9-04}
 \end{subfigure}
  \begin{subfigure}{0.9\textwidth}
   \centering
   \includegraphics[width=\textwidth]{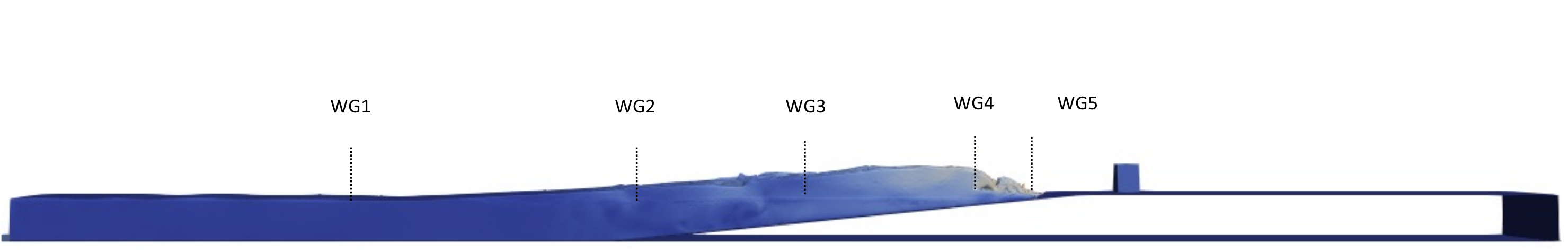}
   \caption{}
   \label{fig:WH9-05}
 \end{subfigure}
  \begin{subfigure}{0.9\textwidth}
   \centering
   \includegraphics[width=\textwidth]{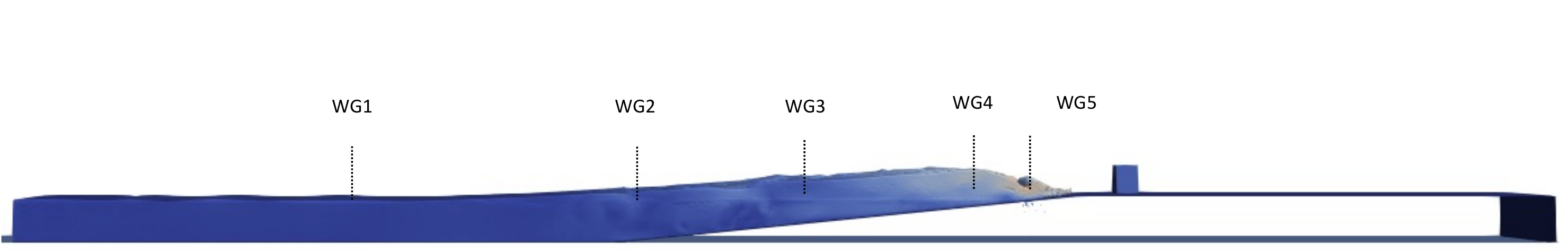}
   \caption{}
   \label{fig:WH9-06}
 \end{subfigure}
   \begin{subfigure}{0.9\textwidth}
   \centering
   \includegraphics[width=\textwidth]{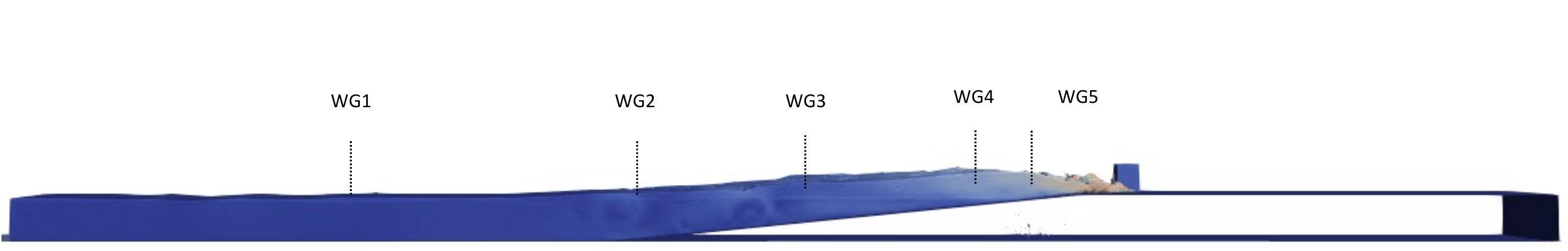}
   \caption{}
   \label{fig:WH9-07}
 \end{subfigure}
    \begin{subfigure}{0.9\textwidth}
   \centering
   \includegraphics[width=\textwidth]{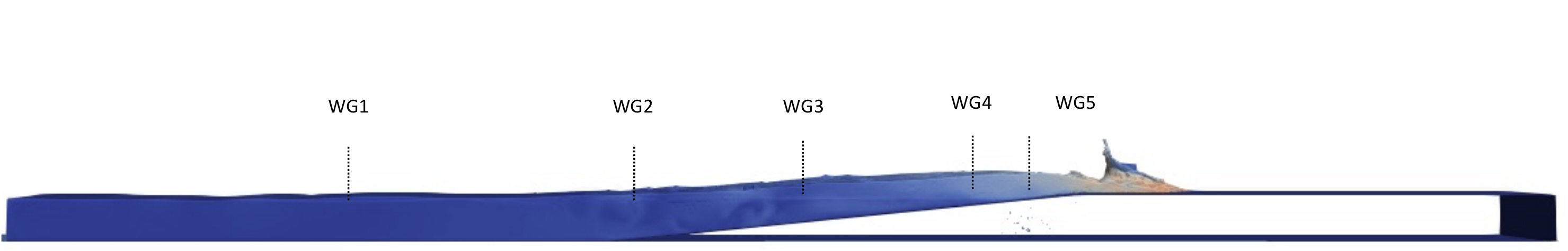}
   \caption{}
   \label{fig:WH9-08}
 \end{subfigure}
 \begin{subfigure}{0.9\textwidth}
   \centering
   \includegraphics[width=0.4\textwidth]{Appendix/images/vel_magnitude}
   \caption{}
   \label{fig:vel-mag09}
 \end{subfigure}
\caption{Propagation of the solitary wave for the initial wave height of 0.90 \unit{\meter}.}
\label{fig:waveheight9}
\end{figure}

The wave contours (side-view) are shown for three wave heights, namely 0.40 \unit{\meter} (\autoref{fig:waveheight4}), 0.60 \unit{\meter} (\autoref{fig:waveheight6}), 0.70 \unit{\meter}  (\autoref{fig:waveheight7}) and 0.90 \unit{\meter} (\autoref{fig:waveheight9}). As discussed in the paper, the overall drag force on the structure reduces from 0.40 - 0.70 \unit{\meter} and further increases to 0.90 \unit{\meter}.

\begin{figure}[!htb]
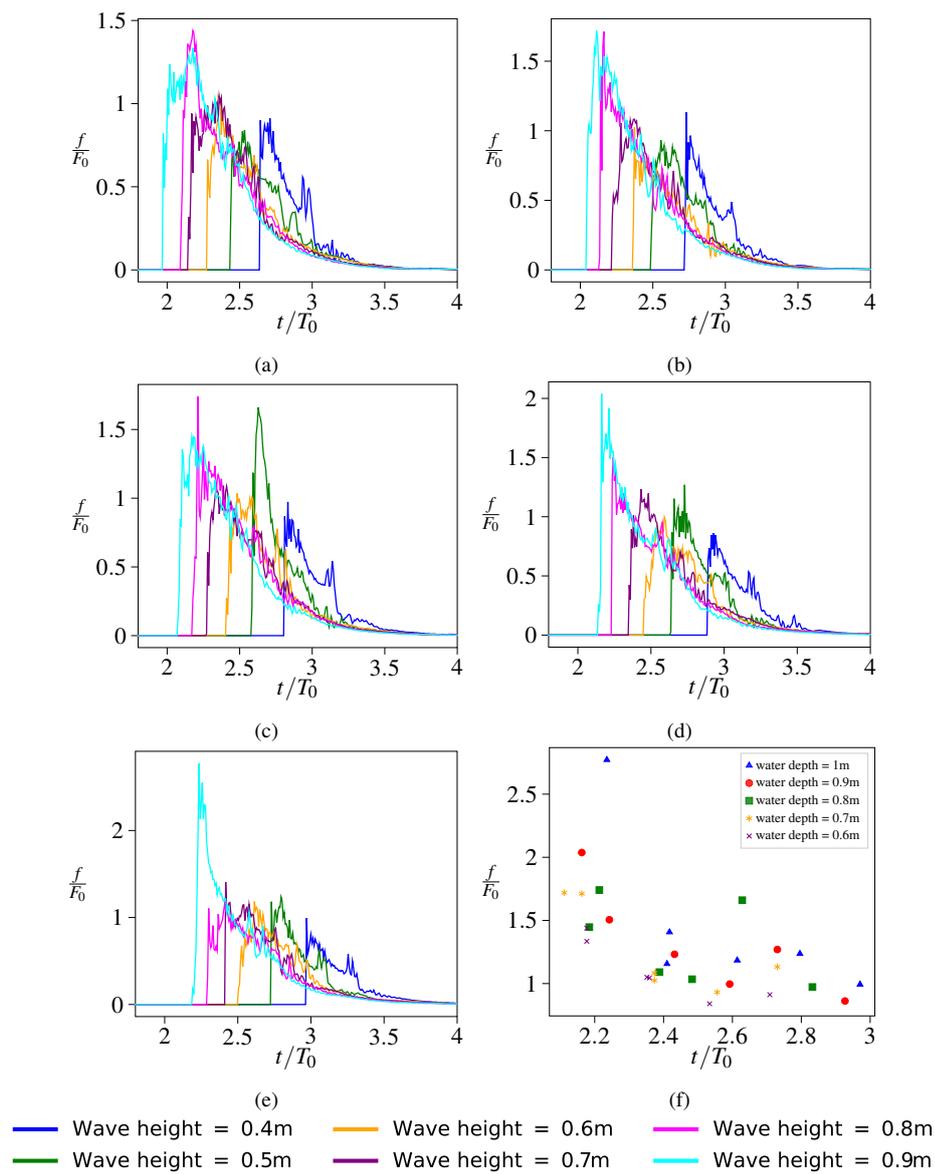

  \centering
   \begin{subfigure}{0.44\textwidth}
    \centering
    \adjustbox{width=\linewidth}{ \input{Appendix/images/Diff_swh_force/swh06_Diffwvh_force_normalized}}
    \caption{}
    \label{fig:force_swh06}
  \end{subfigure} 
  \begin{subfigure}{0.44\textwidth}
    \centering
    \adjustbox{width=\linewidth}{ \input{Appendix/images/Diff_swh_force/swh07_Diffwvh_force_normalized}}
    \caption{}
    \label{fig:force_swh07}
  \end{subfigure} 
  \begin{subfigure}{0.44\textwidth}
    \centering
    \adjustbox{width=\linewidth}{ \input{Appendix/images/Diff_swh_force/swh08_Diffwvh_force_normalized}}
    \caption{}
    \label{fig:force_swh08}
  \end{subfigure}
    \begin{subfigure}{0.44\textwidth}
    \centering
    \adjustbox{width=\linewidth}{\input{Appendix/images/Diff_swh_force/swh09_Diffwvh_force_normalized}}
    \caption{}
    \label{fig:force_swh09}
  \end{subfigure}
    \begin{subfigure}{0.44\textwidth}
    \centering
    \adjustbox{width=\linewidth}{\input{Appendix/images/Diff_swh_force/swh1_Diffwvh_force_normalized}}
    \caption{}
    \label{fig:force_swh1}
  \end{subfigure}
  \begin{subfigure}{0.44\textwidth}
    \centering
    \adjustbox{width=\linewidth}{
\begin{tikzpicture}

\definecolor{darkgray176}{RGB}{176,176,176}
\definecolor{green}{RGB}{0,128,0}
\definecolor{lightgray204}{RGB}{204,204,204}
\definecolor{orange}{RGB}{255,165,0}
\definecolor{purple}{RGB}{128,0,128}

\begin{axis}[
legend cell align={left},
legend style={fill opacity=0.8, draw opacity=1, text opacity=1, draw=lightgray204},
tick align=outside,
tick pos=left,
x grid style={darkgray176},
xlabel={$t/T_0$},
xmin=2.06850259811919, xmax=3.01356858557662,
xtick style={color=black},
y grid style={darkgray176},
ylabel style={rotate=-90},
ylabel={$\frac{f}{F_0}$},
ymin=0.743948811459691, ymax=2.86847920845067,
ytick style={color=black}
]
\addplot [draw=blue, fill=blue, mark=triangle*, only marks]
table{%
x  y
2.97061104069219 0.993176328220288
2.79585846477585 1.23651336965468
2.61383703035982 1.18401121815094
2.41725312397982 1.40730573515784
2.40998790591611 1.15644953092086
2.23523824234854 2.77190964495108
};
\addlegendentry{water depth = 1m}
\addplot [draw=red, fill=red, mark=*, only marks]
table{%
x  y
2.92691889223354 0.862035118137554
2.73033789820232 1.26974047248317
2.59199841902423 0.995899388508921
2.43182360490292 1.23185901401137
2.24250892106621 1.50549541209209
2.16241969378757 2.0375993859741
};
\addlegendentry{water depth = 0.9m}
\addplot [draw=green, fill=green, mark=square*, only marks]
table{%
x  y
2.83226500872936 0.972706148709065
2.62839513380062 1.66001418312806
2.48278679612284 1.03433191109525
2.38813400474945 1.09057557736748
2.21339016587943 1.7405145811036
2.18426230960272 1.44688035016115
};
\addlegendentry{water depth = 0.8m}
\addplot [draw=orange, fill=orange, mark=asterisk, only marks]
table{%
x  y
2.73033680607153 1.13273695107396
2.55560316042223 0.93136210674518
2.37356498000074 1.02512466162662
2.37358318218058 1.08291012825468
2.16242333422354 1.71194022310935
2.11146014300362 1.71906082803022
};
\addlegendentry{water depth = 0.7m}
\addplot [draw=purple, fill=purple, mark=x, only marks]
table{%
x  y
2.70849746664875 0.911802832950647
2.53375581204031 0.840518374959281
2.35173619784226 1.04977982246201
2.35901561360626 1.04483921724
2.17699417919023 1.43859295085899
2.17699053875426 1.33491240797396
};
\addlegendentry{water depth = 0.6m}
\end{axis}

\end{tikzpicture}}
    \caption{}
    \label{fig:peakforce_comparison}
  \end{subfigure}
 \begin{subfigure}{\textwidth}
    \centering
    \includegraphics[width=\columnwidth]{Section4/images/waveheights/wgh_legend.pdf}
    \label{fig:legend_Diff_swh_peak} 
  \end{subfigure}  
  \caption{Normalised force time history obtained from (a) water depth is 0.6m (b) water depth is 0.7m (c) water depth is 0.8m (d) water depth is 0.9m (d) water depth is 1m (e)the summary extracts only the peak values from (a) to (d) results. $f$ represents the force measured/calculated, $F_0$ represent the characteristic force  ($F_0 = 695$ \unit{\newton}), $t$ indicates the simulation time, and $T_0$ represent the characteristic time ($T_0 = 2.747$ \unit{\second}).}
  \label{fig:forceComparison2}
\end{figure}

\begin{figure}[!htb]
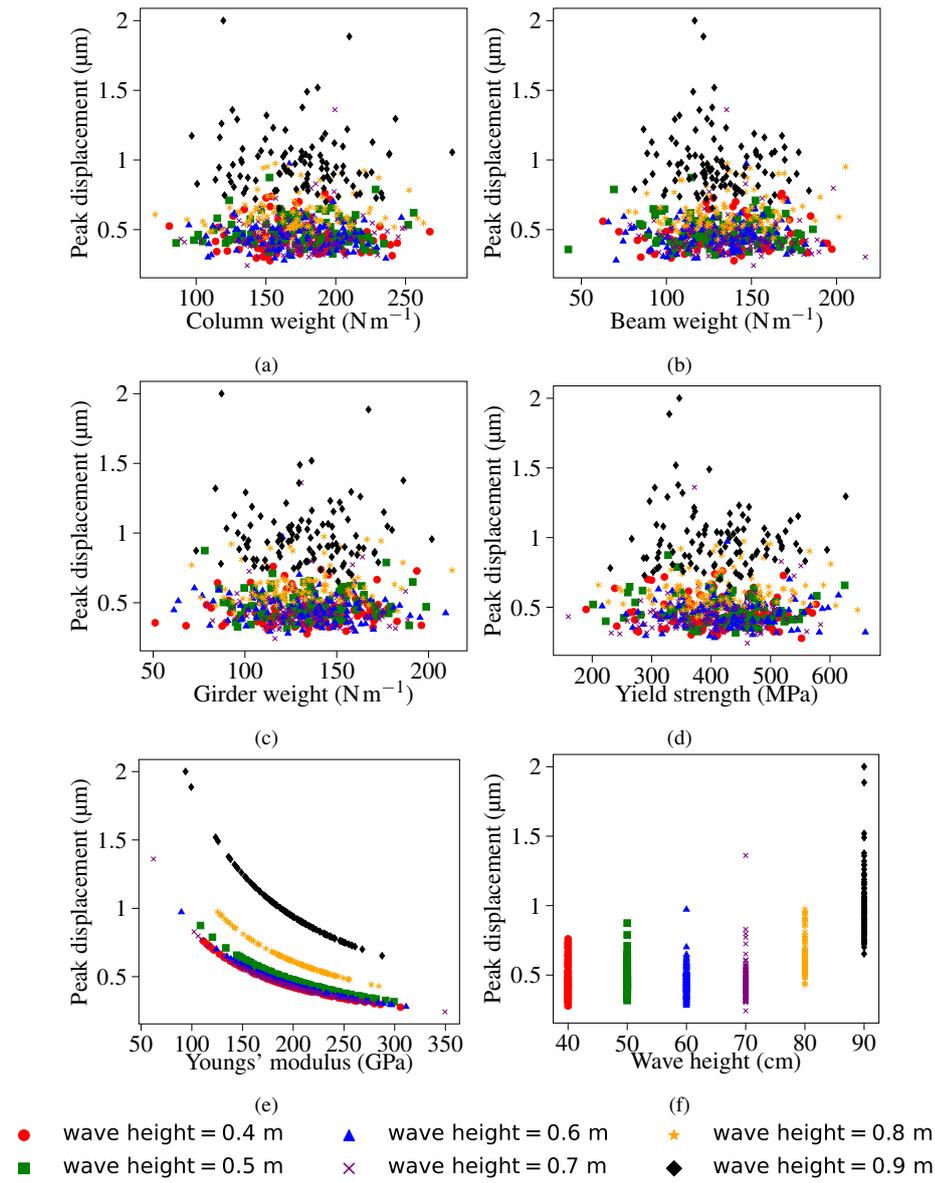

 \centering
 \begin{subfigure}{0.44\textwidth}
   \centering
   \adjustbox{width=\linewidth}{ \input{Section4/images/uq/OpenSees/myQdl_peakdisp}}
   \caption{}
   \label{fig:myQdl_peakdisp}
 \end{subfigure}
 \begin{subfigure}{0.44\textwidth}
   \centering
   \adjustbox{width=\linewidth}{ \input{Section4/images/uq/OpenSees/myQBeam_peakdisp}}
   \caption{}
   \label{fig:myQBeam_peakdisp}
 \end{subfigure}
 \begin{subfigure}{0.44\textwidth}
   \centering
   \adjustbox{width=\linewidth}{ \input{Section4/images/uq/OpenSees/myQGrid_peakdisp}}
   \caption{}
   \label{fig:myQGrid_peakdisp}
 \end{subfigure}
  \begin{subfigure}{0.44\textwidth}
   \centering
   \adjustbox{width=\linewidth}{ \input{Section4/images/uq/OpenSees/myfy_peakdisp}}
   \caption{}
   \label{fig:myfy_peakdisp}
 \end{subfigure}
 \begin{subfigure}{0.44\textwidth}
   \centering
   \adjustbox{width=\linewidth}{ \input{Section4/images/uq/OpenSees/myEs_peakdisp}}
   \caption{}
   \label{fig:myEsl_peakdisp}
 \end{subfigure}
 \begin{subfigure}{0.44\textwidth}
   \centering
   \adjustbox{width=\linewidth}{ 
\begin{tikzpicture}

\definecolor{darkgray176}{RGB}{176,176,176}
\definecolor{green01270}{RGB}{0,127,0}
\definecolor{orange}{RGB}{255,165,0}
\definecolor{purple}{RGB}{128,0,128}

\begin{axis}[
tick align=outside,
tick pos=left,
x grid style={darkgray176},
xlabel={Wave height (\unit{\centi\meter})},
xmin=37.5, xmax=92.5,
xtick style={color=black},
y grid style={darkgray176},
ylabel={Peak displacement (\unit{\micro\meter})},
ymin=0.15463065, ymax=2.08848235,
ytick style={color=black}
]
\addplot [draw=red, fill=red, mark=*, only marks]
table{%
x  y
40 0.37013
40 0.459889
40 0.386845
40 0.499686
40 0.42529
40 0.69456
40 0.344625
40 0.404007
40 0.355605
40 0.455216
40 0.522292
40 0.324302
40 0.576169
40 0.441205
40 0.419912
40 0.665658
40 0.366802
40 0.277246
40 0.490254
40 0.525749
40 0.485507
40 0.428969
40 0.484472
40 0.345182
40 0.464774
40 0.334539
40 0.738777
40 0.413427
40 0.360502
40 0.432699
40 0.460493
40 0.332275
40 0.615312
40 0.44542
40 0.45249
40 0.423363
40 0.365901
40 0.471588
40 0.438961
40 0.391826
40 0.435145
40 0.313172
40 0.700242
40 0.43192
40 0.417006
40 0.423114
40 0.342526
40 0.357445
40 0.296307
40 0.760415
40 0.59776
40 0.440699
40 0.345988
40 0.589663
40 0.717768
40 0.401447
40 0.364316
40 0.475386
40 0.463006
40 0.325923
40 0.530115
40 0.410658
40 0.409926
40 0.427566
40 0.346679
40 0.334747
40 0.522958
40 0.727989
40 0.636015
40 0.302699
40 0.518013
40 0.571355
40 0.401648
40 0.482438
40 0.356155
40 0.516762
40 0.415818
40 0.342608
40 0.473502
40 0.410961
40 0.364584
40 0.45073
40 0.486977
40 0.521376
40 0.323085
40 0.48488
40 0.445622
40 0.642298
40 0.458941
40 0.560145
40 0.392459
40 0.645205
40 0.437882
40 0.397397
40 0.323586
40 0.467944
40 0.394262
40 0.33764
40 0.399875
40 0.418289
};
\addplot [draw=green01270, fill=green01270, mark=square*, only marks]
table{%
x  y
50 0.541127
50 0.641694
50 0.404563
50 0.451876
50 0.591612
50 0.473338
50 0.575836
50 0.369499
50 0.471779
50 0.874243
50 0.502767
50 0.495899
50 0.710461
50 0.413026
50 0.534794
50 0.470709
50 0.355245
50 0.449331
50 0.436774
50 0.57278
50 0.563447
50 0.536881
50 0.403546
50 0.420751
50 0.657652
50 0.498145
50 0.605891
50 0.599575
50 0.420102
50 0.399509
50 0.430126
50 0.413486
50 0.550814
50 0.520112
50 0.608611
50 0.461679
50 0.317191
50 0.424458
50 0.465633
50 0.620075
50 0.615435
50 0.415328
50 0.513479
50 0.531737
50 0.424241
50 0.387272
50 0.468962
50 0.341161
50 0.462588
50 0.351844
50 0.473571
50 0.444871
50 0.632938
50 0.512838
50 0.417422
50 0.411509
50 0.621516
50 0.440246
50 0.538341
50 0.464753
50 0.414922
50 0.376005
50 0.423801
50 0.423499
50 0.44566
50 0.454235
50 0.499581
50 0.358277
50 0.556862
50 0.450714
50 0.390033
50 0.442821
50 0.485189
50 0.397023
50 0.379196
50 0.343696
50 0.401069
50 0.648884
50 0.423022
50 0.647164
50 0.460997
50 0.399914
50 0.477679
50 0.422507
50 0.539331
50 0.391802
50 0.496496
50 0.788683
50 0.337026
50 0.326104
50 0.602497
50 0.576944
50 0.421978
50 0.417863
50 0.463771
50 0.359432
50 0.502093
50 0.581889
50 0.475426
50 0.56351
};
\addplot [draw=blue, fill=blue, mark=triangle*, only marks]
table{%
x  y
60 0.500457
60 0.394608
60 0.375669
60 0.407049
60 0.399378
60 0.452717
60 0.393431
60 0.443266
60 0.51642
60 0.530877
60 0.463018
60 0.518578
60 0.293289
60 0.469541
60 0.647326
60 0.572527
60 0.505524
60 0.36372
60 0.392608
60 0.438436
60 0.488736
60 0.495723
60 0.320967
60 0.447171
60 0.511194
60 0.446952
60 0.34587
60 0.350789
60 0.318008
60 0.701342
60 0.519659
60 0.39422
60 0.583381
60 0.410529
60 0.430987
60 0.406154
60 0.41111
60 0.568053
60 0.509845
60 0.618062
60 0.41658
60 0.55997
60 0.386923
60 0.484589
60 0.336169
60 0.351175
60 0.421462
60 0.307452
60 0.499996
60 0.393072
60 0.371268
60 0.320422
60 0.353331
60 0.50831
60 0.42018
60 0.433668
60 0.397063
60 0.34854
60 0.425738
60 0.532301
60 0.479248
60 0.629461
60 0.521681
60 0.507745
60 0.336997
60 0.457191
60 0.357161
60 0.412015
60 0.480907
60 0.524681
60 0.405977
60 0.358557
60 0.331477
60 0.47879
60 0.586237
60 0.484263
60 0.448124
60 0.491204
60 0.373367
60 0.294983
60 0.554269
60 0.460742
60 0.464107
60 0.604878
60 0.459757
60 0.384357
60 0.301983
60 0.506475
60 0.427098
60 0.514305
60 0.549722
60 0.477156
60 0.972245
60 0.441647
60 0.51228
60 0.35779
60 0.280087
60 0.339923
60 0.591944
60 0.431439
};
\addplot [draw=purple, fill=purple, mark=x, only marks]
table{%
x  y
70 0.501658
70 0.336709
70 0.797297
70 0.51276
70 0.4486
70 0.433922
70 0.772414
70 0.444198
70 0.566675
70 0.393722
70 0.527786
70 0.441657
70 0.242533
70 0.449971
70 0.321199
70 0.417182
70 0.383148
70 0.601629
70 0.355003
70 0.396984
70 0.340671
70 0.581955
70 0.417893
70 0.490686
70 0.349315
70 0.505674
70 0.389748
70 0.496247
70 0.432408
70 0.403145
70 0.331167
70 0.454378
70 0.353041
70 0.480038
70 0.410564
70 0.353214
70 0.525195
70 0.37123
70 0.510315
70 0.540164
70 0.558597
70 0.359903
70 0.441151
70 0.310478
70 0.374243
70 0.386799
70 0.41271
70 0.358417
70 0.466588
70 0.325339
70 0.38291
70 0.488473
70 0.381056
70 1.36094
70 0.388805
70 0.531249
70 0.478214
70 0.650989
70 0.377102
70 0.308998
70 0.469213
70 0.387874
70 0.389112
70 0.459453
70 0.426659
70 0.498952
70 0.535721
70 0.329002
70 0.420804
70 0.608333
70 0.318651
70 0.724881
70 0.473405
70 0.828516
70 0.5476
70 0.407256
70 0.564185
70 0.317542
70 0.605582
70 0.425051
70 0.37587
70 0.400116
70 0.381296
70 0.46013
70 0.500673
70 0.559642
70 0.350302
70 0.451915
70 0.545133
70 0.303817
70 0.379684
70 0.31906
70 0.391302
70 0.439925
70 0.483919
70 0.419699
70 0.508155
70 0.380898
70 0.315359
70 0.428803
};
\addplot [draw=orange, fill=orange, mark=asterisk, only marks]
table{%
x  y
80 0.574977
80 0.95094
80 0.974833
80 0.667175
80 0.506022
80 0.587332
80 0.518893
80 0.5473
80 0.608814
80 0.630545
80 0.563857
80 0.902064
80 0.572376
80 0.921084
80 0.557379
80 0.654163
80 0.536975
80 0.903276
80 0.430281
80 0.524241
80 0.507012
80 0.684125
80 0.481804
80 0.521872
80 0.536175
80 0.659965
80 0.533711
80 0.588527
80 0.857746
80 0.583895
80 0.528915
80 0.644229
80 0.632772
80 0.782885
80 0.739163
80 0.657771
80 0.555285
80 0.565158
80 0.811462
80 0.620346
80 0.61728
80 0.590882
80 0.562084
80 0.788004
80 0.630962
80 0.54047
80 0.620818
80 0.741179
80 0.566223
80 0.548422
80 0.536778
80 0.660542
80 0.523113
80 0.659099
80 0.554561
80 0.5642
80 0.480762
80 0.512781
80 0.532338
80 0.581521
80 0.500092
80 0.891631
80 0.624696
80 0.626465
80 0.560001
80 0.527632
80 0.478866
80 0.942788
80 0.768104
80 0.571753
80 0.760273
80 0.595677
80 0.603363
80 0.655199
80 0.441722
80 0.530899
80 0.601361
80 0.621336
80 0.678763
80 0.510996
80 0.504447
80 0.653267
80 0.954286
80 0.631424
80 0.872789
80 0.571353
80 0.582889
80 0.648941
80 0.706229
80 0.522486
80 0.587904
80 0.732687
80 0.614445
80 0.674269
80 0.656528
80 0.769458
80 0.853945
80 0.49023
80 0.537545
80 0.638259
};
\addplot [draw=black, fill=black, mark=diamond*, only marks]
table{%
x  y
90 0.935114
90 0.884022
90 0.737487
90 0.8999
90 0.922045
90 0.79517
90 1.06569
90 0.946191
90 0.983649
90 1.51907
90 0.852915
90 1.04761
90 0.90628
90 0.9672
90 1.0922
90 0.864892
90 0.905914
90 1.48984
90 0.828588
90 1.21499
90 1.16192
90 0.954885
90 0.872893
90 0.653097
90 0.9658
90 0.806014
90 1.08339
90 0.700574
90 2.00058
90 0.948938
90 0.727284
90 0.890918
90 0.831139
90 0.971845
90 0.762998
90 1.0272
90 1.23096
90 0.74682
90 1.29558
90 0.938889
90 0.95604
90 1.18875
90 0.900518
90 0.753289
90 0.719921
90 1.03291
90 1.02425
90 0.816324
90 1.00488
90 1.09695
90 1.37732
90 1.22
90 0.918301
90 1.0458
90 1.8864
90 1.29158
90 0.789678
90 0.98027
90 0.8431
90 0.870978
90 0.912153
90 0.806717
90 1.12099
90 0.729767
90 1.04018
90 0.757279
90 1.26106
90 0.814596
90 0.973652
90 1.03831
90 1.14987
90 0.781145
90 1.05491
90 1.05539
90 1.32013
90 0.858211
90 0.796146
90 0.853962
90 1.12798
90 0.904152
90 0.899067
90 1.0221
90 1.0006
90 0.76168
90 1.0604
90 0.882627
90 0.897237
90 0.780176
90 0.939745
90 1.02113
90 1.17305
90 0.811695
90 0.77926
90 1.08014
90 0.990927
90 0.854194
90 0.888517
90 1.35882
90 1.15406
90 0.748724
};
\end{axis}

\end{tikzpicture}}
   \caption{}
   \label{fig:MultipleEvent_peakdisp}
 \end{subfigure}
 \begin{subfigure}{\textwidth}
    \centering
    \includegraphics[width=\columnwidth]{Section4/images/uq/OpenSees/RMSA_legend.pdf}
    \label{fig:legend_PFD} 
  \end{subfigure}
\caption{Two-storey structure's peak displacements obtained using OpenSees vs. (a) W-section weight per length for columns (Column weight), (b) W-section weight per length for beams (Beam weight), (c) W-section weight per length for girder (Girder weight), (d) Yield strength, (e) Youngs' modulus, (f) Wave height.}
\label{fig:Peak displacements_2stories}
\end{figure}

\end{document}